\documentclass[a4paper,UKenglish]{scrreport}
\usepackage{lmodern}
\usepackage[T1]{fontenc}
\usepackage[UKenglish]{babel}
\usepackage{amsmath}
\usepackage{amsfonts}
\usepackage{amssymb}
\usepackage{amsthm}
\usepackage{mathtools}
\usepackage{hyperref}
\usepackage{bookmark}
\usepackage{booktabs}
\usepackage[outer=25mm,inner=40mm,top=30mm,bottom=30mm]{geometry}

\usepackage{verbatim}
\usepackage{faktor}
\usepackage{graphicx}
\usepackage{framed}
\usepackage{tikz}
\usepackage{mathrsfs}
\usepackage{tikz-cd}
\usepackage{enumerate}
\usepackage{faktor}
\usepackage{caption}
\usepackage{dsfont}
\usepackage[utf8]{inputenc}
\usepackage{csquotes}
\usepackage{stmaryrd}
\usepackage[backend=biber,sortlocale=en_GB,url=false,doi=false,isbn=false,giveninits=true,style=alphabetic]{biblatex}
\bibliography{../jabref.bib}
\usepackage[letterspace=150]{microtype}
\usepackage{listings}
\hypersetup{colorlinks=false, linkbordercolor={white}, pdfauthor={Felix Schremmer}, hidelinks,
pdftitle={Master's Thesis -- Weighted projective lines and Hochschild cohomology}}
\title{Weighted projective lines and Hochschild cohomology}
\author{Felix Schremmer}
\date{27 August 2019}
\usetikzlibrary{calc}
\usetikzlibrary{arrows.meta}
\usepackage{fancyhdr}
\DeclareMathOperator{\Aut}{Aut}
\DeclareMathOperator{\Alg}{Alg}
\DeclareMathOperator{\Vect}{Vect}
\DeclareMathOperator{\lcm}{lcm}
\DeclareMathOperator{\Proj}{Proj}
\DeclareMathOperator{\Spec}{Spec}
\DeclareMathOperator{\Frac}{Frac}
\DeclareMathOperator{\SHom}{\mathcal{H}om}
\DeclareMathOperator{\im}{im}
\DeclareMathOperator{\rk}{rk}
\DeclareMathOperator{\incl}{incl}
\DeclareMathOperator{\HH}{HH}
\DeclareMathOperator{\Hom}{Hom}

\let\H\Hoperator

\DeclareMathOperator{\Z}{Z}
\DeclareMathOperator{\Ext}{Ext}
\DeclareMathOperator{\Tor}{Tor}
\DeclareMathOperator{\id}{id}
\DeclareMathOperator{\Mod}{Mod}
\DeclareMathOperator{\catmod}{mod}

\DeclareMathOperator{\Ch}{Ch}

\DeclareMathOperator{\dgMod}{dgMod}
\DeclareMathOperator{\Fun}{Fun}

\DeclareMathOperator{\Ob}{Ob}

\DeclareMathOperator{\End}{End}
\DeclareMathOperator{\rad}{rad}
\DeclareMathOperator{\coker}{coker}
\DeclareMathOperator{\coh}{coh}
\DeclareMathOperator{\Qcoh}{Qcoh}
\DeclareMathOperator{\PGL}{PGL}
\DeclareMathOperator{\RHom}{\mathbf{R}Hom}
\DeclareMathOperator{\fieldchar}{char}
\def\env{{\mathrm e}}
\def\bounded{{\mathrm b}}
\def\op{{\mathrm{op}}}
\def\dg{{\mathrm{dg}}}
\def\D{{\mathrm{D}}}
\def\T{{\mathrm{T}}}
\def\LMod{\text{-}\hspace{-0.3ex}\Mod}
\def\LVect{\text{-}\hspace{-0.3ex}\Vect}
\def\LAlg{\text{-}\hspace{-0.3ex}\Alg}
\def\Lmod{\text{-}\hspace{-0.3ex}\catmod}

\def\Lotimes{\otimes^{\mathbf L}}
\def\rightqed{\pushQED{\qed}\qedhere\popQED}

\def\nextvec#1{
\def\param{#1}
\ifx\param\empty
\end{smallmatrix}\right)
\else
\\#1
\expandafter\nextvec
\fi}
\def\colvec#1{\left(\begin{smallmatrix} #1
\nextvec}
\def\Nextvec#1{
\def\param{#1}
\ifx\param\empty
\end{pmatrix}
\else
\\#1
\expandafter\Nextvec
\fi}
\def\Colvec#1{\begin{pmatrix} #1
\Nextvec}
\def\blank{-}
\def\abs#1{{\left\lvert{}{#1}\right\rvert}}

\newtheorem{lemma}{Lemma}[section]
\newtheorem*{lemma*}{Lemma}
\newtheorem{theorem}[lemma]{Theorem}
\newtheorem{alphatheorem}{Theorem}

\newtheorem*{theorem*}{Theorem}
\newtheorem{corollary}[lemma]{Corollary}
\newtheorem*{corollary*}{Corollary}
\newtheorem{proposition}[lemma]{Proposition}
\newtheorem*{proposition*}{Propositoin}
\theoremstyle{definition}
\newtheorem{definition}[lemma]{Definition}
\newtheorem*{definition*}{Definition}
\theoremstyle{remark}
\newtheorem*{example*}{Example}
\newtheorem{example}[lemma]{Example}
\newtheorem{remark}[lemma]{Remark}
\newtheorem*{remark*}{Remark}

\newtheorem{notation*}{Notation}
\newtheorem{convention}[lemma]{Convention}
\def\citestacks#1{\cite[\href{https://stacks.math.columbia.edu/tag/#1}{#1}]{stacks-project}}
\begin{document}
\maketitle
\addcontentsline{toc}{chapter}{Abstract}
\begin{abstract}
We describe the dimensions of Hochschild (co)homology groups of weighted projective curves over complex numbers. Surprisingly, all but one of those numbers depend only on the genus of the underlying non-weighted curve and the number of exceptional points. Our proof involves revising a classical representation-theoretic argument of Happel together with more recent results of Lenzing and Arinkin, Căldăraru and Hablicsek. We give concrete realizations of a large class of weighted projective lines as quotient stacks. This paper conicides with the author's master's thesis submitted to the University of Bonn in 2019.
\end{abstract}
\chapter{Introduction}
The Hochschild homology and cohomology groups for an algebra over a field (as they were originally defined \cite{Hochschild1945})
give an easily definable yet crucial invariant, that e.g.\ controls the deformation theory of such an algebra.
They have been generalized to algebras over commutative rings, schemes, stacks and differential graded categories.
In this thesis, we will focus in the Hochschild (co)homology of a special family of stacks, namely \emph{weighted projective lines}.
Discussing the various definitions of Hochschild (co)homology in different contexts will constitute the first chapter of this thesis. The main result of this thesis is the following description of the dimensions of Hochschild (co)homology groups of weighted projective curves.
\begin{alphatheorem}[{Cf.\ Theorem~\ref{thm:summaryRootStack}}]
Let $X$ be a smooth irreducible projective curve, $x_1,\dotsc,x_m\in X$ be distinct points and $e_1,\dotsc,e_m\in\mathbb Z_{\geq 2}$.
Let $\mathcal X$ be the iterated root stack obtained from $X$ by attaching the weight $e_i$ to $x_i$ for $i=1,\dotsc,m$.

Denote by $g$ the genus of $X$ and define
\begin{align*}
e := e_1+\cdots + e_m,\quad a:=\begin{cases}3,&g=0\\1,&g=1\\0,&g\geq 2\end{cases},\quad
r:=\min\{a,m\},\quad d := \begin{cases}0,&g=0\\1,&g=1\\3g-3,&g\geq 2\end{cases}.
\end{align*}
The Hochschild (co)homology of $\mathcal X$ in degree $i\in \mathbb Z$ has dimension
\begin{align*}
&\dim\HH_i(\mathcal X)=\begin{cases}2+e-m,&i=0\\
g,&i=\pm 1\\
0,&\text{otherwise}\end{cases}\\
&\dim\HH^i(\mathcal X)=\begin{cases}1,&i=0\\
g+a-r,&i=1\\
d+m-r,&i=2\\
0,&\text{otherwise}.\end{cases}
\end{align*}
\end{alphatheorem}
We remark that $d$ denotes the dimension of the coarse moduli space of curves of genus $g$ and $a$ is the dimension of the automorphism group of a generic curve in it.

A weighted projective line can be understood as a geometric object similar to the projective line $\mathbb P^1$
with a finite number of \emph{exceptional points}, which come with come with a \emph{weight} $> 1$ each. There
is a direct construction of weighted projective lines and their module category in \cite{Geigle1987}.
An alternative approach is the \emph{iterated root stack construction},
which can be used to iteratively attach some weight $>1$ to a point of $\mathbb P^1$
or a weighted projective line. We will discuss both constructions and some properties of weighted projective
lines in the second chapter.

The category of coherent sheaves of a weighted projective line is intimely linked (namely derived equivalent) to the representation theory of a certain
finite-dimensional algebra, a so-called \emph{canonical algebra}. Since the Hochschild (co)homology groups are derived invariants,
one approach to compute the Hochschild (co)homology groups of weighted projective lines is to do this for the corresponding
canonical algebras. This has been done by Happel in \cite{Happel1998}. In the third chapter, we will describe
Happel's computation of the Hochschild (co)homology groups of canonical algebras. We also present a slightly modified version,
in order to shed some light on the question \emph{why does the number of exceptional points alone
determine the Hochschild cohomology?}

The Hochschild (Co)homology groups of smooth projective schemes can be described in terms of their differential forms resp.\ polyvector fields,
this is the famous \emph{Hochschild-Kostant-Rosenberg decomposition} \cite{Hochschild1962}.
This has been generalized to a certain class of stacks, namely orbifold quotients of smooth projective varieties by a finite group action,
by Arinkin, Căldăraru and Hablicsek
\cite{Arinkin2014}.
In the fourth chapter, we recall this result
and use it to calculate the Hochschild (co)homology of quotient stacks of curves by finite group actions.
There is a description due to Lenzing \cite{Lenzing2016}, which weighted projective curves can be realized as the orbifold quotient
of a projective curve by a finite group. He gives no explicit procedure how the curve and the group may be chosen,
but plenty of examples.
Both results combined with Happel's result of the third chapter will entail a full description of the Hochschild (co)homology of weighted projective curves.
Moreover, we discuss how to apply the result of \cite{Arinkin2014} to some of Lenzing's examples. Particularly, we study a class of weighted projective \emph{lines}.
\begin{alphatheorem}[{Cf.\ Corollary~\ref{cor:realizingWPL}}]\label{thm:introRealizing}
Let $x_1,\dotsc,x_m$ be distinct points in $\mathbb P^1$ and $a_1,\dotsc,a_m\in\mathbb Z_{\geq 2}$. Assume for $i=1,\dotsc,m$ that
\begin{align*}
a_i\text{ divides }\lcm(a_1,a_2,\dotsc,a_{i-1},a_{i+1},\dotsc,a_m).
\end{align*}
Let $X$ denote the weighted projective line obtained from $\mathbb P^1$ by attaching the weight $e_i$ to the point $x_i$ for $i=1,\dotsc,m$. Then we can give an explicit curve $C$ together with an action of a finite abelian group $G$ acting on $C$ such that $X \cong [C/G]$.
\end{alphatheorem}
We classify which weighted projective lines have positive resp.\ zero Euler characteristic in Sections~\ref{sec:nonNegEuler}. Note that all cases of Euler characteristic zero can be realized as in Theorem~\ref{thm:introRealizing}.

The subsequent pages of this paper are taken with very minor changes from an unpublished thesis written in 2019 for the sole purpose of supporting the my graduation. To my great pleasure, I was informed six years later that its contents might be of interest beyond this single occasion. With this new purpose of my old thesis and hopefully grown maturity in mathematical writing, I have slightly reworked the introduction in 2025.
\pdfbookmark[chapter]{Contents}{content}
\tableofcontents\thispagestyle{empty}
\chapter{Hochschild (co)homology groups}
Hochschild (co)homology is a homological invariant, which can be associated to a number of objects:
Associative algebras, schemes, stacks, differential graded categories\ldots\
In this chapter, we collect and compare the various definitions.
\section{Towards Hochschild (co)homology of algebras}
For this section, let $k$ denote a commutative ring and $A$ a $k$-algebra. Throughout this thesis,
algebras and rings are understood to be unital and associative, but not necessarily commutative.
If not denoted otherwise, all modules occuring will be left modules.
\begin{definition}
\begin{enumerate}[(a)]
\item The \emph{opposite algebra} of $A$, denoted $A^\op$, has the same $k$-module structure but swapped multiplication:
\begin{align*}
\forall a,b\in A:\qquad a\cdot_{A^\op} b = b\cdot_A a.
\end{align*}
\item The \emph{enveloping algebra} of $A$ is $A^\env := A\otimes_k A^\op$, i.e.\ as $k$-vector space, it is $A\otimes_k A$
with multiplication
\begin{align*}
\forall a,b,c,d\in A:\qquad (a\otimes b)(c\otimes d)=(ac)\otimes (db).
\end{align*}
\end{enumerate}
\end{definition}
\begin{remark}
\begin{itemize}
\item An $A^\env$-module is the same as an $A$-bimodule. For us, the most important example is $A$ itself,
via $(a\otimes b)x = axb$. The structure of this module is closely related to the structure
of the algebra; e.g., submodules of $A$ (as bimodule) correspond to two-sided ideals in $A$ and
direct summands of $A$ as bimodule are direct factors of $A$ as algebra.
\item
Hochschild (co)homology is a homological property of the $A^\env$-module $A$. It is, e.g., important for studying
\emph{deformations} of $A$ (see \cite{Gerstenhaber1988}).
\end{itemize}
\end{remark}
\begin{definition}
Let $M$ be an $A$-bimodule and $n\geq 0$.
\begin{enumerate}[(a)]
\item The \emph{$n$-th Hochschild cohomology group of $A$ with coefficients in $M$} is
\begin{align*}\HH^n(A,M) := \Ext^n_{A^\env}(A,M).\end{align*} If $A=M$, we speak of the \emph{$n$-th Hochschild cohomology group of $A$}, $\HH^n(A)=\HH^n(A,A)$.
\item The \emph{$n$-th Hochschild homology group of $A$ with coefficients in $M$} is
\begin{align*}\HH_n(A,M) := \Tor_n^{A^\env}(A,M).\end{align*} If $A=M$, we speak of the \emph{$n$-th Hochschild homology group of $A$}, $\HH_n(A)=\HH_n(A,A)$.
\end{enumerate}
\end{definition}
\begin{remark}
\begin{itemize}
\item Hochschild (co)homology is a derived invariant of $A$, although it might not appear so from the definition given here.
\item Hochschild (co)homology groups are (at least in the cases interesting to us) not just a sequence of $k$-modules, they
come with various algebraic structures (e.g.\ $\HH^\bullet(A)$ is a graded-commutative $k$-algebra).
The structures won't play an important role for us, so we won't explain them.
\item There is a \enquote{canonical} resolution of $A$ as $A^\env$-module, which is important when working abstractly with Hochschild
(co)homology (it is less important for computations due to its enormous size):
\end{itemize}
\end{remark}
\begin{lemma}[Bar complex]
\begin{enumerate}[(a)]
\item
There exists a complex of $A^\env$-modules $\overline{\mathbb B_\bullet}$
defined as follows:
\begin{align*}
&\overline{\mathbb B_n}=0\text{ if }n<-1\\
&\overline{\mathbb B_n} = A^{\otimes(n+2)}\text{ if }n\geq -1\\
&\overline{\mathbb B_n}\xrightarrow{d_n}\overline{\mathbb B_{n-1}}:
a_0\otimes\cdots\otimes a_{n+1}\mapsto \sum_{i=0}^n (-1)^i a_0\otimes\cdots\otimes a_{i-1}\otimes a_ia_{i+1}\otimes a_{i+2}\otimes\cdots\otimes 
a_{n+1}.
\end{align*}
\item $\overline{\mathbb B_\bullet}$ is, as complex of $A$-modules, contractible with contraction
\begin{align*}
\overline{\mathbb B_{n-1}}\xrightarrow{s_n}\overline{\mathbb B_{n}}: a_0\otimes\cdots\otimes a_n\mapsto 1\otimes a_0\otimes\cdots\otimes a_n.
\end{align*}
\item If $\mathbb B_{\bullet}$ denotes the truncation $\left(\overline{\mathbb B_{\bullet}}\right)_{\geq 0}$,
$\mathbb B_{\bullet}$ together with $d_0 : \mathbb B_{\bullet}\rightarrow A$ is a resolution
of $A$ as $A^\env$-module. $\mathbb B_{\bullet}$ is a free/projective/flat resolution
of $A$ if $A$ is free/projective/flat as $k$-module.\end{enumerate}
\end{lemma}
\begin{proof}
A simple inductive proof checks (a) and (b) together; then (c) follows easily (use that $\mathbb B_n = A^\env\otimes_k A^{\otimes n}$).
\end{proof}
\begin{remark}
When computing $\HH^\bullet(A)$, we thus have to take the cohomology
of the complex $\Hom_{A^\env}(\mathbb B_\bullet,A)$, where $\Hom_{A^\env}(\mathbb B_n,A)\cong \Hom_k(A^{\otimes n},A)$
and the differential becomes
\begin{align*}
&\Hom_k(A^{\otimes n},A)\xrightarrow{\partial^n}\Hom_k(A^{\otimes(n+1)},A),\quad
f\mapsto \partial^n(f)
\\&
\partial^n(f)(a_0\otimes\cdots\otimes a_{n+1})=
a_0f(a_1\otimes\cdots\otimes a_{n+1})
\\&\quad+\sum_{i=1}^n (-1)^if(a_0\otimes\cdots\otimes a_ia_{i+1}\otimes\cdots\otimes a_{n+1})
\\&\quad+(-1)^{n+1} f(a_0\otimes\cdots\otimes a_{n})a_{n+1}.
\end{align*}
In case $A=k$, the Hochschild cohomology is the cohomology of the cochain complex
\begin{align*}
0\rightarrow k\xrightarrow 0k\xrightarrow 1k\xrightarrow 0k\xrightarrow 1k\rightarrow \cdots,
\end{align*}
so $\HH^0(k)=k$ and $\HH^n(k)=0$ for $n\neq 0$.

Moreover, $\HH^0(A) = \ker(\Hom_k(k,A)\xrightarrow{\partial^0} \Hom_k(A,A))$
consists of all $f:k\rightarrow A$ such that $af(1)-f(1)a=0$ for all $a\in A$,
i.e.\ $\HH^0(A) = \Hom_k(k,Z(A))\cong Z(A)$ is the center of $A$.

We will come back to more efficient calculations of Hochschild cohomology of algebras in the third chapter.
\end{remark}
\section{Hochschild (co)homology of differential graded categories}
Hochschild (co)homology was historically first defined for algebras. It is possible to extend the definition to a much broader class of objects, namely
\emph{differential graded categories}. We explain the theory of Hochschild (co)homology of differential graded categories,
mostly following \cite{Keller2003}.

Throughout this section, let $k$ be a commutative ring. Recall the definition of a differential graded category:
\begin{definition}
\begin{itemize}
\item A \emph{$k$-linear category} is a category $\mathcal C$ such that, for any two objects $X,Y\in \mathcal C$,
$\Hom_{\mathcal C}(X,Y)$ is a $k$-module and
the composition maps
\begin{align*}
\Hom_{\mathcal C}(X,Y)\times \Hom_{\mathcal C}(Y,Z)\rightarrow \Hom_{\mathcal C}(X,Z)
\end{align*}
are $k$-bilinear.
\item If $\mathcal C$ is a $k$-linear category
and moreover for each $X,Y\in \mathcal C$, $\Hom_{\mathcal C}(X,Y)$
is a cochain complex over $k$ such that the composition maps
\begin{align*}
\Hom_{\mathcal C}(X,Y)\otimes \Hom_{\mathcal C}(Y,Z)\rightarrow \Hom_{\mathcal C}(X,Z)
\end{align*}
are morphisms of cochain complexes (for the usual tensor product of cochain complexes)
and such that $\id_X\in \Z^0(\Hom_{\mathcal C}(X,X))$ for all objects $X\in\mathcal C$,
we say $\mathcal C$ is a differential graded category.
\end{itemize}
\end{definition}
\begin{example}\label{example_chaincomplex}The following are examples of dg categories:
\begin{itemize}
\item If $A$ is a $k$-algebra (i.e.\ a $k$-linear category with one object), it defines
a dg category with one object by setting the differential to be zero.
A dg category with one object is called \emph{dg algebra}.
\item Recall that for a $k$-linear category $\mathcal C$,
we may define the dg category of cochain complexes $\Ch^\dg(\mathcal C)$,
whose objects are cochain complexes in $\mathcal C$.
For cochain complexes $C^\bullet,D^\bullet$ over $\mathcal C$, we set
$\Hom_{\mathcal C}(C^\bullet,D^\bullet)$ to be the cochain complex where
$\Hom_{\mathcal C}(C^\bullet,D^\bullet)^n$ consists of all families $f^\bullet$
of morphisms $(f^m : C^m\rightarrow D^{n+m})_{m\in\mathbb Z}$.
The differential of $\Hom_{\mathcal C}(C^\bullet,D^\bullet)$ sends $f^\bullet \in \Hom_{\mathcal C}(C^\bullet,D^\bullet)^n$ to
\begin{align*}
d_D\circ f + (-1)^{n+1} f\circ d_C.
\end{align*}
\end{itemize}
\end{example}
The last example is important for defining the module category of a dg category:
\begin{definition}
For a dg category $\mathcal C$, we define its \emph{module category} $\dgMod_{\mathcal C}$ to be
the category of dg functors
$\Fun^\dg(\mathcal C, \Ch^\dg(k)^\op)$.
It is a dg category where the morphisms of degree $n$ are those natural transformations which are pointwise of degree $n$, and sums and differentials
are computed pointwise.
\end{definition}
\begin{example}
\begin{enumerate}[(a)]
\item
If $A$ is a $k$-algebra, interpreted as dg category with one object $\ast$ and
$\Hom(\ast,\ast)=A$, concentrated in degree zero, a functor $F:A\rightarrow \Ch^\dg(k)^\op$
consists of a choice of object $C^\bullet:=F\ast$ in $\Ch^\dg(k)$ and a $k$-algebra homomorphism
$A\rightarrow Z^0(\End_{\Ch^\dg}(C^\bullet))$. i.e.\ to each $a\in A$
and $n\in\mathbb Z$, we associate a $k$-linear endomorphism $F(a) : C^n\rightarrow C^n$
that is compatible with the differentials of $C^\bullet$.
This precisely describes $C^\bullet$ as a chain complex of $A$-modules.
We see that $\dgMod_A \cong \Ch^\dg(A\LMod)$ as dg categories.
\item
If $c\in\mathcal C$ is an object, the functor
\begin{align*}
\mathcal C(\blank,c) : \mathcal C\rightarrow\Ch^\dg(k)^\op
\end{align*}
is a dg module over $\mathcal C$. The induced functor $\mathcal C\rightarrow \dgMod_{\mathcal C}$
is called the \emph{Yoneda embedding}. The dg modules in its image are called the \emph{representable} dg modules.
\end{enumerate}
\end{example}
In order to define the Hochschild homology and cohomology of a dg category, we introduce its derived category:
\begin{definition}
Let $\mathcal C$ be a dg category.
\begin{itemize}
\item The \emph{underlying category} $\Z^0(\mathcal C)$ is the category with the same objects as $\mathcal C$, where for $X,Y\in\Ob(\mathcal C)$,
\begin{align*}
\Hom_{\Z^0(\mathcal C)}(X,Y) = \Z^0\left(\Hom_{\mathcal C}(X,Y)\right).
\end{align*}
Composition is defined by restricting the composition of $\mathcal C$. $\Z^0(\mathcal C)$ is $k$-linear. One similarly defines $\H^0(\mathcal C)$.
\item Consider the category $\Z^0(\dgMod_{\mathcal C})$. A morphism
$f\in \Hom_{\Z^0(\dgMod_{\mathcal C})}(M,N)$ is called \emph{quasi-isomorphism}
if, for all $x\in \mathcal C$, the induced morphism
\begin{align*}
M(x)\xrightarrow{f(x)} N(x)
\end{align*}
is a quasi-isomorphism in $\Ch^{\dg}(k)$.
The \emph{derived category} of $\mathcal C$, denoted $\D(\mathcal C)$, is the localization of $\Z^0(\dgMod_{\mathcal C})$ at the quasi-isomorphisms.
\end{itemize}
\end{definition}
\begin{example}
By Example \ref{example_chaincomplex}, $\Z^0(\Ch^\dg(k)) \cong \Ch(k)$.

Let $A$ be a $k$-algebra, interpreted as dg category.
Since $\dgMod_A\cong\Ch^\dg(A\LMod)$, we see $\Z^0(\dgMod_A)\cong \Z^0(\Ch^\dg(A\LMod))\cong \Ch(A\LMod)$
and hence the derived category of $A$ as dg category and the derived category of $A$ as algebra coincide.
\end{example}
\begin{remark}The derived category of a dg category is triangulated (shifts and cones of modules
are defined by taking shifts and cones in $\Ch(k)$). Similarly, the homotopy category $\H^0(\dgMod_{\mathcal C})$
is triangulated.
\end{remark}
\begin{remark}
For a $k$-algebra $A$, \begin{align*}
&\HH^\bullet(A) \cong \Ext_{A^\env}^\bullet(A,A) \cong \Hom_{\D(A^\env)}(A,A[\bullet]) \cong \H_\bullet(\RHom_{A^\env}(A,A)).
\\&\HH_\bullet(A)\cong \Tor^{A^\env}_\bullet(A,A)\cong \H_\bullet(A\Lotimes_{A^\env} A).
\end{align*}
We want to generalise the latter expressions to dg categories in order to define their Hochschild (co)homology.
\end{remark}
\begin{definition}
Let $\mathcal C$ and $\mathcal D$ be dg categories over $k$.
\begin{itemize}
\item We define the dg category $\mathcal C^\op$ by $\Hom_{\mathcal C^\op}(X,Y) := \Hom_{\mathcal C}(Y,X)$ as cochain complexes over $k$,
composition defined as in $\mathcal C$.
\item We define the dg category $\mathcal C\otimes \mathcal D:= \mathcal C\otimes_k\mathcal D$ by
$\Ob(\mathcal C\otimes \mathcal D) = \Ob(\mathcal C)\times \Ob(\mathcal D)$
and for $X,X'\in \Ob(\mathcal C), Y,Y'\in \Ob(\mathcal D)$,
we define $\Hom_{\mathcal C\otimes \mathcal D}( (X,Y), (X',Y'))$ to be the cochain complex tensor product
$\Hom_{\mathcal C}(X,X')\otimes_k \Hom_{\mathcal D}(Y,Y')$.
Composition of morphisms and identity morphisms are calculated for each tensor factor individually.
\item A $\mathcal C$-$\mathcal D$ dg bimodule is, by definition, a dg module over $\mathcal C\otimes_k \mathcal D^\op$.
\end{itemize}
\end{definition}
\begin{remark}
Given a $\mathcal C$-$\mathcal D$ dg bimodule $M$ and a $\mathcal D$-$\mathcal E$ dg bimodule $N$,
we can define the tensor product $M\otimes N$ similar to the construction for algebras, cf.\ \cite[section~6.1]{Keller1994}.

The tensor product and the Hom-functor can be derived as for algebras, yielding
\begin{align*}
&\text{-}\Lotimes \text{-} : \D(\mathcal C\otimes \mathcal D^\op)\times \D(\mathcal D\otimes \mathcal E^\op)\rightarrow \D(\mathcal C\otimes \mathcal E^\op).
\\&\RHom:\D(\mathcal C\otimes \mathcal E^\op)^\op\times \D(\mathcal D\otimes \mathcal E^\op)\rightarrow \D(C^\op\otimes \mathcal D).
\end{align*}
\end{remark}
\begin{definition}[Hochschild (co)homology of dg categories]
Let $\mathcal C$ be a dg category over the commutative ring $k$.
\begin{itemize}
\item The \emph{diagonal bimodule} is the $\mathcal C\text{-}\mathcal C\text{-}$bimodule, denoted (by abuse of notation) by $\mathcal C$, is defined by
\begin{align*}
\mathcal C:\mathcal C\otimes \mathcal C^\op\rightarrow \Ch^\dg(k)^\op,
(c,d)\mapsto \mathcal C(c,d) = \Hom_{\mathcal C}(c,d).
\end{align*}
\item The \emph{Hochschild cohomology} of $\mathcal C$ in degree $n\geq 0$, $\HH^n(\mathcal C)\in k\LMod$, is defined to be the $n$-th cohomology group of
\begin{align*}
\RHom_{\mathcal C\otimes \mathcal C^\op}(\mathcal C,\mathcal C)\in \D(k).
\end{align*}
\item The \emph{Hochschild homology} of $\mathcal C$ in degree $n\geq 0$, $\HH_n(\mathcal C)\in k\LMod$, is defined to be the $n$-th homology group of
\begin{align*}
\mathcal C\Lotimes_{\mathcal C\otimes \mathcal C^\op}\mathcal C\in \D(k).
\end{align*}
\end{itemize}
\end{definition}
\begin{example}
If $A$ is a $k$-algebra, we can interpret it as dg category over $k$ and get the same notion of Hochschild (co)homology.
\end{example}
\begin{remark}\label{rem:kellerHochschild}
Hochschild (co)homology of differential graded categories is a derived invariant:
Suppose $\mathcal C$ and $\mathcal D$ are dg categories and $M$ is a $\mathcal C\text{-}\mathcal D\text{-}$bimodule
such that $M\Lotimes\blank : \D(\mathcal D)\rightarrow \D(\mathcal C)$ is an equivalence of categories.
Then $M$ induces an isomorphism of Hochschild (co)homology groups of $\mathcal C$ and $\mathcal D$
(follows from \cite{Keller2003}).

Such equivalences e.g.\ arise from tilting objects:
\end{remark}
\begin{definition}
Let $\mathcal T$ be a $k$-linear triangulated category allowing arbitrary direct sums.
\begin{enumerate}[(a)]
\item An object $X\in\mathcal T$ is called \emph{compact} if the functor
\begin{align*}
\mathcal T(X,\blank) : \mathcal T\rightarrow\mathcal T
\end{align*}
commutes with arbitrary direct sums.
\item A class of objects $C\subseteq \Ob(\mathcal T)$ is called \emph{generating}
if the smallest full triangulated subcategory of $\mathcal T$ closed under arbitrary direct sums, direct summands
and isomorphisms containing $C$ is equal to $\mathcal T$.
\item An object $X\in\mathcal T$ is called \emph{tilting} if $X$ is compact, $\{X\}$ is generating and
$\mathcal T(X,X[n])=0$ for all $n\in\mathbb Z\setminus\{0\}$.
\end{enumerate}
\end{definition}
Tilting objects give derived equivalences:
\begin{proposition}\label{prop:tiltingDerivedEquivalence}
Let $\mathcal C$ be a dg category over the field $k$ and let $X$ be a dg module over $\mathcal C$:
\begin{enumerate}[(a)]
\item Suppose that $X$ is compact in $\D(\mathcal C)$
and that $\{X\}$ is generating in $\D(\mathcal C)$.
Let $A$ be the dg algebra $\RHom(X,X)$.

There exists a $A\text{-}\mathcal C\text{-}$bimodule $M$ such that
$\D(\mathcal C)\xrightarrow{M\Lotimes\blank}\D(A)$ is an equivalence of categories.
\item Suppose that $X$ is tilting in $\D(\mathcal C)$
and let $A := \End_{\D(\mathcal C)}(X)$. There exists an $A\text{-}\,\mathcal C\text{-}$bimodule
$M$ such that $\D(\mathcal C)\xrightarrow{M\Lotimes \blank}\D(A)$ is an equivalence of categories.
\end{enumerate}
\end{proposition}
\begin{proof}
(a) follows from \cite[Example~6.6]{Keller1994}. Using (a), (b) follows since
$\RHom(X,X)$ is concentrated in degree $0$, where it is $\End_{\D(\mathcal C)}(X)$.
\end{proof}
Compact objects can moreover be described as follows:
\begin{proposition}\label{prop:compactObjDG}
Let $\mathcal C$ be a dg category. The full subcategory of compact objects in $\D(\mathcal C)$
equals the smallest full triangulated subcategory of $\D(\mathcal C)$ closed under direct summands
and isomorphisms which contains all representable dg modules over $\mathcal C$.
\end{proposition}
\begin{proof}
This follows from \cite[Section~5.3]{Keller1994}.
\end{proof}
We will apply these results when dealing with dg enhancements.
\section{Hochschild (co)homology of abelian categories}
Throughout this section, let $k$ be a commutative ring an $\mathcal A$
be an abelian $k$-linear category. We would like to define the notion
of Hochschild (co)homology of $\mathcal A$. We do this by replacing $\mathcal A$
by a suitable dg category.
\begin{definition}
\begin{enumerate}[(a)]
\item A dg category $\mathcal C$ over $k$ is called \emph{pretriangulated}
if the image of $\H^0(\mathcal C)$ in $\H^0(\dgMod_{\mathcal C})$ (under the Yoneda embedding)
is a triangulated subcategory.
\item A \emph{dg enhancement} of $\D^\bounded(\mathcal A)$ consists of a pretriangulated dg category
$\mathcal C$ together with an equivalence of triangulated categories $\D^\bounded(\mathcal A)\cong \H^0(\mathcal C)$.
\item A cochain complex $I^\bullet$ over $\mathcal A$ is called \emph{homotopy injective}
if $\Hom_{K(\mathcal A)}(C^\bullet,I^\bullet)=0$ for every acyclic complex $C^\bullet$ over $\mathcal A$.
\end{enumerate}
\end{definition}
\begin{proposition}\label{prop:hinjEnhancement}
Suppose that $\mathcal A$ is Grothendieck abelian.
Ley $\mathcal C\subseteq \Ch^\dg(\mathcal A)$ denote the full subcategory consisting of those homotopy injective complexes
which have bounded cohomology. Then $\mathcal C$ gives rise to a dg enhancement of $\H^0(\mathcal C)\cong \D^\bounded(\mathcal A)$.
\end{proposition}
\begin{proof}
This follows from the existence of homotopy injective resolutions, cf.\ \citestacks{079P}
using a proof similar to \cite[Proposition~3.5.43]{Zimmermann2014}.
\end{proof}
We use this dg enhancement to define Hochschild (co)homology of the abelian category $\mathcal A$.
\begin{definition}
Let $n\in\mathbb Z$ and $\mathcal C$ be the dg category from Proposition~\ref{prop:hinjEnhancement}.
Suppose that $\mathcal A$ is Grothendieck abelian.
We define the $n$-th Hochschild (co)homology group of $\mathcal A$ to be the $n$-th Hochschild (co)homology
group of $\mathcal C$.
\end{definition}
This definition satisfies the following derived invariance property:
\begin{proposition}\label{prop:abelianTilting}
Suppose that $k$ is a field and $\mathcal A$ Grothendieck abelian.
Let $T\in \D^\bounded(\mathcal A)$ be an object
such that
\begin{itemize}
\item $\D^\bounded(\mathcal A)$ is the smallest full triangulated subcategory of $\D^\bounded(\mathcal A)$ closed under direct summands and isomorphisms
that contains $T$.
\item $\Hom_{\D^\bounded(\mathcal A)}(T,T[n])=0$ for $n\in\mathbb Z\setminus\{0\}$.
\end{itemize}
$\D^\bounded(\mathcal A)$ and $\D^\bounded(A)$ are equivalent as triangulated categories.
The Hochschild (co)homology of the abelian category $\mathcal A$ is isomorphic to the Hochschild (co)homology of the algebra $\End_{\D^\bounded(\mathcal A)}(T)$.
\end{proposition}
\begin{proof}
Observe that the image of $T$ in $\D(\mathcal C)$ is tilting with endomorphism algebra $\End_{\D^\bounded(\mathcal A)}(T)$,
apply Proposition~\ref{prop:tiltingDerivedEquivalence} and Remark~\ref{rem:kellerHochschild}.
\end{proof}
\begin{example}
If $A$ is a $k$-algebra ($k$ again a field) and $\mathcal A=A\LMod$,
$A$ is a tilting object in $\mathcal A$, so that $A\LMod$ and $A^\op$
have the same Hochschild (co)homology. It is simple to see (from the definition) that the Hochschild (co)homology groups
of $A$ and $A^\op$ are isomorphic. Hence the Hochschild (co)homology groups of $A$ and $A\LMod$ are isomorphic.
\end{example}
\begin{example}
Let $k$ be a field and consider the path algebras for the quiver
\begin{align*}
Q : 1\rightarrow 2\rightarrow 3.
\end{align*}
Let, for $v\in Q_0 = \{1,2,3\}$, $P(v) := kQ \cdot v$ denote the projective module for the vertex $v$,
and $S(v) = P(v)/\rad P(v)$ the simple module. The projective modules have composition series
given as follows:
\begin{align*}
\begin{tikzcd}[ampersand replacement=\&]
P(1)\ar[d,no head,"S(1)"]\& P(2)\ar[d,no head,"S(2)"]\& P(3)\ar[d, no head, "S(3)"]\\
0\&
P(1)\ar[d,no head,"S(1)"]\& P(2)\ar[d,no head,"S(2)"]\\
\&0\&
P(1)\ar[d,no head,"S(1)"]
\\\&\&0
\end{tikzcd}
\end{align*}
Let $T := S(2) \oplus P(2) \oplus P(3)$.
\begin{itemize}
\item
We claim that $\Ext^1_{kQ}(T,T)=0$:
We only have to show that $\Ext^1_{kQ}(S(2),M)=0$ for $M\in \{S(2),P(2),P(3)\}$.
Note that a projective resolution of $S(2)$ is given by \begin{align*}0\rightarrow P(1)\rightarrow P(2)\rightarrow S(2)\rightarrow 0.\end{align*}
Hence for any $kQ$-module $M$,
\begin{align*}
\Ext^1_{kQ}(S(2),M) \cong \coker\left(\Hom_{kQ}(P(2),M)\rightarrow \Hom_{kQ}(P(1),M)\right).
\end{align*}
Now for each $M\in \{S(2),P(2),P(3)\}$, one quickly checks that every morphism $P(1)\rightarrow M$
factors over $P(1)\rightarrow P(2)$, showing that $\Ext^1_{kQ}(S(2),M)=0$.
\item We claim that the image of $T$ in $\D^\bounded(kQ)$ is tilting:
By the previous point, we have $\Hom_{\D^\bounded(kQ)}(T,T[n])=\Ext^n_A(T,T)=0$ for $0\neq n\in\mathbb Z$.
Moreover, the projective modules $P(1)$, $P(2)$ and $P(3)$ are direct summands of $A$,
hence all of them are compact (using Proposition~\ref{prop:compactObjDG}). As the short exact sequence
\begin{align*}
0\rightarrow P(1)\rightarrow P(2)\rightarrow S(2)\rightarrow 0
\end{align*}
gives rise to a distinguished triangle in $\D^\bounded(kQ)$, also $S(2)$ is compact. We conclude that $T$ is compact.

We still have to show that $T$ generates $\D^\bounded(kQ)$. Let $\mathcal C\subseteq \D^\bounded(kQ)$
be the smallest full triangulated subcategory that contains $T$ and is closed under isomorphisms and direct summands.
It suffices to show that $A\in \mathcal C$. Certainly $S(2), P(2)$ and $P(3)$ are in $\mathcal C$.
By the distinguished triangle $P(1)\rightarrow P(2)\rightarrow S(2)\rightarrow P(1)[1]$ considered earlier,
also $P(1)\in \mathcal C$. Thus $kQ=P(1)\oplus P(2)\oplus P(3)$ is in $\mathcal C$.
\item We are interested in the endomorphism ring of $T$ (in $kQ\Lmod$ or $\D^\bounded(kQ)$, they are the same):
We calculate
\begin{align*}
\begin{array}{ccc}
\Hom_{kQ}(S(2),S(2)) = k \id&\Hom_{kQ}(S(2),P(2)) = 0&\Hom_{kQ}(S(2),P(3)) = 0\\[1em]
\dim \Hom_{kQ}(P(2),S(2))=1&\Hom_{kQ}(P(2),P(2)) = k \id&\dim\Hom_{kQ}(P(2),P(3))=1\\[1em]
\Hom_{kQ}(P(3),S(2))=0&\Hom_{kQ}(P(3),P(2)) = 0&\Hom_{kQ}(P(3),P(3))=k \id\\
\end{array}
\end{align*}

From this description, it is easy to see that the endomorphism algebra of $T$ is isomorphic to
the path algebra over the quiver
\begin{align*}
\bullet\leftarrow\bullet\rightarrow \bullet.
\end{align*}
We see that $kQ$ is derived equivalent to the path algebra $k(\bullet\rightarrow\bullet\leftarrow\bullet)$,
but rather obviously not isomorphic and thus not Morita equivalent (as both are basic algebras).
Moreover, we conclude that the Hochschild (co)homology of these algebras will coincide, although we do not yet have
the tools to compute them for any of the two algebras.
\end{itemize}
\end{example}
\section{Hochschild (co)homology of schemes}
Let $k$ be a commutative ring. An affine $k$-scheme is the same as a commutative $k$-algebra,
hence comes with a notion of Hochschild (co)homology groups.
We can expand this definition to arbitrary $k$-schemes $X$, where we assume for simplicity that
$k$ is a field and $X$ noetherian.
\begin{definition}
Let $X$ be a noetherian $k$-scheme.
Let $\coh X$ denote the category of coherent sheaves on $X$.
We define the Hochschild cohomology and Hochschild homology groups of $X$ to be
the Hochschild (co)homology groups of this abelian category.
\end{definition}
\begin{remark}
This definition has the advantage that it allows us to use the results on Hochschild (co)homology of abelian categories
from the previous section. However, it may seem unnatural to define the Hochschild (co)homology of a scheme $X$
by referring to the entire category $\coh X$, just like we defined the Hochschild (co)homology of an algebra without referring
to referring to the entire module category.

There are other ways of defining Hochschild (co)homology of schemes, which do not refer to such huge objects.
e.g.\ \cite{Kuznetsov2009} uses
\begin{align*}
&\HH^\bullet(X) = \Ext^\bullet_{X\times X}(\Delta_\ast \mathcal O_X,\Delta_\ast \mathcal O_X),
\\&\HH_\bullet(X) = \Ext^\bullet_{X\times X}(\Delta_\ast \mathcal O_X,\Delta_\ast \omega_X[\dim X]),
\end{align*}
where $\Delta : X\rightarrow X\times X$ is the embedding of the diagonal.
This has its advantages (e.g.\ it is easier to prove finiteness results from this definition than from the one used here).

We chose to define the Hochschild (co)homology using the detour of dg categories and abelian categories
in order to show how tilting invariance is achieved.
\end{remark}
\begin{example}
If $X=\Spec k$ for a field $k$, $\coh X\cong k\Lmod$, so the object $k$ is tilting in $\D^\bounded(\coh X) =: \D^\bounded(X)$.
We conclude $\HH^n(X)\cong \HH^n(k)$ and $\HH_n(X)\cong \HH_n(k)$ are isomorphic to $k$ for $n=0$ and $0$ for $n\neq 0$.
\end{example}
\begin{definition}[Exceptional sequence, {\cite[Definition~1.57]{Huybrechts2006}}]\label{def:exceptional}
Let $\mathcal D$ be a $k$-linear triangulated category.
\begin{enumerate}[(a)]
\item An object $E\in\mathcal D$ is called \emph{exceptional}
if 
\begin{align*}
\Hom_{\mathcal D}(E,E[\ell])\cong\begin{cases}k,&\ell=0\\0,&\ell\neq 0\end{cases}
\end{align*}
for all $\ell\in\mathbb Z$.
\item A sequence $E_1,\dotsc,E_n$ of objects in $\mathcal D$ is called \emph{exceptional} if 
each $E_i$ is exceptional and $\Hom_{\mathcal D}(E_i,E_j[\ell])=0$ if $i<j$ and $\ell\in\mathbb Z$.
\item An exceptional sequence $E_1,\dotsc,E_n$ is called \emph{full} if the smallest full triangulated subcategory
of $\mathcal D$ closed under direct summands and isomorphisms containing $E_1,\dotsc,E_n$ is $\mathcal D$.
\end{enumerate}
\end{definition}
\begin{proposition}\label{prop:exceptionalTilting}
Let $E_1,\dotsc,E_n$ be a full exceptional sequence in $\D^\bounded(X)$.
Then $X$ is derived equivalent to the algebra
\begin{align*}
\End_{\D^\bounded(X)}(E_1\oplus\cdots\oplus E_n)\cong \begin{pmatrix}k&0&\hdots&0\\
\Hom_{\D^\bounded(X)}(E_2,E_1)&k&\ddots&\vdots\\
\vdots&\ddots&\ddots&0\\
\Hom_{\D^\bounded(X)}(E_n,E_1)&\hdots&\Hom_{\D^\bounded(X)}(E_n,E_{n-1})&k
\end{pmatrix}.
\end{align*}
The Hochschild (co)homology groups of $X$ and this algebra are isomorphic.
\end{proposition}
\begin{proof}
Observe that $E_1\oplus\cdots\oplus E_n$ satisfies the conditions of Proposition~\ref{prop:abelianTilting}.
For the last isomorphism, use that $\End_{\D^\bounded(X)} \cong \left(\Hom_{\D^\bounded(X)}(E_i,E_j)\right)_{i,j=1}^n$
and simplify using the definition of an exceptional sequence.
\end{proof}
\begin{example}\label{ExampleDerivedEquivalenceProjectiveSpace}
Let $k$ be algebraically closed and $X=\mathbb P^n$.
Beilinson proved in \cite{Beilinson1978} that
\begin{align*}
\mathcal O_{\mathbb P^n}, \mathcal O_{\mathbb P^n}(1),\dotsc,\mathcal O_{\mathbb P^n}(n)
\end{align*}
is a full exceptional sequence in $\D^\bounded(\mathbb P^n)$.

Since the homomorphisms $\mathcal O_{\mathbb P^n}(r)\rightarrow \mathcal O_{\mathbb P^n}(s)$ are in one-to-one correspondence
with homogeneous polynomials in $k[X_0,\dotsc,X_n]$ of degree $s-r$, one easily calculates
that the algebra $A$ of Proposition~\ref{prop:exceptionalTilting} is isomorphic to
the bound quiver algebra defined by
\begin{align*}
\begin{tikzpicture}
\draw (0,0) node {$0$};
\draw (2,0) node {$1$};
\draw (4,0) node {$2$};
\draw (6,0) node {$\hdots$};
\draw (8,0) node {$n-1$};
\draw (10,0) node {$n$};
\draw[->] (0.4,1.2) -- (1.6,1.2) node[pos=.5,above] {$X_0$};
\draw[->] (0.4,.5) -- (1.6,.5) node[pos=.5,above] {$X_1$};
\draw (1,.3) node {$\vdots$};
\draw[->] (0.4,-.5) -- (1.6,-.5) node[pos=.5,above] {$X_{n-1}$};
\draw[->] (0.4,-1.2) -- (1.6,-1.2) node[pos=.5,above] {$X_{n}$};
\draw[->] (2.4,1.2) -- (3.6,1.2) node[pos=.5,above] {$X_0$};
\draw[->] (2.4,.5) -- (3.6,.5) node[pos=.5,above] {$X_1$};
\draw (3,.3) node {$\vdots$};
\draw[->] (2.4,-.5) -- (3.6,-.5) node[pos=.5,above] {$X_{n-1}$};
\draw[->] (2.4,-1.2) -- (3.6,-1.2) node[pos=.5,above] {$X_{n}$};
\draw[->] (4.4,1.2) -- (5.6,1.2) node[pos=.5,above] {$X_0$};
\draw[->] (4.4,.5) -- (5.6,.5) node[pos=.5,above] {$X_1$};
\draw (5,.3) node {$\vdots$};
\draw[->] (4.4,-.5) -- (5.6,-.5) node[pos=.5,above] {$X_{n-1}$};
\draw[->] (4.4,-1.2) -- (5.6,-1.2) node[pos=.5,above] {$X_{n}$};
\draw[->] (6.4,1.2) -- (7.6,1.2) node[pos=.5,above] {$X_0$};
\draw[->] (6.4,.5) -- (7.6,.5) node[pos=.5,above] {$X_1$};
\draw (7,.3) node {$\vdots$};
\draw[->] (6.4,-.5) -- (7.6,-.5) node[pos=.5,above] {$X_{n-1}$};
\draw[->] (6.4,-1.2) -- (7.6,-1.2) node[pos=.5,above] {$X_{n}$};
\draw[->] (8.4,1.2) -- (9.6,1.2) node[pos=.5,above] {$X_0$};
\draw[->] (8.4,.5) -- (9.6,.5) node[pos=.5,above] {$X_1$};
\draw (9,.3) node {$\vdots$};
\draw[->] (8.4,-.5) -- (9.6,-.5) node[pos=.5,above] {$X_{n-1}$};
\draw[->] (8.4,-1.2) -- (9.6,-1.2) node[pos=.5,above] {$X_{n}$};
\end{tikzpicture}
\end{align*}
(so it has $(n+1)$ vertices and $n(n+1)$ arrows),
modulo the relations that
\begin{align*}
\left(a \xrightarrow{X_i} (a+1) \xrightarrow{X_j} (a+2)\right) = \left(a\xrightarrow{X_j} (a+1) \xrightarrow{X_i} (a+2)\right)
\end{align*}
for all vertices $a\in \{0,\dotsc,n-2\}$ and $i,j\in \{0,\dotsc,n\}$.

Hence $\D^\bounded(X)\cong \D^\bounded(A^\op)\cong \D^\bounded(A)$.

For $n=1$, we see that the category of coherent sheaves on $\mathbb P^1$ is
derived equivalent to the \emph{Kronecker algebra} $A = k\left(\begin{tikzcd}[column sep=1em]\bullet\ar[r,shift left=.3em]\ar[r,shift right=.3em]&\bullet\end{tikzcd}\right)$.
\end{example}

To finish this chapter, we present one of the gems of the theory of Hochschild (co)homology
for schemes: The \emph{Hochschild-Kostant-Rosenberg decomposition}.
\begin{theorem}[Hochschild-Kostant-Rosenberg, \cite{Hochschild1962}]\label{thm:HKR}
Let $k$ be an algebraically closed field of characteristic zero, and $X$ a smooth projective $k$-scheme.
Denote by $\T_X$ the tangent bundle of $X$ and by $\Omega^q$ the $q$-th exterior power of the cotangent bundle. Then
\begin{align*}
&\HH^i(X) \cong \bigoplus_{p+q=i} \H^p(X,\bigwedge^q \T_X)
\\&\HH_i(X) \cong \bigoplus_{q-p=i} \H^p(X,\Omega^q_X).
\end{align*}
\end{theorem}
We give two applications of this theorem:
\begin{corollary}\label{cor:hkrVanishingAndSymmetry}
If $k,X$ are as above and $\dim X = n$, we have
\begin{itemize}
\item
$\HH^i(X)=0$ if $i\notin [0,2n]$.
\item $\HH_i(X) =\HH_{-i}(X)$, and if $i\notin [-n,n]$, this is moreover $0$.
\end{itemize}
\end{corollary}
\begin{proof}
Recall that $\Omega^q=0=\bigwedge^q\T_X$ for $q\notin [0,n]$, and $\H^p(X,\mathcal F)=0$ for $p\notin[0,n]$ and $\mathcal F$ quasi-coherent, this yields the vanishing assertions.

Let $X = \bigcup_{i=1}^m X_i$ be the decomposition of $X$ into irreducible components (since $X$ is smooth, these are also the connected components).
Let $n_i := \dim X_i$. By Serre duality, $\H^p(X_i,\Omega^q_{X_i}) = \H^{n_i-p}(X_i,\Omega^{n_i-q}_{X_i})$
(using the perfect pairing $\Omega^{n_i-q}_{X_i}\otimes \Omega^q_{X_i}\rightarrow \omega_{X_i}, (u\otimes v)\mapsto u\wedge v$).
Thus
\begin{align*}
\HH_\bullet(X) \cong&\bigoplus_{q-p=\bullet}\bigoplus_{i=1}^m \H^p(X_i,\Omega^q_{X_i})
\cong\bigoplus_{i=1}^m \bigoplus_{q-p=\bullet} \H^{n_i-p}(X_i,\Omega^{n_i-q}_{X_i})
\\\cong&\bigoplus_{i=1}^m \bigoplus_{p'-q'=\bullet} \H^{p'}(X_i,\Omega^{q'}_{X_i})
\cong \HH_{-\bullet}(X).
\end{align*}
This finishes the proof.
\end{proof}
\begin{example}
Let $k$ be an algebraically closed field of characteristic zero and
$X=\mathbb P^1$. We would like to descibe $\HH^\bullet(X)$ and $\HH_\bullet(X)$.
Recall that if $p\neq 0,1$, $\H^p(X,\mathcal F)=0$ for any quasi-coherent sheaf $\mathcal F$. Thus
\begin{align*}
&\HH^0(X) = \H^0(X,\mathcal O_X)\cong k\oplus 0=k
\\&\HH^1(X)\cong \H^1(X,\mathcal O_X)\oplus \H^0(X,\T_X)\cong 0\oplus k^3 \cong k^3
\\&\HH^2(X)\cong \H^1(X,\T_X) =0.\\
&\HH_{-1}(X)\cong \HH_1(X) \cong \H^0(X,\Omega^1_X)\cong 0
\\&\HH_0(X)\cong \H^0(X,\mathcal O_X)\oplus \H^1(X,\Omega^1_X) \cong k\oplus k^3 \cong k^4.
\end{align*}
\end{example}
\chapter{Weighted projective lines}
Throughout this chapter, we work over an algebraically closed field $k$ of characteristic zero, which we usually will omit from our notation.

In this chapter, we introduce weighted projective lines. A weighted projective line roughly
looks like the projective line $\mathbb P^1$ where each point gets attached a \emph{weight},
i.e.\ an integer $\geq 1$, such all but finitely many points have weight $1$.

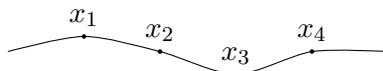
\begin{center}
\begin{tikzpicture}
\draw plot[smooth] coordinates { (0,0) (1,0.2) (2.0,0) (3,-0.3) (4,0) (5,0) };
\fill (1,0.2) circle[radius=0.03] node[above]{$x_1$};
\fill (2,0) circle[radius=0.03] node[above]{$x_2$};
\fill (3,-.3) circle[radius=0.03] node[above]{$x_3$};
\fill (4,0) circle[radius=0.03] node[above]{$x_4$};
\end{tikzpicture}
\end{center}

We fix finitely many distinct closed points $x_1,\dotsc,x_n\in\mathbb P^1$,
and attach to each $x_i$ a weight $a_i\in\mathbb Z_{\geq 1}$.

The points of weight $>1$ are called \emph{exceptional}, whereas the points of weight $1$ are called \emph{ordinary}.
We denote this weighted projective line by $\mathbb P^1\langle x_1,\dotsc,x_n; a_1,\dotsc,a_n\rangle$.
This is not a scheme (unless $n=0$), but an \emph{algebraic stack}. There are several ways of approaching
this stack, of which we will present two.

Since the weighted projective line $X=\mathbb P^1\langle x_1,\dotsc,x_n; a_1,\dotsc,a_n\rangle$ is a (sufficiently well-behaved)
algebraic stack, it comes with a category of coherent sheaves of modules, $\coh X$.
This is an abelian category, so we define the Hochschild cohomology and Hochschild homology
of $X$ to be the Hochschild (co)homology of the abelian category $\coh X$.

There are several occasions in mathematics where weighted projective lines appear naturally
(for a nice collection, see \cite{Lenzing2016}). We only mention the following one:
\begin{theorem}[\cite{Happel2001}]\label{HappelClassificationThm}
Let $\mathcal H$ be an abelian $k$-linear category such that
\begin{itemize}
\item $\mathcal H$ is \emph{hereditary}, i.e.\ $\Ext^2(\blank,\blank)$, defined as Yoneda-Ext, vanishes.
\item $\Hom(X,Y)$ and $\Ext^1(X,Y)$ are finite-dimensional for all objects $X,Y$ of $\mathcal H$.
\item $\mathcal H$ is \emph{connected}, i.e.\ it is not possible to decompose $\mathcal H = \mathcal H_1\oplus \mathcal H_2$
for non-zero abelian categories $\mathcal H_1, \mathcal H_2$.
\item $\mathcal H$ contains a \emph{tilting object}, i.e.\ some $T\in\mathcal H$ with $\Ext^1(T,T)=0$
and for all $X\in\mathcal H$, $\Hom(T,X)=0=\Ext^1(T,X)$ implies $X=0$.
\end{itemize}
Then $\mathcal H$ is derived equivalent to $H\Lmod$ for some hereditary algebra $H$
or to the category of coherent sheaves $\coh X$ for some weighted projective
line $X$.
\end{theorem}
\section{The approach of Geigle and Lenzing}\label{sec:GeigleLenzing}
Weighted projective lines were first defined in \cite{Geigle1987}. Their description does not use the terminology of algebraic stacks
and many steps appear to be rather ad-hoc, but this is sufficient for our purposes.
In this section, we summarize the (for our purposes) most important steps of \cite{Geigle1987},
omitting most of the proofs. 

Fix distinct points $x_1,\dotsc,x_n\in \mathbb P^1$ and integers $a_1,\dotsc,a_n\in\mathbb Z_{\geq 1}$. Put $\mathbf a := (a_1,\dotsc,a_n)$ and $\mathbf x:=(x_1,\dotsc,x_n)$.

In order to define the weighted projective line with weight $a_i$ at the point $x_i$, we may freely add or remove points
with weight $1$ to or from our list. We therefore assume $n\geq 2$. Moreover, in case $n>2$, we remove points of weight $1$.
In other words, we may assume the following convention:
\begin{convention}\label{conv:wpl}
We assume that one of the following holds:
\begin{itemize}
\item
$n=2$
\item $n>2$ and $a_i>1$ for $i=1,\dotsc,n$
\end{itemize}
Since $\PGL_2$ acts triply transitive on $\mathbb P^1$, we may moreover choose coordinates
on $\mathbb P^1$ such that $x_1=\infty$, $x_2=0$ and, in case $n\geq 3$, $x_3=1$.
\end{convention}
\begin{definition}\label{def:GeigleLenzing}The weighted projective line $X$ can be described as follows:
\begin{enumerate}[(i)]
\item The group $L(\mathbf a)$ is defined as the free abelian group in generators $\ell_1,\dotsc,\ell_n$ modulo the relations
\begin{align*}
a_1 \ell_1 = a_2 \ell_2 = \cdots = a_n \ell_n.
\end{align*}
The subset of elements of the form $\sum\limits_{i=1}^n \lambda_i \ell_i$ with all $\lambda_i\geq 0$ is denoted $L_+(\mathbf a)$.
We turn $L(\mathbf a)$ into an ordered group by saying $\ell\geq 0$ iff $\ell\in L_+(\mathbf a)$ for all $\ell\in L(\mathbf a)$.
\item We denote by
\begin{align*}
S(\mathbf a) = k[X_1,\dotsc,X_n]
\end{align*}
the $L(\mathbf a)$-graded algebra with $\deg X_i = \ell_i$. The quotient
\begin{align*}
S(\mathbf a,\mathbf x) := S(\mathbf a)/I(\mathbf a,\mathbf x),\qquad I(\mathbf a,\mathbf x) = (X_i^{a_i} - X_2^{a_2} + x_i X_1^{a_1}\text{ for }i=3,\dotsc,n)
\end{align*}
is also $L(\mathbf a)$-graded.
\item Denote
\begin{align*}
&G(\mathbf a) := \{ (t_1,\dotsc,t_n) \in (k^\times)^{n} \mid t_1^{a_1} = \cdots = t_n^{a_n}\}.
\\&F(\mathbf a,\mathbf x) := \{ (z_1,\dotsc,z_n) \in \mathbb A_n\setminus\{0\} \mid z_i^{a_i} = z_2^{a_2} - x_i z_1^{a_1}\text{ for }i=3,\dotsc,n\}.
\end{align*}
$G(\mathbf a)$ acts on $F(\mathbf a,\mathbf x)$ by $(t_1,\dotsc,t_n).(z_1,\dotsc,z_n) := (t_1z_1,\dotsc,t_n z_n)$.
Denote the (set-theoretic) quotient $F(\mathbf a,\mathbf x) / G(\mathbf a)$ by $X$.
\item For an $L(\mathbf a)$-homogeneous element $f\in S(\mathbf a)$, denote by
\begin{align*}
V(f) := \{[x]\in X \mid f(x)=0\}
\end{align*}
the vanishing set of $f$. As each element of $I(\mathbf a,\mathbf x)$ vanishes on $X$,
we have the same definition of $f\in S(\mathbf a,\mathbf x)$.
We equip $X$ with a topology where the closed subsets are those of the form $V(f)$ for some $L(\mathbf a)$-homogeneous
element $f\in S(\mathbf a,\mathbf x)$.
\item For an $L(\mathbf a)$-homogeneous element $f\in S(\mathbf a,\mathbf x)$, denote
by $S(\mathbf a,\mathbf x)_f$ the localization, which is again $L(\mathbf a)$-graded.
Denote by $\mathcal O_X$ the sheaf on $X$ attached to the presheaf
\begin{align*}
(X\setminus V(f)) \mapsto S(\mathbf a,\mathbf x)_f
\end{align*}
of $L(\mathbf a)$-graded $k$-algebras.
\item Without changing the situation much, we may as well add a \emph{generic point} $\xi$ to the topological space $X$.
The stalk of $\mathcal O_X$ is the localization of $S(\mathbf a,\mathbf x)$ at all $L(\mathbf a)$-homogeneous nonzerodivisors
of $S(\mathbf a,\mathbf x)$. The resulting space is called \emph{weighted projective line} with weights $\mathbf a$
at the points $\mathbf x$.
\item Denote by $\Mod^{L(\mathbf a)}_X$ the category of $L(\mathbf a)$-graded sheaves of $\mathcal O_X$-modules on $X$.
Each $\ell\in L(\mathbf a)$ induces a shift on $\Mod^{L(\mathbf a)}_X$, where
for $M\in \Mod^{L(\mathbf a)}_X$, the \emph{shifted module} $M(\ell)$ is given by $M(\ell)_\mu = M_{\ell+\mu}$ for all $\mu\in L(\mathbf a)$.
\item A sheaf of modules $M\in \Mod^{L(\mathbf a)}_X$ if called \emph{quasi-coherent} if $X$ can be covered by open subsets $U$ such that
there exist families $(l_i)_{i\in I}, (l_j)_{j\in J}$ of elements in $L(\mathbf p)$ and an exact sequence
\begin{align*}
\bigoplus_{i\in I} \mathcal O_X(l_i)\vert_U
\rightarrow
\bigoplus_{j\in J} \mathcal O_X(l_j)\vert_U
\rightarrow
M\vert_U\rightarrow 0
\end{align*}
of sheaves of modules on $U$.
$M$ is called \emph{coherent} if $X$ can be covered by such open subsets $U$ such that a presentation
as above exists with $I$ and $J$ being finite sets. We denote by $\coh X$ the category of coherent sheaves on $X$.
\end{enumerate}
\end{definition}
It turns out that $\coh X$ satisfies all assumptions of Theorem \ref{HappelClassificationThm}.
In particular $\coh X$ is an abelian hereditary category. It satisfies a version of Serre duality \cite[Prop.~2.2]{Geigle1987}.
However, most crucial for us is the existence of a tilting object.
\begin{proposition}[{\cite[Proposition~4.1]{Geigle1987}}]\label{prop:wplTiltingObject}
Denote $c := a_1 \ell_1 = \cdots = a_n \ell_n\in L(\mathbf a)$. Then\footnote{In Definition~\ref{def:exceptional},
we required the elements of the exceptional collection to be totally ordered. In this situation,
$\leq$ defines only a partial order. The proposition holds for \emph{any} total
order on $\{\mathcal O_X(\ell):0\leq \ell\leq c\}$ refining $\leq$.}
\begin{align*}
\{\mathcal O_X(\ell) : 0\leq \ell\leq c\}
\end{align*}
is an exceptional collection of $\D^\bounded(\coh X)$.
In particular,
\begin{align*}
T := \bigoplus_{0\leq \ell \leq c} \mathcal O_X(\ell)
\end{align*}
is a tilting object in $\D^\bounded(\coh X)$.\rightqed\end{proposition}
Put $A := \End_{\coh X}(T)$. By Proposition \ref{prop:abelianTilting},
the Hochschild (co)homology of the weighted projective line $X$ is
the Hochschild (co)homology of the finite-dimensional algebra $A$.
Moreover, we can describe $A$ explicitly:
\begin{proposition}\label{proposition-canonical}
$A$ is isomorphic to the bound quiver algebra $kQ / I$, where the quiver $Q$ and the ideal $I$ are defined as follows:

$Q$ has as set of vertices
\begin{align*}
Q_0 = \{m \ell_i : 0\leq m \leq a_i, \quad i=1,\dotsc,n\}
\end{align*}
(recall $0\ell_i = 0 \ell_j$ and $a_i \ell_i = a_j \ell_j$ for all $i,j$).

For vertices $v,w\in Q_0$, there is
\begin{itemize}
\item An arrow from $v$ to $w$ labelled $X_i$ in case $w=v+\ell_i$ (for $i\in\{1,\dotsc,n\}$).
\item No arrow from $v$ to $w$ otherwise.
\end{itemize}

$Q$ can be depicted as follows:
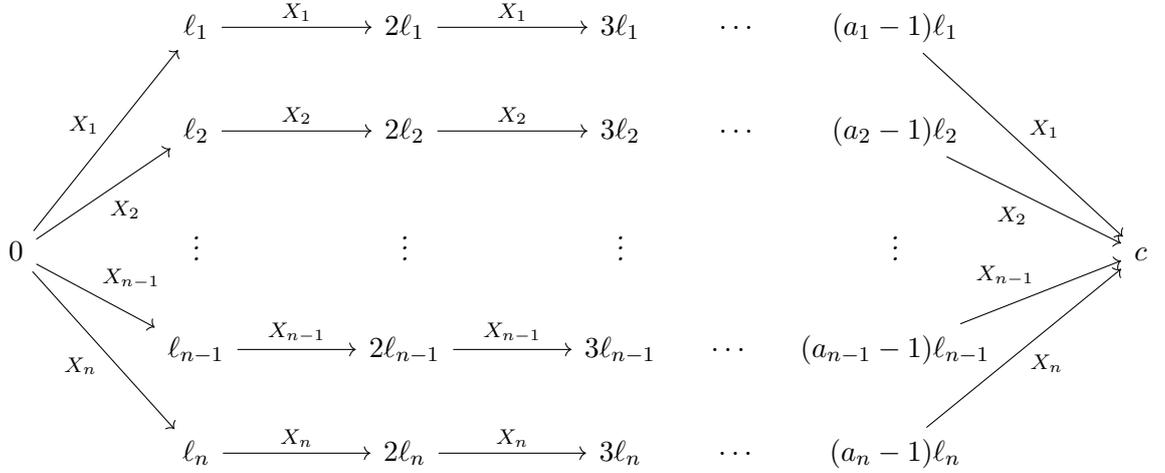
\begin{align*}
\begin{tikzcd}[ampersand replacement=\&,column sep=4em]
\&\ell_1\ar[r,"X_1"]\&2\ell_1\ar[r,"X_1"]\&3\ell_1\ar[r,phantom,"\hdots"]\&(a_1-1) \ell_1\ar[ddr,"X_1"]\\
\&\ell_2\ar[r,"X_2"]\&2\ell_2\ar[r,"X_2"]\&3\ell_2\ar[r,phantom,"\hdots"]\&(a_2-1)\ell_2\ar[dr,"X_2"']\\
0\ar[ruu,"X_1"]\ar[ru,"X_2"']\ar[dr,"X_{n-1}"]\ar[ddr,"X_n"']\&\vdots\&\vdots\&\vdots\&\vdots\&c\\
\&\ell_{n-1}\ar[r,"X_{n-1}"]\&2 \ell_{n-1}\ar[r,"X_{n-1}"]\&3 \ell_{n-1}\ar[r,phantom,"\hdots"]\&(a_{n-1}-1)\ell_{n-1}\ar[ur,"X_{n-1}"]\\
\&\ell_n\ar[r,"X_n"]\&2 \ell_n\ar[r,"X_n"]\&3 \ell_n\ar[r,phantom,"\hdots"]\&(a_n-1) \ell_n\ar[uur,"X_n"']
\end{tikzcd}
\end{align*}
It is for pictorial reasons called \emph{hammock quiver}.

In case $n\leq 2$, the ideal $I\subseteq kQ$ is defined to be zero, so $A\cong kQ$.
Otherwise, we denote by $X_i^{a_i}$ the path in $Q$ obtained by composing all the $a_i$ arrows with label $X_i$ (this is a path
starting at $0$ and ending in $c$).
Then $I$ is generated by the elements
\begin{align*}
X_i^{a_i} - X_2^{a_2} + x_i X_1^{a_1}
\end{align*}
for $i=3,\dotsc,n$.
\end{proposition}
\begin{definition}
The \emph{canonical algebra} associated to the points $\mathbf x$ and the weights $\mathbf a$
is by definition $kQ/I$, where $Q$ and $I$ are defined as in Proposition \ref{proposition-canonical}.
\end{definition}
\begin{corollary}\label{cor:wplCanonicalHochschild}
The weighted projective line and the canonical algebra associated to the points $\mathbf x$
and weights $\mathbf a$ have isomorphic Hochschild homology and cohomology groups.\rightqed
\end{corollary}
\begin{example}
Consider the case $n=2$ and $a_1=a_2=1$. We would expect to get the usual $\mathbb P^1$ back.

Indeed, $L(\mathbf a) = \mathbb Z\ell_1 = \mathbb Z\ell_2$ is isomorphic to $\mathbb Z$,
so that $S(\mathbf a)\cong k[X,Y]$ with the normal monomial grading. Our construction above yields exactly $\mathbb P^1$
and $\Qcoh(\mathbb P^1)$ resp. $\coh(\mathbb P^1)$. The hammock quiver is given by
\begin{align*}\begin{tikzcd}[ampersand replacement=\&]\bullet\ar[r,shift left=.3em]\ar[r,shift right=.3em]\&\bullet\end{tikzcd}
\end{align*}
and there are no relations for the canonical algebra. We recognize this from Example \ref{ExampleDerivedEquivalenceProjectiveSpace}.

If we have precisely one exceptional point, i.e.\ $X=\mathbb P^1\langle \infty; a_1\rangle$,
the corresponding canonical algebra has the quiver
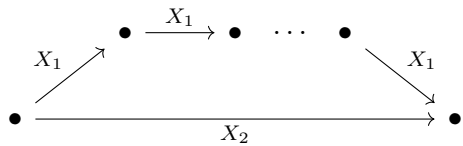
\begin{align*}
\begin{tikzcd}[ampersand replacement=\&]
\&\bullet\ar[r,"X_1"]\&\bullet\ar[r,phantom,"\hdots"]\&\bullet\ar[dr,"X_1"]\\
\bullet\ar[ru,"X_1"]\ar[rrrr,"X_2"']\&\&\&\&\bullet
\end{tikzcd}
\end{align*}
and no relations.

As a final example, if $X=\mathbb P^1\langle \infty,0,1; 2,2,3\rangle$,
its canonical algebra is defined by the quiver
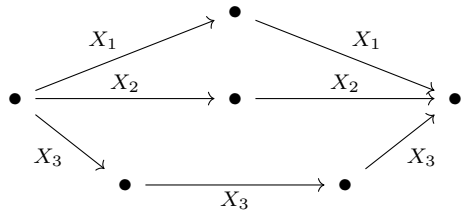
\begin{align*}
\begin{tikzcd}[ampersand replacement=\&]
\&\&\bullet\ar[drr,"X_1"]\\
\bullet\ar[urr,"X_1"]\ar[rr,"X_2"]\ar[dr,"X_3"']\&\&\bullet\ar[rr,"X_2"]\&\&\bullet\\
\&\bullet\ar[rr,"X_3"']\&\&\bullet\ar[ur,"X_3"']
\end{tikzcd}
\end{align*}
modulo the relation
\begin{align*}
X_3^3 - X_2^2 + X_1^2.
\end{align*}
\end{example}
While the approach of Geigle and Lenzing allows to easily write down what a weighted projective
line \enquote{is}, its nature is rather ad-hoc and allows little geometric intuition.
In the following section, we sketch how one can construct weighted projective lines
using the language of stacks.
\section{Weighted projective lines as iterated root stacks}
The idea of root stacks is that given a stack $\mathcal X$ (e.g.\ a scheme) and an effective Cartier divisor $D\subseteq \mathcal X$,
we want to add an automorphism group cyclic of order $r\geq 1$ to $D$ and not change $\mathcal X$ outside of $D$.

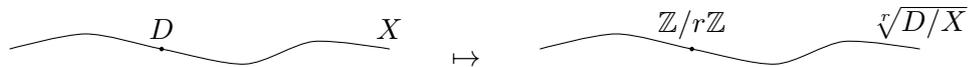
\begin{center}\begin{tikzpicture}
\draw plot[smooth] coordinates { (0,0) (1,0.2) (2.0,0) (3.1,-0.2) (4,0.1) (5,0) };
\draw (5,0) node[above]{$X$};
\fill (2,0) circle[radius=0.03] node[above]{$D$};
\end{tikzpicture}
\quad$\mapsto$\qquad
\begin{tikzpicture}
\draw plot[smooth] coordinates { (0,0) (1,0.2) (2.0,0) (3.1,-0.2) (4,0.1) (5,0) };
\draw (5,0) node[above]{$\sqrt[r]{D/X}$};
\fill (2,0) circle[radius=0.03] node[above]{$\mathbb Z/r\mathbb Z$};
\end{tikzpicture}\end{center}

If $\mathcal X$ is one-dimensional, this means that we may attach the weight $r\geq 1$ to some point $D\in \mathcal X$
and leave the other points of $\mathcal X$ the same.
Starting with $\mathcal X=\mathbb P^1$, one can get any weighted projective line by iterating this construction.

It is beyond the scope of this thesis to give a detailled account to the notion of algebraic stacks;
we refer the reader to \cite{Olsson2016} for a good introduction.

Root stacks and iterated root stacks were introduced in \cite{Cadman2007}.
We won't repeat his definition and the necessary prerequisites here for space reasons.
Instead, we will illustrate some properties of this construction.
Namely, a root stack comes with an interesting semi-orthogonal decomposition.
\begin{definition}[{\cite[Def.\ 1.42 and 1.57]{Huybrechts2006}}]
Let $\mathcal D$ be a triangulated category.
\begin{enumerate}[(a)]
\item A full triangulated subcategory $\mathcal D'\subseteq \mathcal D$ is called \emph{admissible} if
the inclusion functor $\mathcal D'\hookrightarrow \mathcal D$ has a right adjoint.
\item A sequence $\mathcal D_1,\dotsc,\mathcal D_n$ of admissible subcategories of $\mathcal D$
is called \emph{semi-orthogonal} if for all $i<j$ and $X\in \mathcal D_i,Y\in\mathcal D_j$,
we have $\Hom_{\mathcal D}(Y,X)=0$.
\item A semi-orthogonal sequence $\mathcal D_1,\dotsc,\mathcal D_n$ is called \emph{semi-orthogonal decomposition}
of $\mathcal D$ if the smallest full triangulated subcategory of $\mathcal D$ containing all $\mathcal D_i$ is $\mathcal D$.

In this case we write
\begin{align*}
\mathcal D\cong\langle \mathcal D_1,\dotsc,\mathcal D_n\rangle.
\end{align*}
\end{enumerate}
\end{definition}
\begin{remark}
Full exceptional sequences (cf.\ Definition~\ref{def:exceptional})
give semi-orthogonal decompositions by \cite[Example~1.60]{Huybrechts2006}.
\end{remark}
\begin{proposition}[{\cite[Theorem~4.7]{Bergh2016}}]\label{prop:rootStackSemiOrthogonal}
Let $X$ be a stack, $D\subseteq X$ an effective Cartier divisor, $r\geq 1$ an integer
and $\tilde X$ be the root stack.

The derived category $\D^\bounded(\tilde X)$ allows a semi-orthogonal decomposition
of the form
\begin{align*}
\D^\bounded(\tilde X) =\Bigl\langle \D^\bounded(X),\underbrace{\D^\bounded(D),\dotsc,\D^\bounded(D)}_{r-1\text{ times}}\Bigr\rangle.
\end{align*}
(i.e.\ a semi-orthogonal decomposition into subcategories equivalent to $\D^\bounded(X)$ resp.\ $\D^\bounded(D)$).\rightqed
\end{proposition}
\begin{example}\label{ex:squidOneLeg}
Consider the case $X=\mathbb P^1$. $\mathbb P^1$ is derived equivalent to the Kronecker algebra
$k\left(\begin{tikzcd}[column sep=1em]\bullet\ar[r,shift left=.3em]\ar[r,shift right=.3em]&\bullet\end{tikzcd}\right)$.
If $E$ is the Cartier divisor associated to a closed point $x\in \mathbb P^1$, the resulting root stack (a weighted projective
line where $x$ has weight $r$ and all other points have weight $1$) is derived equivalent to
the path algebra over the quiver
\begin{align*}
\begin{tikzcd}[ampersand replacement=\&]
\bullet\ar[r,shift left=.3em]\ar[r,shift right=.3em]\&\bullet\ar[r,"x_1"]\&\bullet\&\hdots\&\bullet\ar[r,"x_{r-1}"]\&\bullet.
\end{tikzcd}
\end{align*}
We identify the left end as $\D^\bounded(\mathbb P^1)$ and the $r-1$ following points each as $\D^\bounded(E)$.
\end{example}
\begin{proposition}\label{prop:squid}
Let $X$ be the weighted projective line assocated to the points $x_1,\dotsc,x_n\in\mathbb P^1$
and the weights $a_1,\dotsc,a_n\in\mathbb Z_{\geq 1}$. Then $X$ is derived equivalent
to the algebra $A=kQ/I$ defined as follows:
The quiver $Q$ is given by
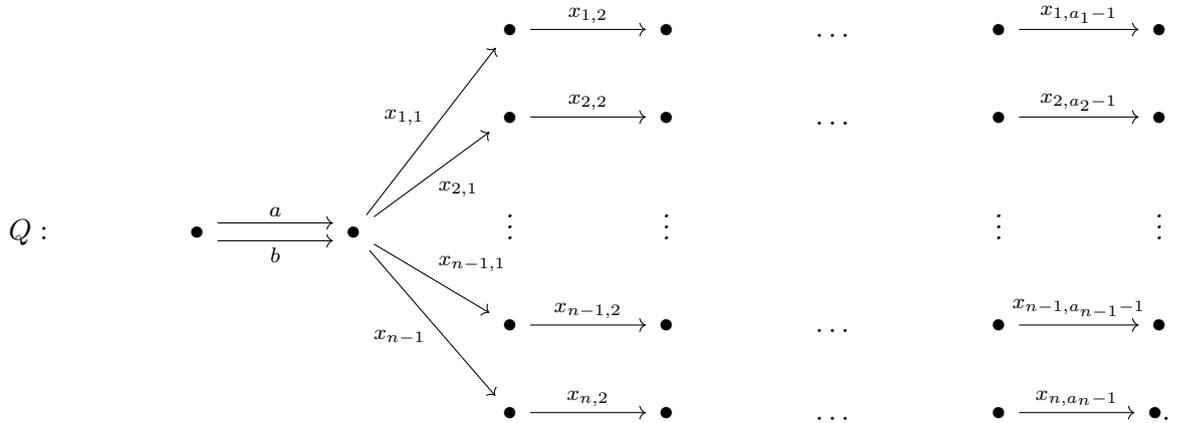
\begin{align*}
\begin{tikzcd}[ampersand replacement=\&,column sep=3.9em]
\&\&\&\bullet\ar[r,"x_{1,2}"]\&\bullet\&\hdots\&\bullet\ar[r,"x_{1,a_1-1}"]\&\bullet\\
\&\&\&\bullet\ar[r,"x_{2,2}"]\&\bullet\&\hdots\&\bullet\ar[r,"x_{2,a_2-1}"]\&\bullet\\
Q:\&\bullet\ar[r,"a",shift left=0.3em]\ar[r,"b"',shift right=.3em]\&\bullet\ar[ruu,"x_{1,1}"]\ar[ru,"x_{2,1}"']\ar[rd,"x_{n-1,1}"]\ar[rdd,"x_{n-1}"']\&\vdots\&\vdots\&\&\vdots\&\vdots\\
\&\&\&\bullet\ar[r,"x_{n-1,2}"]\&\bullet\&\hdots\&\bullet\ar[r,"x_{n-1,a_{n-1}-1}"]\&\bullet\\
\&\&\&\bullet\ar[r,"x_{n,2}"]\&\bullet\&\hdots\&\bullet\ar[r,"x_{n,a_n-1}"]\&\bullet.
\end{tikzcd}
\end{align*}
The ideal $I$ is generated by
\begin{align*}
ax_{1,1},\quad bx_{2,1}, \quad (x_i a - b) x_{i,1}\text{ for }i=3,\dotsc,n.
\end{align*}
\end{proposition}
\begin{proof}
This can be shown as in Example~\ref{ex:squidOneLeg} by iterating the root stack construction,
or a tiltling approach similar to Proposition~\ref{prop:wplTiltingObject}.
The latter is explained e.g.\ in \cite[Example~3.2]{Plamondon2016}.
\end{proof}
For pictorial reasons, the algebra $A$ in Proposition~\ref{prop:squid} is called \emph{squid algebra}.
The leftmost \enquote{head} of the squid  $\begin{tikzcd}[column sep=1em]\bullet\ar[r,shift left=.3em,"a"]\ar[r,shift right=.3em,"b"']&\bullet\end{tikzcd}$
comes from $\D^\bounded(\mathbb P^1)$, and the other vertices correspond to the exceptional points.
\chapter{Happel's method}\label{chapter-3}
By Corollary~\ref{cor:wplCanonicalHochschild}, the Hochschild (co)homology of a weighted projective line can be calculated
by calculating the Hochschild (co)homology of the associated canonical algebra\footnote{Of course, one can equivalently
study squid algebras. The methods described in this chapter work similarly for squid algebras.}.

In \cite{Happel1989}, Happel
proves important facts on the Hochschild cohomology of finite-dimensional associative algebras, which are very useful
in practical computations. However, it is not a full procedure to compute the Hochschild cohomology. It can easily be adapted for Hochschild homology.

In the sequel \cite{Happel1998}, Happel gives a calculation of the Hochschild cohomology of canonical algebras.

In this chapter, we reproduce Happel's work and give a procedure to compute the Hochschild cohomology of finite-dimensional
algebras of global dimension $\leq 2$.

Throughout this chapter, let $k$ be any field.
\section{The minimal resolution}\label{sec:minimalResolution}
In this section, we fix a finite quiver $Q$ and an admissible ideal $I\subseteq kQ$.
Denote by $A:=kQ/I$ the bound quiver algebra\footnote{Our conventions are the following:
All modules are, unless stated otherwise, left modules. If $p$ and $q$ are paths in $kQ$,
we denote by $pq$ the path described by first following $p$ and then $q$; it is zero if the end vertex of $p$
is distinct from the starting vertex of $q$. Note that this is the opposite composition from \cite{Assem2006},
so we translate the results cited from there accordingly.}. Recall that finite-dimensional algebra is Morita
equivalent to a bound quiver algebra where the quiver $Q$ and the ideal $I$ are almost uniquely determined \cite[Section~II.3]{Assem2006}.

Enumerate the vertices of $Q$ by $e_1,\dotsc,e_n$. For each $i\in\{1,\dotsc,n\}$,
let $P(i) := Ae_i$ be the associated projective and $S(i) = P(i)/\rad(P(i))$ the associated simple module.

Denote the set of arrows in $Q$ by $Q_1$. For each path $p$ in $Q$, let $s(p)$ resp.\ $e(p)$
denote the starting resp.\ end vertex of $p$.

By definition, $\HH^\bullet(A) = \Ext^\bullet_{A^\env}(A,A)$. We want to compute this by using a projective resolution
of $A$ as $A^\env$-module. Recall that by \cite[Corollary~I.5.10]{Assem2006}, $A$ allows a minimal projective resolution, unique up to isomorphism\footnote{Uniqueness
is not stated explicitly there, but follows using \cite[Theorem~I.5.8]{Assem2006}.}.
We can describe it as follows:
\begin{proposition}[{\cite[Lemma~1.5]{Happel1989}}]\label{prop:happel89}
$A$ allows a minimal projective resolution as $A^\env\text{-}$module
\begin{align*}
\cdots\rightarrow P_m\xrightarrow{d_m} P_{m-1}\rightarrow\cdots\rightarrow P_1\xrightarrow{d_1} P_0\xrightarrow{\varepsilon} A\rightarrow 0,
\end{align*}
where for $n\geq 0$,
\begin{align*}
P_m = \bigoplus_{i,j=1}^n (A e_j\otimes_k e_i A)^{\dim \Ext^m_A(S(i),S(j))}.
\end{align*}
$Ae_i\otimes_ke_jA$ becomes an $A$-bimodule by multiplication from the left resp.\ right.
\end{proposition}
This describes precisely the objects in the projective resolution using only the Ext-groups of the simple modules.
However, there is no description of the differentials.

Before we use this to compute Hochschild cohomology groups, let us note the following consequence:
\begin{corollary}\label{cor:cohomologyVanishesAboveGlDim}
If $A$ has finite global dimension $g$, and $i>g$, then $\HH^i(A)=0$.
\end{corollary}
\begin{proof}
By Proposition~\ref{prop:happel89}, the projective dimension of $A$ as $A^\env$-module is $\leq g$,
so in particular $\Ext^i_{A^\env}(A,A)=0$.
\end{proof}
We can find a description of the differentials in small degrees as follows:
We will construct an exact sequence of $A^\env$-modules
\begin{align*}
C_1\xrightarrow{c_1} C_0\xrightarrow{c_0} A\rightarrow 0
\end{align*}
with $C_1\cong P_1, C_0\cong P_0$ and an explicit description of $c_0,c_1$. By uniqueness
of the minimal projective resolutions, we find a minimal projective
resolution as in Proposition~\ref{prop:happel89} with $P_0=C_0, P_1=C_1, \varepsilon=c_0$ and $d_1=c_1$.

First, we recall some technical preliminaries.
\begin{lemma}\label{lem:trivialFactsFDAlgebras}
Let $i,j\in\{1,\dotsc,n\}$.
\begin{enumerate}[(a)]
\item For any $A$-bimodule $M$, there is an isomorphism
\begin{align*}
&\Hom_{A^\env}(Ae_i\otimes_k e_j A, M)\xrightarrow\cong e_i M e_j,
\\& \varphi\mapsto \varphi(e_i\otimes e_j)
\\& \left[ ae_i\otimes e_j b\mapsto a m b\right]\mapsfrom m
\end{align*}
Moreover, there is an isomorphism
\begin{align*}
&(Ae_i\otimes_k e_j A)\otimes_{A^\env} M\xrightarrow\cong e_j M e_i
\\&(a\otimes b)\otimes m\mapsto bma
\\&(e_i\otimes e_j)\otimes m\mapsfrom m.
\end{align*}
\item It holds that
\begin{align*}
\dim\Ext^0_A(S(i),S(j)) =&\begin{cases}1,&i=j\\
0,&i\neq j,\end{cases}
\\\dim\Ext^1_A(S(i),S(j)) =&\#\{\alpha\in Q_1\mid s(\alpha)=e_j, e(\alpha)=e_i\}.
\end{align*}
\item Denote by $Q_0$ the subquiver of $Q$ which contains all vertices but no arrows.
Then the canonical projection $kQ\rightarrow A$ has a $kQ_0$-bilinear section.
\end{enumerate}
\end{lemma}
\begin{proof}
(a) is trivial. (b) follows from Schur's lemma \cite[Lemma~I.3.1]{Assem2006} and \cite[Lemma~III.2.12]{Assem2006}.

(c) follows since $kQ_0$ and therefore $(kQ_0)^\env$ is semi-simple \cite[Proposition~II.1.10]{Assem2006}.
\end{proof}
\begin{proposition}\label{prop:algebraHomology}
Suppose that the quiver $Q$ is acylic\footnote{Which is equivalent to $A$ having finite global dimension, but we won't use this fact.}.
The Hochschild homology of $A$ in degree $i\in\mathbb Z$ is given by
\begin{align*}
\dim\HH_i(A) = \begin{cases}n,&i=0\\0,&i\neq 0.\end{cases}
\end{align*}
\end{proposition}
\begin{proof}
Consider the projective resolution of $A$ as $A^\env$-module specified in Proposition~\ref{prop:happel89}.
We claim that $P_m\otimes_{A^\env} A=0$ for $m>0$:

Consider a simple module $S(i)$ and a minimal projective resolution as left $A$-module
\begin{align*}
\cdots\rightarrow T_s\rightarrow T_{s-1}\rightarrow\cdots\rightarrow T_0\rightarrow S(i)\rightarrow 0.
\end{align*}
Let $V$ be the set of vertices $v\in Q_0$ such that there exists a path from $v$ to $e_i$ in $Q$.
By induction, using minimality of the projective resolution, one shows that each $T_s$
has the form $T_s\cong \bigoplus_{v\in V} (Av)^{a_{s,v}}$.

In particular, if $e_j\notin V$, we have $\Hom_A(T_s,S(j))=0$.
Hence $\Ext^s_A(S(i),S(j))=0$ unless there is a path from $j$ to $i$.
Moreover, if $s>0$, we also have $\Ext^s_A(S(i),S(i))=0$.

By Lemma~\ref{lem:trivialFactsFDAlgebras},
\begin{align*}
P_m\otimes_{A^\env} A
\cong\bigoplus_{i,j=1}^n (e_i A e_j)^{\dim\Ext^m_A(S(i),S(j))}.
\end{align*}
Note that if $m>0$ and $\Ext^m_A(S(i),S(j))\neq 0$,
there exists a path from $i$ to $j$ in $Q$ which has positive length.
By assumption, there will be no path from $i$ to $j$ in $Q$, so $e_i A e_j=0$.

This shows $P_m\otimes_{A^\env} A = 0$ for $m>0$. Moreover,
\begin{align*}
P_0\otimes_{A^\env} A \cong\bigoplus_{i=1}^n e_i A e_i =\bigoplus_{i=1}^n k e_i,
\end{align*}
which has dimension $n$.
\end{proof}
\begin{definition}
Define $C_0$ and $C_1$
to be the $A^\env$-modules
\begin{align*}
C_0 :=&\bigoplus_{i=1}^n A e_i\otimes_k e_i A,
\\C_1 :=&\bigoplus_{\alpha\in Q_1} A s(\alpha)\otimes_k e(\alpha) A.
\end{align*}
For $i\in \{1,\dotsc,n\}$, let $\incl_{e_i} : A e_i\otimes_k e_i A\hookrightarrow C_0$ be the canonical inclusion.
Similarly for $\alpha\in Q_1$, let $\incl_\alpha : As(\alpha)\otimes_ke(\alpha)A\hookrightarrow C_1$ be the canonical inclusion.

We define the $A^\env$-linear maps
\begin{align*}
C_1\xrightarrow{c_1}C_0\xrightarrow{c_0} A
\end{align*}
as follows: For $i\in\{1,\dotsc,n\},\alpha\in Q_1$ and $a,b\in A$, we define
\begin{align*}
&c_0(\incl_{e_i}(a\otimes b)) := ab
\\& c_1(\incl_\alpha(a\otimes b)) := \incl_{s(\alpha)}(a\otimes\alpha b) - \incl_{e(\alpha)}(a\alpha \otimes b).
\end{align*}
\end{definition}
\noindent{}Note that by Lemma~\ref{lem:trivialFactsFDAlgebras}, $C_0\cong P_0$ and $C_1\cong P_1$.
\begin{lemma}\label{lem:minProjective1}~
\begin{enumerate}[(a)]
\item The $k$-linear map
\begin{align*}
\partial : kQ\rightarrow C_1,
\end{align*}
sending a path $p=\alpha_1\cdots\alpha_r$ to
\begin{align*}
\sum_{i=1}^r \incl_{\alpha_i}(\alpha_1\cdots\alpha_{i-1}s(\alpha_i)\otimes e(\alpha_i)\alpha_{i+1}\cdots\alpha_n)
\end{align*}
satisfies the Leibnitz rule
\begin{align*}
\forall p,q\in kQ: \partial(pq) = p\partial(q)+q\partial(p).
\end{align*}
\item For a path $p$ in $Q$, we have
\begin{align*}
c_1 \partial(p) = \incl_{s(p)}(s(p)\otimes p) - \incl_{e(p)}(p\otimes e(p)).
\end{align*}
\item The morphism of left $A$-modules $s_{-1} : A\rightarrow C_0$, defined
by $s_{-1}(e_i) := e_i\otimes e_i$ and $A$-linear continuation,
is a splitting of $c_0$.
\item Fix a $kQ_0$-linear splitting $\varphi : A\rightarrow kQ$ as in Lemma~\ref{lem:trivialFactsFDAlgebras}
and define the morphism of left $A$-modules
\begin{align*}
s_0 : C_0\rightarrow C_1
\end{align*}
for $i\in\{1,\dotsc,n\}$, $a\in Ae_i$ and $b\in e_i A$ by
\begin{align*}
s_0(\incl_{e_i}(a\otimes b)) := a\cdot\partial(\varphi(b)).
\end{align*}
Then the sequence of left $A$-modules
\begin{align*}
C_1\xrightarrow{c_1}C_0\xrightarrow{c_0}A\rightarrow 0
\end{align*}
is contractible via the homotopy $s_{-1},s_0$.
\end{enumerate}
\end{lemma}
\begin{proof}
For (a) and (b), one uses a simple induction on path lengths. (c) is trivial.

For (d), we have to show $s_{-1} c_0 + c_1 s_0 = \id_{C_0}$. Since both sides are morphisms
of left $A$-modules, we may check equality on elements of the form $\incl_{e_i}(e_i\otimes b)$ for $b\in e_i A$.
Write $\varphi(b) := \sum\limits_{i=1}^r \lambda_i p_i$ for paths $p_1,\dotsc,p_r$ in $Q$ (all of them starting in $e_i$). Then
\begin{align*}
c_1 s_0(\incl_{e_i}(e_i\otimes b) =&c_1(\partial(\varphi(b))) =\sum_{i=1}^r \lambda_i c_1(\partial(p_i))
\\\underset{\text{(b)}}=&\sum_{i=1}^r \lambda_i \left(\incl_{s(p_i)} (s(p_i)\otimes p_i) - \incl_{e(p_i)}(p_i\otimes e(p_i))\right)
\\\underset{\incl_\bullet\text{ is $A^\env$-linear}}=&\incl_{e_i}(e_i\otimes\sum_{i=1}^r \lambda_i p_i)
-\sum_{i=1}^r \lambda_ip \incl_{e(p)}(e(p)\otimes e(p))
\\\underset{\text{Def. }s_{-1}}=&\incl_{e_i}(e_i\otimes b) - \sum_{i=1}^r s_{-1}(\lambda_i p_i)
\\=&\incl_{e_i}(e_i\otimes b) - s_{-1}(b) = \incl_{e_i}(e_i\otimes b) - s_{-1} c_0(\incl_{e_i}(e_i\otimes b)).\qedhere
\end{align*}
\end{proof}
\begin{corollary}\label{cor:minProjective1}
$A$ has a minimal projective resolution as $A^\env$-module
\begin{align*}
\cdots \rightarrow P_m\xrightarrow{d_m}P_{m-1}\rightarrow\cdots \rightarrow P_2\rightarrow C_1\xrightarrow{c_1}C_0\xrightarrow{c_0}A\rightarrow 0
\end{align*}
\end{corollary}
\begin{proof}
From Proposition~\ref{prop:happel89}, Lemma~\ref{lem:minProjective1} and uniqueness of minimal projective resolutions.
\end{proof}
\begin{corollary}[{\cite[Proposition~1.6]{Happel1989}}]\label{cor:hochschildCohomologyHereditary}~
\begin{enumerate}[(a)]
\item
$\dim\HH^0(A)$ is the number of path-connected components of $Q$.
\item
For each arrow $\alpha\in Q_1$, let $\nu(\alpha):=\dim_k s(\alpha)Ae(\alpha)$. Then
\begin{align*}
\dim\HH^1(A)\leq \dim \HH^0(A)-n + \sum_{\alpha\in Q_1} \nu(\alpha).
\end{align*}
If $A$ has global dimension $\leq 1$, we have equality.
\end{enumerate}
\end{corollary}
\begin{proof}
Apply $\Hom_{A^\env}(\blank,A)$ to the minimal projective resolution from Corollary~\ref{cor:minProjective1}.
We get the sequence
\begin{align*}
0\rightarrow \Hom_{A^\env}(C_0,A)\xrightarrow{c_1^\ast} \Hom_{A^\env}(C_1,A)\rightarrow \Hom_{A^\env}(P_2,A).
\end{align*}
By Lemma~\ref{lem:trivialFactsFDAlgebras}, we can identify this with the sequence
\begin{align*}
0\rightarrow \bigoplus_{i=1}^n e_i A e_i\xrightarrow\phi \bigoplus_{\alpha\in Q_1} s(\alpha) A e(\alpha)\rightarrow \Hom_{A^\env}(P_2,A).
\end{align*}
Note that $e_i A e_i=ke_i$. Here, $\phi$ is the map which sends $e_i$ to
\begin{align*}
\sum_{s(\alpha)=e_i} \alpha - \sum_{e(\alpha)=e_i} \alpha.
\end{align*}
Thus if $Q = V_1\cup\cdots\cup V_r$ is the decomposition of $Q$ into its path-connected components,
a basis of $\ker(\phi)$ is given by
\begin{align*}
\Bigl\{\sum_{\substack{i=1\\e_i\in V_j}}^n e_i \mid j=1,\dotsc,r\Bigr\}.
\end{align*}
Thus $\dim\HH^0(A)=\dim\ker(\phi)=r$. 

Moreover,
\begin{align*}
&\dim \HH^1(A) = \dim\ker(\Hom_{A^\env}(C_1,A)\rightarrow\Hom_{A^\env}(P_2,A)) - \dim\im(\phi)
\\\leq&\dim\Hom_{A^\env}(C_1,A) - (\dim\Hom_{A^\env}(C_0,A)-\dim\ker(\phi))
\\=&\left(\sum_{\alpha\in Q_1}\nu(\alpha)\right) - n + \dim\HH^0(A).
\end{align*}
If $A$ has global dimension $\leq 2$, $P_2=0$ and equality holds.
\end{proof}
We can use this in the case of canonical algebras for $2$ points.
\begin{proposition}\label{prop:cohomologyCanonical2}
Let $A$ be the canonical algebra associated to the points $\mathbf x=(\infty,0)$ and the weights $\mathbf a = (a_1,a_2)$.
The Hochschild (co)homology of $A$ in degree $i\in\mathbb Z$ is given by
\begin{align*}
&\dim\HH_i(A) = \begin{cases}a_1+a_2,&i=0\\
0,~&\text{otherwise}\end{cases}
\\&\dim\HH^0(A)=1
\\&\dim\HH^1(A)=\begin{cases}3,&a_1=a_2=1\\
2,&a_1=1\neq a_2\text{ or }a_1\neq 1=a_2\\
1,&a_1\neq 1\neq a_2\end{cases}.
\\&\dim\HH^i(A)=0\text{ for }i\notin\{0,1\}.
\end{align*}
\end{proposition}
\begin{proof}
Recall that with the given assumptions, $A$ is the path algebra over the quiver
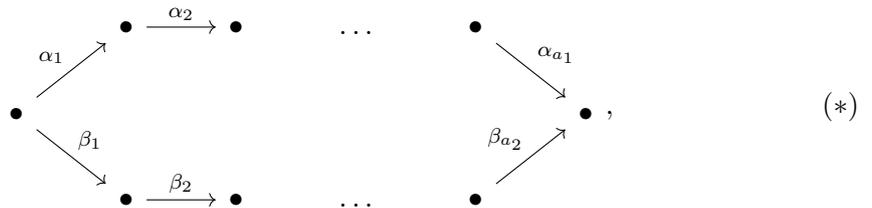
\begin{align*}
\begin{tikzcd}[ampersand replacement=\&]
\&\bullet\ar[r,"{\alpha_2}"]\&\bullet\&\hdots\&\bullet\ar[dr,"{\alpha_{a_1}}"]\\
\bullet\ar[ur,"{\alpha_1}"]\ar[dr,"{\beta_1}"]\&\&\&\&\&\bullet\\
\&\bullet\ar[r,"{\beta_2}"]\&\bullet\&\hdots\&\bullet\ar[ur,"{\beta_{a_2}}"]
\end{tikzcd},\tag{$\ast$}
\end{align*}
where the top string contains $a_1$ arrows and the bottom one $a_2$ arrows.
In particular, $A$ has global dimension equal to $1$.

The claim on $\HH_\bullet(A)$ follows immediately from Proposition~\ref{prop:algebraHomology} and simple counting.
The claim on $\HH^0(A)$ follows immediately from Corollary~\ref{cor:hochschildCohomologyHereditary}.
Vanishing of $\HH^{\geq 2}(A)$ follows from Corollary~\ref{cor:cohomologyVanishesAboveGlDim}.

By Corollary~\ref{cor:hochschildCohomologyHereditary}, we have
\begin{align*}
\dim \HH^1(A) = 1-(a_1+a_2) + \sum_{i=1}^{a_1} \nu(\alpha_i) + \sum_{i=1}^{a_1} \nu(\beta_i).
\end{align*}
If $a_i>1$, we have $\nu(\alpha_i)=1$ for $i=1,\dotsc,a_1$, so that $-a_1 + \sum\limits_{i=1}^{a_i}\nu(\alpha_i)=0$.
If $a_i=1$, we have $\nu(\alpha_1)=2$ so that $-a_1 + \sum\limits_{i=1}^{a_i}\nu(\alpha_i)=1$.

The same argument shows
\begin{align*}
-a_2 + \sum_{i=1}^{a_2} \nu(\beta_i)=\begin{cases}1,&a_2=1\\
0,&a_2>1.\end{cases}
\end{align*}
Then the claim on $\HH^1(A)$ follows.
\end{proof}
This method works only for hereditary algebras. The canonical algebra associated
to a weighted projective line with at least three exceptional points will not be hereditary,
so this method needs a refinement in order to deal with these cases.
We present two solutions for this.
\section{One-point extensions}
Throughout this section, let $A$ be a finite-dimensional algebra over $k$.
\begin{definition}{\cite[Section~5.1]{Happel1989}}
Let $M$ be a left $A$-module. The \emph{one-point extension} of $A$ by $M$ is
\begin{align*}
A[M] = \begin{pmatrix}A&M\\0&k\end{pmatrix},
\end{align*}
i.e.\ elements are matrices $\begin{pmatrix}a&m\\0&\lambda\end{pmatrix}$
with $a\in A,m\in M, \lambda\in k$ and multiplication comes from the usual rules for matrix multiplication.
\end{definition}
\begin{theorem}[{\cite[Theorem~5.3]{Happel1989}}]\label{thm:happelLES}If $M$ is a left $A$-module, there exists a long exact sequence
\begin{align*}
&0\rightarrow \HH^0(A[M])\rightarrow \HH^0(A)\rightarrow \Hom_A(M,M)/k\rightarrow \HH^1(A[M])\rightarrow \HH^1(A)\rightarrow\Ext^1_A(M,M)\rightarrow
\\&\qquad\cdots\rightarrow \Ext^i_A(M,M)\rightarrow \HH^{i+1}(A[M])\rightarrow\HH^{i+1}(A)\rightarrow\Ext^{i+1}_A(M,M)\rightarrow\cdots.\rightqed
\end{align*}
\end{theorem}
Note that this theorem is generalized in \cite[Lemma~4.5]{Keller2003}.

In \cite{Happel1998}, this is applied to canonical algebras.
Consider a vector of points $\mathbf x=(x_1,\dotsc,x_n)$ and weights $\mathbf a=(a_1,\dotsc,a_n)$ defining a canonical algebra.
\begin{definition}\label{def:concreteOnePointExtension}
We denote by $Q$ the quiver from Proposition~\ref{proposition-canonical}, and label it as
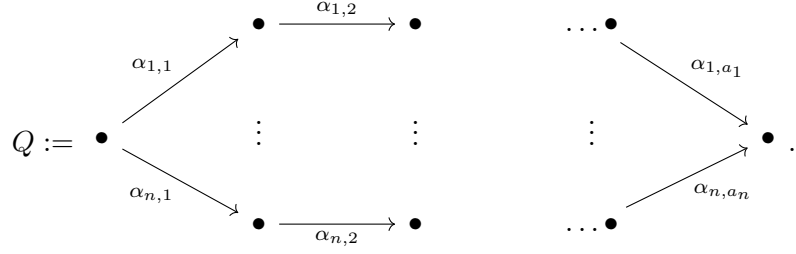
\begin{align*}
Q:=\begin{tikzcd}[ampersand replacement=\&,column sep=4em]
\&\bullet\ar[r,"{\alpha_{1,2}}"]\&\bullet\&\hdots\bullet\ar[dr,"{\alpha_{1,a_1}}"]\\
\bullet\ar[ru,"{\alpha_{1,1}}"]\ar[rd,"{\alpha_{n,1}}"']\&\vdots\&\vdots\&\vdots\&\bullet\\
\&\bullet\ar[r,"{\alpha_{n,2}}"']\&\bullet\&\hdots\bullet\ar[ur,"{\alpha_{n,a_n}}"']\\
\end{tikzcd}.
\end{align*}
Let $I\subseteq kQ$ be the ideal such that $C:=kQ/I$ is the canonical algebra from Proposition~\ref{proposition-canonical}.
Denote by $c$ the unique sink in $Q$ and by $s$ the unique source.

Let $Q'\subseteq Q$ be the subquiver by removing the arrows $\alpha_{1,a_1},\dotsc,\alpha_{n,a_n}$
and the vertex $c$.

Let $A:= kQ'$ and $M\subseteq C$ be the $A$-module
generated by all paths in $Q$ which end in $c$ and have positive length.
\end{definition}
\begin{lemma}
$C\cong A[M]$.
\end{lemma}
\begin{proof}
This is a simple observation.
\end{proof}
\begin{lemma}\label{lem:hochschildCohomologyAuxillary}
$\HH^0(A)\cong k$ and $\HH^i(A)=0$ for $i\neq 0$.
\end{lemma}
\begin{proof}
This follows again from Corollaries~\ref{cor:cohomologyVanishesAboveGlDim} ($A$ has global dimension one)
and \ref{cor:hochschildCohomologyHereditary} (all arrows $\alpha$ in $Q'$ have $\nu(\alpha)=1$).
\end{proof}
\begin{lemma}\label{lem:Mresolution}
For $i=1,\dotsc,n$, let $v_i:=s(\alpha_{i,a_i})$ be the $i$-th sink of $Q'$.
\begin{enumerate}[(a)]
\item
$P:=\bigoplus\limits_{i=1}^n Av_i$ is a projective $A$-submodule of $A$.
\item
There is an epimorphism
\begin{align*}
P\rightarrow M, x \mapsto x(\alpha_{1,a_1}+\cdots + \alpha_{n,a_n}).
\end{align*}
\item The kernel of this epimorphism $K$ is generated by the elements
\begin{align*}
k_i := \alpha_{i,1}\cdots \alpha_{i,a_i-1} - \alpha_{2,1}\cdots \alpha_{2,a_2-1} + x_i \alpha_{1,1}\cdots \alpha_{1,a_1-1}
\end{align*}
for $i=3,\dotsc,n$.
\end{enumerate}
\end{lemma}
\begin{proof}
\begin{enumerate}[(a)]
\item
Recall that $A$ is, as $A$-module, $\bigoplus\limits_{v\text{ vertex in }Q'} (Av)$.
\item The image of $v_i$ is $\alpha_{i,a_i}$, and these elements clearly generate $M$ as $A$-module.
\item These generators correspond to the generators of the ideal in Proposition~\ref{proposition-canonical}.
\qedhere\end{enumerate}
\end{proof}
\begin{lemma}\label{lem:homExtAuxillary}
Suppose $n\geq 3$. Then
$\dim\Hom_A(M,M)=1$ and $\dim\Ext^1_{A}(M,M)=n-3$.
\end{lemma}
\begin{proof}
We apply $\Hom_A(\blank,M)$ to Lemma~\ref{lem:Mresolution} and get a short exact sequence
\begin{align*}
0\rightarrow\Hom_A(M,M)\rightarrow \Hom_A(P,M)\rightarrow \Hom_A(K,M)\rightarrow \Ext^1_A(M,M)\rightarrow 0
\end{align*}
since $P$ is projective.

As $P=\bigoplus\limits_{i=1}^n Av_i$, we have $\Hom_A(P,M)\cong \bigoplus\limits_{i=1}^n v_i M = \bigoplus\limits_{i=1}^n k \alpha_{i,a_i}$.
We can identify this as a subspace of $M$.

$K$ is a projective $A$-module (as $A$ is hereditary). It is isomorphic to the projective module $(As)^{n-2}$ (with $s$ being the vertex from
Definition~\ref{def:concreteOnePointExtension}).
Hence $\Hom_A(K,M)\cong (sM)^{n-2}$. Note that $sM$ corresponds to all paths in $Q$ from $s$ to $c$
modulo the ideal defining the canonical algebra, so $sM$ two-dimensional with basis
\begin{align*}
\alpha_{1,1}\cdots \alpha_{1,a_1},\quad \alpha_{2,1}\cdots\alpha_{2,a_2}.
\end{align*}

With these identifications, the map $\Hom_A(P,M)\rightarrow \Hom_A(K,M)$ corresponds to
\begin{align*}
\bigoplus_{i=1}^n k\alpha_{i,a_i}\rightarrow (sM)^{n-2},
x\mapsto (k_3 x,\dotsc,k_n x)
\end{align*}
with $k_3,\dotsc,k_n$ as in Lemma~\ref{lem:Mresolution}.
We claim that the rank of this map is $n-1$:

One the one hand, the rank is at least $n-1$, since the elements $\alpha_{2,a_2},\dotsc,\alpha_{n,a_n}$ are sent
to linearly independent elements in $(sM)^{n-2}$ (if $i\geq 3$, $\alpha_{i,a_i}$ is sent
to $(0,\dotsc,0,\alpha_{i,1}\cdots\alpha_{i,a_i},0,\dotsc,0)$ with only the $i$-th entry non-zero).

On the other hand, the kernel is non-trivial, since $\alpha_{1,a_1}+\cdots+\alpha_{n,a_n}$ is sent to zero.

Since $\dim\bigoplus\limits_{i=1}^n k\alpha_{i,a_i}=n$, the map must have rank equal to $n-1$.

Hence $\Hom_A(M,M)$, identified with the kernel, is one-dimensional. $\Ext^1_A(M,M)$, identified as the cokernel,
is has dimension $n-3$.
\end{proof}
\begin{proposition}[{\cite[Theorem~2.4]{Happel1998}}]\label{prop:cohomologyCanonical3}
Suppose $n\geq 3$. The Hochschild cohomology of the canonical algebra $C=A[M]$ in degree $i\in\mathbb Z$ is given by
\begin{align*}
\dim \HH^i(C)=\begin{cases}1,&i=0\\
0,&i=1\\
n-3,&i=2\\
0,&\text{otherwise}.\end{cases}
\end{align*}
\end{proposition}
\begin{proof}
We use the long exact sequence of Theorem~\ref{thm:happelLES}.
Since $\dim\Hom_A(M,M)=1$ by Lemma~\ref{lem:homExtAuxillary}, exactness of
\begin{align*}
0\rightarrow \HH^0(C)\rightarrow \HH^0(A)\rightarrow \Hom_A(M,M)/k
\end{align*}
yields $\dim\HH^0(C)=\dim\HH^0(A)=1$ by Lemma~\ref{lem:hochschildCohomologyAuxillary}.

Moreover, exactness of
\begin{align*}
0=\Hom_A(M,M)/k\rightarrow \HH^1(C)\rightarrow \HH^1(A)
\end{align*}
yields $\HH^1(C)=0$ since $\HH^1(A)=0$ by Lemma~\ref{lem:hochschildCohomologyAuxillary}.

Exactness of
\begin{align*}
0=\HH^1(A)\rightarrow \Ext^1_A(M,M)\rightarrow \HH^2(C)\rightarrow \HH^2(A)=0
\end{align*}
shows, using Lemma~\ref{lem:homExtAuxillary}, that $\dim\HH^2(C)=n-3$.

If $i\geq 3$, exactness of
\begin{align*}
\Ext^{i-1}_A(M,M)\rightarrow \HH^i(C)\rightarrow \HH^i(A)=0
\end{align*}
shows $\HH^i(C)=0$ since $A$ has global dimension $1<i-1$.
\end{proof}
We conclude:
\begin{corollary}\label{cor:happelWPL}
The Hochschild (co)homology of a weighted projective line $X$ with $n$ exceptional points
of weight $a_1,\dotsc,a_n\in\mathbb Z_{\geq 2}$
in degree $i\in\mathbb Z$ is given by
\begin{align*}
&\dim\HH_i(X)=\begin{cases}2+\sum\limits_{i=1}^n a_i-1,&i=0\\
0,&i\neq 0\end{cases}
\\&\dim\HH^i(X)=\begin{cases}1,&i=0\\
\max(3-n,0),&i=1\\
\max(n-3,0),&i=2\\
0,~&\text{otherwise}.\end{cases}
\end{align*}
\end{corollary}
\begin{proof}
We know that the Hochschild (co)homology groups of $X$ and the associated
canonical algebra coincide by Corollary~\ref{cor:wplCanonicalHochschild}.

The Hochschild homology is due to Proposition~\ref{prop:algebraHomology} and some simple counting.

Hochschild cohomology for $n\leq 2$ is Proposition~\ref{prop:cohomologyCanonical2}.
For $n\geq 3$, it is Proposition~\ref{prop:cohomologyCanonical3}.
\end{proof}
\section{Alternative proof}
Happel's method using one-point extensions gives a quick calculation of the Hochschild cohomology
of a canonical algebra in case $n\geq 3$. However, it does not make clear \emph{why} the Hochschild
cohomology depends only on the number of exceptional points, and not the weights $a_i$ nor the coordinates $x_i$.

We want to shed some light on this question by presenting a different method to compute the Hochschild cohomology.

Let $A=kQ/I$ be as in section~\ref{sec:minimalResolution}. Suppose for this section that $Q$ is acyclic.

Recall that the differentials in Proposition~\ref{prop:happel89} are implicit. We were able to make
$P_1\xrightarrow{d_1} P_0\xrightarrow{\varepsilon} A$ explicit in section~\ref{sec:minimalResolution}.
In this section, we also make $P_2\xrightarrow{d_2} P_1$ explicit, allowing to directly
write down a complex computing the Hochschild cohomology of an algebra of global dimension $\leq 2$.

Let $e_1,\dotsc,e_n$ be again the vertices of $Q$, $S(i)$ and $P(i)$ the corresponding
simple resp.\ projective modules.
\begin{lemma}\label{lem:minimalRelations}~
\begin{enumerate}[(a)]
\item There exist elements $r_1,\dotsc,r_m\in I$ such that
\begin{itemize}
\item For each $r_i$, there exist uniquely determined vertices $s(r_i),e(r_i)\in Q_0$, calles \emph{start} resp.\ \emph{end} of $r_i$,
such that $r_i \in s(r_i) (kQ) e(r_i)$. i.e.\ $r_i$ is a linear combination of paths starting in $s(r_i)$ and ending in $e(r_i)$.
\item $r_1,\dotsc,r_m$ is a minimal generating set of $I$. i.e.\ $\{r_1,\dotsc,r_m\}$ generates $I$ as ideal, but no proper subset does.
\end{itemize}
\item For any such generating system $r_1,\dotsc,r_m$ and for  $i,j\in\{1,\dotsc,n\}$, $\dim\Ext^2_A(S(i),S(j))$ is the number of elements $r_\ell$
such that $s(r_\ell)=e_j, e(r_\ell)=e_i$.
\end{enumerate}
\end{lemma}
\begin{proof}
\begin{enumerate}[(a)]
\item
By \cite[Corollary~II.2.9]{Assem2006}, there exists a generating set $r_1,\dotsc,r_m\in I$
such that each $r_i$ is in $s(r_i)(kQ)e(r_i)$ for uniquely determined vertices $s(r_i),e(r_i)\in Q$.
Replacing $\{r_1,\dotsc,r_m\}$ by a minimal generating subset, we may assume
that $\{r_1,\dotsc,r_m\}$ is a minimal generating set of $I$.
\item
Recall that $\rad(kQ)$ is the ideal of $kQ$ generated by all arrows in $Q$.
The fact that $\{r_1,\dotsc,r_m\}$ is a minimal generating system of $I$ is equivalent
to $r_1+\rad(I),\dotsc,r_m+\rad(I)$ being a basis of $I/\rad(I)$
(generating $I/\rad(I)$ is as vector space is equivalent to generating $I$ by Nakayama's lemma,
and a minimal generating system of a vector space is a basis).

We would like to calculate $\Ext^2_A$ of simple modules.
Let $i\in\{1,\dotsc,n\}$ and denote the arrows ending in $e_i$ by $\pi := \{\alpha \in Q_1 : e(\alpha)=e_i\}$.
Then there exists a minimal projective presentation of $S(i)$ of the form
\begin{align*}
\bigoplus_{\alpha\in \pi} As(\alpha)\rightarrow P(i)\rightarrow S(i)\rightarrow 0.
\end{align*}
Let us continue this to a minimal projective resolution
\begin{align*}
\cdots\rightarrow P_3\rightarrow P_2\rightarrow 
\bigoplus_{\alpha\in \pi} As(\alpha)\rightarrow P(i)\rightarrow \rightarrow 0.
\end{align*}
Now $\Ext^\bullet_A(S(i),S(j))$ can be computed by applying $\Hom_A(\blank,S(j))$ to this resolution and taking cohomology.

For each differential $P_s\xrightarrow{d_s} P_{s-1}$ of this resolution, the resolution being
minimal implies that $\ker(d_s)\subseteq \rad(P_s)$ and $\im(d_s)\subseteq \rad(P_{s-1})$.
Therefore, $\Hom_A(d_s,S(j))=0$ for all $s$.

We conclude
\begin{align*}
&\Ext^2_A(S(i),S(j))\cong \Hom_A(P_2,S(j))\cong \Hom_A(P_2/\rad(P_2),S(j))
\\&\underset{\ker(d_2)\subseteq \rad(P_2)}\cong \Hom_A(\im(d_2)/\rad(\im(d_2)),S(j))
\\&\cong \Hom_A(\im(d_2),S(j)) \underset{\text{exactness}}\cong \Hom_A(\ker(d_1),S(j)).
\end{align*}
Put $K := \ker(d_1) = \ker\Bigl(\bigoplus\limits_{\alpha\in \pi} As(\alpha)\rightarrow P(i)\Bigr)$.

Note that for each projective module $P(i) = Ae_i$, there is a short exact sequence of $kQ$-modules
$0\rightarrow Ie_i\rightarrow (kQ)e_i\rightarrow Ae_i\rightarrow 0$.
We therefore get a commutative diagram:
\begin{align*}
\begin{tikzcd}[ampersand replacement=\&,row sep=1em]
\&\&0\ar[d]\&0\ar[d]\&0\ar[d]\\
R_1:\&0\ar[r]\&\bigoplus_{\alpha\in \pi}Is(\alpha)\ar[d]\ar[r]\&Ie_i\ar[r]\ar[d]\&0\ar[r]\ar[d]\&0\\
R_2:\&0\ar[r]\&\bigoplus_{\alpha\in \pi}(kQ)s(\alpha)\ar[d]\ar[r]\&(kQ)e_i\ar[r]\ar[d]\&S(i)\ar[d,equals]\ar[r]\&0\\
R_3:\&0\ar[r]\&\bigoplus_{\alpha\in \pi}As(\alpha)\ar[d]\ar[r]\&P(i)\ar[r]\ar[d]\&S(i)\ar[d]\ar[r]\&0\\
\&\&0\&0\&0
\end{tikzcd}
\end{align*}
Here, the map $\bigoplus\limits_{\alpha\in \pi}(kQ)s(\alpha)\rightarrow (kQ)e_i$ sends an element in $x\in (kQ)s(\alpha)$ to $x\alpha\in (kQ)e_i$.
Note that all columns and the row $R_2$ are exact. This yields a short exact sequence of chain complexes
$0\rightarrow R_1\rightarrow R_2\rightarrow R_3\rightarrow 0$. The long exact homology sequence
together with exactness of $R_2$ shows $\H_\bullet(R_3)=\H_{\bullet+1}(R_1)$.
Hence 
\begin{align*}
K=\H_0(R_3) \cong \H_1(R_1)\cong \frac{Ie_i}{\bigoplus_{\alpha\in \pi} I\alpha}=\frac{Ie_i}{\rad(Ie_i)},
\end{align*}
where $\rad(Ie_i)$ denotes the radical of the left $kQ$-module $Ie_i$.
Now
\begin{align*}
&\Hom_A((Ie_i)/\rad(Ie_i),S(j)) \cong \Hom_{A/\rad(A)} ((Ie_i)/\rad(Ie_i),S(j))
\\&\cong e_j(Ie_i)/\rad(Ie_i)\cong e_j(I/\rad(I))e_i.
\end{align*}
Since $r_1+\rad(I),\dotsc,r_m+\rad(I)$ is a basis of $I/\rad(I)$, the dimansion of $e_j(I/\rad(I))e_i$
is equal to the number of relations $r_\ell$ with $s(r_\ell)=e_j, e(r_\ell)=e_i$.
\qedhere\end{enumerate}
\end{proof}
\begin{corollary}
With relations $r_1,\dotsc,r_m$ as in part (a) of Lemma~\ref{lem:minimalRelations}
and
\begin{align*}
\cdots\rightarrow P_m\xrightarrow{d_m} P_{m-1}\rightarrow\cdots\rightarrow P_1\xrightarrow{d_1} P_0\xrightarrow{\varepsilon} A\rightarrow 0
\end{align*}
being the minimal projective resolution from Proposition~\ref{prop:happel89},
we have
\begin{align*}
P_2\cong \bigoplus_{i=1}^m As(r_i)\otimes_k e(r_i)A.
\end{align*}
\end{corollary}
\begin{proof}
This follows from Proposition~\ref{prop:happel89} and the observation
\begin{align*}
\forall i,j\in \{1,\dotsc,n\}: (Ae_j\otimes_k e_i A)^{\dim\Ext^2_A(S(i),S(j))}\cong\bigoplus_{\substack{\ell=1\\s(r_\ell)=e_j\\e(r_\ell)=e_i}}^m
As(r_\ell)\otimes_k e(r_\ell)A.
\end{align*}
which follows from Lemma~\ref{lem:minimalRelations} part (b).
\end{proof}
\noindent Throughout this section, fix relations $r_1,\dotsc,r_m$ as in Lemma~\ref{lem:minProjective1} part (a).
\begin{example}\label{ex:relationsCanonical}
We are mainly interested in the case of canonical algebras.
In the setting of Proposition~\ref{proposition-canonical}, the $n-2$
explicit generators $r_1,\dotsc,r_{n-2}$ where
\begin{align*}
r_i = X_i^{a_i} - X_2^{a_2} + x_i X_1^{a_1}
\end{align*}
satisfy the conditions of lemma~\ref{lem:minProjective1} part (a).
\end{example}
\begin{definition}
Let $\partial$ be the map from Lemma~\ref{lem:minProjective1}.

We define
\begin{align*}
&C_2 := \bigoplus_{\ell=1}^m As(r_\ell)\otimes_k e(r_\ell) A
\\&C_2\xrightarrow{c_2} C_1: a \otimes b \mapsto a \partial(r_\ell) b
\end{align*}
for $\ell\in \{1,\dotsc,m\}$, $a\in As(r_\ell)$ and $b\in e(r_\ell)A$.
\end{definition}
\begin{lemma}\label{lem:paritalIinimc_2}
With $\partial,c_2$ as above, we have $\partial(I)\subseteq \im(c_2)$.
\end{lemma}
\begin{proof}
Since $\partial$ is $k$-linear and $I$ is generated, as $k$-vector space,
by elements of the form
\begin{align*}
a r_\ell b
\end{align*}
with $a,b$ paths in $Q$ such that $e(a)=s(r_\ell), e(r_\ell)=s(b)$, it suffices
to show that $\partial(ar_\ell b)\in\im(c_2)$ for all these elements.

By Lemma~\ref{lem:minProjective1} part (a),
\begin{align*}
&\partial(ar_\ell b)=\partial(a)r_\ell b + a\partial(r_\ell)b + ar_\ell\partial(b)
\underset{r_\ell=0\text{ in }A}=a\partial(r_\ell)b = c_2(a\otimes b).\qedhere
\end{align*}
\end{proof}
\begin{lemma}\label{lem:minProjective2}
The sequence
\begin{align*}
C_2\xrightarrow{c_2}C_1\xrightarrow{c_1}C_0
\end{align*}
is exact.
\end{lemma}
\begin{proof}
We first show that $c_1 c_2=0$:
Let $\ell\in\{1,\dotsc,m\}$. Then
\begin{align*}
c_1c_2(s(r_\ell)\otimes e(r_\ell)) =&c_1(\partial(r_\ell))\underset{\text{L\ref{lem:minProjective1}}}=
\incl_{s(r_\ell)}(s(r_\ell)\otimes r_\ell) - 
\incl_{e(r_\ell)}(r_\ell\otimes e(r_\ell))
\\\underset{r_\ell=0\text{ in }A}=&0.
\end{align*}
Now we have to show that the kernel of $c_1$ is contained in the image of $c_2$.

Let $C_0\xrightarrow{s_0}C_1$ be the homotopy from Lemma~\ref{lem:minProjective1}.
We have $\ker(c_1)\subseteq \im(\id_{C_1} - s_0c_1)$, because if $c_1(x)=0$, $x=(\id_{C_1}-s_0c_1)(x)$.

It hence suffices to show that the image of $\id_{C_1} - s_0c_1$ is contained in the image of $c_2$.
Note that $\id_{C_1} - s_0c_1$ and $c_2$ are morphisms of left $A$-modules, and that
$C_1$ is generated, as left $A$-module, by elements of the form
\begin{align*}
\incl_\alpha(1\otimes p)
\end{align*}
for $\alpha\in Q_1$ and $p$ a path in $Q$ starting at $e(\alpha)$. So it is enough to show that 
\begin{align*}
(\id_{C_1}-s_0c_1)(\incl_\alpha(1\otimes p))\in \im(c_2).
\end{align*}

We evaluate $s_0c_1$ as follows:
\begin{align*}
s_0c_1(\incl_\alpha(1\otimes p)) = s_0\left(\incl_{s(\alpha)}(\alpha p) - \incl_{e(\alpha)}(\alpha\otimes p)\right)
=\partial(\varphi(\alpha p)) - \alpha \partial(\varphi(p)).
\end{align*}
Note that $\varphi(\alpha p)-\alpha\varphi(p)\in kQ$ is an element in the kernel of the projection $\pi:kQ\rightarrow A$
(since $\pi$ is an algebra homomorphism and $\pi\varphi=\id_A$). Hence $\varphi(\alpha p)-\alpha\varphi(p)\in I$.
By Lemma~\ref{lem:paritalIinimc_2}, $\partial(\varphi(\alpha p)-\alpha\varphi(p))\in \im(c_2)$.

By Lemma~\ref{lem:minProjective1} part (a),
\begin{align*}
\partial(\varphi(\alpha p)-\alpha\varphi(p))
=\partial(\varphi(\alpha p)) - \partial(\alpha)\varphi(p) - \alpha \partial(\varphi(p))
\end{align*}
By definition of $\partial$ and the observation that the image of $\varphi(p)$ in $A$ is $p$,
we conclude \begin{align*}\partial(\alpha)\varphi(p)=\incl_\alpha(1\otimes 1)p = \incl_\alpha(1\otimes p).\end{align*}

Therefore,
\begin{align*}
(\id-s_0c_1)(\incl_\alpha(1\otimes p)) 
=&\incl_\alpha(1\otimes p) - \left(\partial(\varphi(\alpha p)) - \alpha\partial(\varphi(p))\right)
\\=&-\partial(\varphi(\alpha p)-\alpha\varphi(p)) \in \partial(I)\subseteq \im(c_2).
\end{align*}
This finishes the proof.
\end{proof}
\begin{corollary}\label{cor:cohomologyGlDim2}
Suppose that the global dimension of $A$ is $\leq 2$.
$\HH^\bullet(A)$ is isomorphic to the cohomology of the cochain complex
\begin{align*}
0\rightarrow k[Q_0]\rightarrow \bigoplus_{\alpha\in Q_1} s(\alpha)Ae(\alpha)\rightarrow \bigoplus_{\ell=1}^m s(r_\ell)Ae(r_\ell)\rightarrow 0
\end{align*}
with the two maps defined as follows:
\begin{enumerate}[(a)]
\item For $v\in Q_0$ and $\alpha\in Q_1$,
the image of $v$ in $s(\alpha)A e(\alpha)$ is
\begin{align*}
\begin{cases}\alpha,&s(\alpha)=v\\
-\alpha,&e(\alpha)=v\\
0,~&\text{otherwise}.\end{cases}
\end{align*}
\item For $\alpha\in Q_1$ and $\ell\in\{1,\dotsc,m\}$,
the map $s(\alpha)Ae(\alpha)\rightarrow s(r_\ell)Ae(r_\ell)$ is the composition
\begin{align*}
&s(\alpha)Ae(\alpha)\xrightarrow{\blank\otimes\partial(r_\ell)}
s(\alpha)Ae(\alpha)\otimes\bigoplus_{\beta\in Q_1} (s(r_\ell)As(\beta)\otimes e(\beta)Ae(r_\ell))
\\&\xrightarrow{\id\otimes\text{projection}}
s(\alpha)Ae(\alpha)\otimes (s(r_\ell)As(\alpha)\otimes e(\alpha) Ae(r_\ell))
\xrightarrow{(a\otimes b\otimes c)\mapsto bac} s(r_\ell)Ae(r_\ell).
\end{align*}
\end{enumerate}
\end{corollary}
\begin{proof}
Since $A$ has global dimension $\leq 2$,
the minimal projective resolution of $A$ as $A^\env$-module from Proposition~\ref{prop:happel89} is, by Lemmas~\ref{lem:minProjective1} and \ref{lem:minProjective2},
isomorphic to
\begin{align*}
0\rightarrow C_2\xrightarrow{c_2}C_1\xrightarrow{c_1}C_0\rightarrow 0.
\end{align*}
Applying $\Hom_{A^\env}(\blank,A)$ gives, with little effort and Lemma~\ref{lem:trivialFactsFDAlgebras}, the desired cochain complex.
\end{proof}
\begin{remark}\label{rem:cohomologyGlDim2Rewrite}
The morphism in part (b) of Corollary~\ref{cor:cohomologyGlDim2}
can be made more explicit as follows:

Pick a morphism $\alpha\in Q_1$ and a relation $r_\ell$. Write
\begin{align*}
r_\ell = \sum_{i=1}^N \lambda_i p_i
\end{align*}
for scalars $\lambda_i$ and paths $p_i$ in $Q$. We may order the paths such that
$\alpha$ is contained in $p_1,\dotsc,p_M$ and not contained in $p_{M+1},\dotsc,p_N$.

For $i\leq M$, $\alpha$ occurs precisely once in $p_i$ (as $Q$ is acyclic),
so we may write $p_i = q_i \alpha q_i'$ for uniquely determined paths $q_i,q_i'$.

The projection of $\partial(r_\ell)$ to the $\alpha$-summand of $C_1$ equals
\begin{align*}
\sum_{i=1}^M \lambda_iq_i\otimes q_i'.
\end{align*}
Hence the induced map $s(\alpha) A e(\alpha) \rightarrow s(r_\ell) A e(r_\ell)$
sends $x\in s(\alpha) A e(\alpha)$ to
\begin{align*}
\sum_{i=1}^M \lambda_iq_i x q_i'.
\end{align*}
In particular, this map sends $\alpha\in s(\alpha) A e(\alpha)$
to $\sum\limits_{i=1}^M \lambda_i p_i\in s(r_\ell) A e(r_\ell)$.
\end{remark}
\begin{example}
Consider a weighted projective line with $n\geq 3$
exceptional points. Let $a_1,\dotsc,a_n\in\mathbb Z_{\geq 2}$ be the weights
and $x_1=\infty,x_2=0,x_3=1,x_4,\dotsc,x_n\in \mathbb P^1$ be the coordinates of the exceptional points.

Let $A=kQ/I$ be the canonical algebra associated to this weighted projective line
as in Proposition~\ref{proposition-canonical}. We will use the same labels for the arrows and vertices
as this proposition. We would like to compute the Hochschild cohomology using Corollary~\ref{cor:cohomologyGlDim2}.

Let, as in Example~\ref{ex:relationsCanonical}, $r_3,\dotsc,r_{n}$ be the relations
\begin{align*}
r_i = X_i^{a_i}-X_2^{a_2} + x_i X_1^{a_1}.
\end{align*}
Note that for each morphism $\alpha$ in $Q_1$, $\alpha$
is the only path in $Q$ starting at $s(\alpha)$ and ending in $e(\alpha)$
(since $a_i\geq 2$ for all $i$). Thus $s(\alpha)Ae(\alpha)=k\alpha$.

For each relation $r_\ell$, $s(r_\ell)=0$ and $e(r_\ell)=c$.
Therefore $s(r_\ell)Ae(r_\ell)$ is two-dimensional
with a basis given by $\rho_1:=X_1^{a_1}, \rho_2:=X_2^{a_2}$.
Hence a basis of $\bigoplus\limits_{\ell=3}^{n} s(r_\ell) A e(r_\ell)$
is given by elements of the form
\begin{align*}
\rho_1^3,\rho_2^3,\rho_1^4,\rho_2^4\dotsc,\rho_1^{n},\rho_2^{n}.
\end{align*}

We identify the cochain complex of Corollary~\ref{cor:cohomologyGlDim2} as
\begin{align*}
0\rightarrow k[Q_0]\xrightarrow fk[Q_1]\xrightarrow g\bigoplus_{i=3}^{n} s(r_\ell)A e(r_\ell)\rightarrow 0.
\end{align*}
Here, $f$ sends a vertex $v$ to
\begin{align*}
\sum_{\substack{\alpha\in Q_1\\s(\alpha)=v}} \alpha - 
\sum_{\substack{\alpha\in Q_1\\e(\alpha)=v}}\alpha.
\end{align*}
$g$ sends (using Remark~\ref{rem:cohomologyGlDim2Rewrite}) an arrow labelled $X_1$
to
\begin{align*}
x_3\rho_1^3+x_4\rho_1^4+\cdots+x_n\rho_1^{n}.
\end{align*}
Arrows with the label $X_2$ are sent to
\begin{align*}
-\rho_2^3-\rho_2^4-\cdots-\rho_2^{n}.
\end{align*}
Finally, an arrow labelled $X_i$ for $i\geq 3$ is mapped to
\begin{align*}
-\rho_2^i + x_i \rho_1^i.
\end{align*}
It is straightforward to see that $g$ has rank $n-1$.
By Corollary~\ref{cor:hochschildCohomologyHereditary}, the rank of $f$ is $\#Q_0-1$.
It follows that
\begin{align*}
&\dim \HH^1(A)=\dim k[Q_1] - \rk(f)-\rk(g) = \#Q_1 - \#Q_0 + n-2 = 0,
\\&\dim\HH^2(A) = \dim\left(\bigoplus_{\ell=3}^n s(r_\ell)Ae(r_\ell)\right)-\rk(g) = 2(n-2) - (n-1) = n-3.
\end{align*}
We get the same result as in Proposition~\ref{prop:cohomologyCanonical3}.
\end{example}
\begin{remark}
\begin{enumerate}[(a)]
\item We see again that the number of exceptional points alone
determines the Hochschild cohomology, independent of the weights and coordinates.
We now may give a first answer, why this is the case:

Indeed, the weights $a_i$ \enquote{cancel out} in the calculation of $\HH^1$ (since $\#Q_0-\#Q_1 = n-2$ holds independently
of the weights) and the weights play no role when computing $\HH^2$.
The coordinates $x_i$ are irrelevant since we are only interested in the rank of $g$, which happens
to be $n-1$ regardless of the coordintes $x_i$.
\item
When working with finite-dimensional algebras, it is also possible to use computer algebra methods
to calculate specific examples. We use the package \cite{QPA} of the computer algebra system \cite{GAP}.
Our routines to compute Hochschild cohomology and to create canonical algebras are printed and explained in Appendix~\ref{chap:gap_routines}.
\end{enumerate}
\end{remark}.
\nocite{stacks-project}
\chapter{Quotient stacks}
We work over an algebraically closed field $k$ of characteristic zero, which we usually will omit from our notation.

In this chapter, we develop a method which can be used to compute Hochschild (co)homology groups
of one-dimensional iterated root stacks.
\section{The decomposition of Arinkin-C\u{a}ld\u{a}raru-Hablicsek}
Let $X$ be a smooth projective $k$-scheme and $G\leq \Aut(X)$ a finite group of $k$-scheme-au\-to\-mor\-phisms.
Similar to the HKR-decomposition (Theorem~\ref{thm:HKR}), 
\cite{Arinkin2014} gives a decomposition of the Hochschild (co)homology groups
of the quotient stack $[X/G]$.

We use the following notation, where $g\in G$:
\begin{itemize}
\item $X^g$ is the fixed locus of $g$, i.e.\ the equalizer of
$\begin{tikzcd}X\ar[r,"\id",shift left=0.5ex]\ar[r,"g"',shift right=0.5ex]&X\end{tikzcd}$ in the category of $k$-schemes.
\item $c_g$ is the codimension of $X^g$ in $X$.
\item $N_{X^g/X}$ is the normal bundle of $X^g$ in $X$ and $\omega_g := \bigwedge^{c_g} N_{X^g/X}$ its highest exterior power.
\item $\T_g$ is the vector bundle on $X^g$ obtained by restricting $\T_X$ to $X^g$ and taking coinvariants
with respect to the action of $g$.
\item $\Omega_g^j$ the dual bundle of $\bigwedge^j \T_g$ (as bundles on $X^g$).
\end{itemize}
\begin{theorem}[Arinkin-C\u{a}ld\u{a}raru-Hablicsek, {\cite[Corollary~1.17]{Arinkin2014}}]\label{thm:ACH}
If $X$ is smooth projective and $G\leq \Aut(X)$ is finite, we have
\begin{align*}
\HH_\bullet([X/G]) =&\left(\bigoplus_{g\in G} \bigoplus_{q-p=\bullet} \H^p(X^g, \Omega_g^q)\right)^G.
\\\HH^\bullet([X/G]) =&\left(\bigoplus_{g\in G} \bigoplus_{p+q=\bullet} \H^{p-c_g}(X^g, \bigwedge^q \T_g\otimes \omega_g)\right)_G.
\end{align*}
Here $(\cdots)^G$ denote invariants and $(\cdots)_G$ coinvariants for the group $G$.

The $G$-action is defined as follows: For $h\in G$, $h$ acts on the direct sum of 
the  $\H^p(X^g, \Omega_g^q)$ resp.\ $\H^{p-c_g}(X^g, \bigwedge^q \T_g\otimes \omega_g)$
by the isomorphism $X^g\xrightarrow hX^{hgh^{-1}}$ which identifies
$\Omega_g^q$ with $\Omega_{hgh^{-1}}^q$, $\T_g$ with $\T_{hgh^{-1}}$ and $\omega_g$ with $\omega_{hgh^{-1}}$.
\end{theorem}
\begin{example}
Let us check that we obtain Theorem~\ref{thm:HKR} in case $G=1$:
Indeed, for $g=1$, we have $X^g=X$, $c_g=0$ so that $\omega_g=\mathcal O_X$.
Hence $\T_g = \T_X$ and $\Omega^q_g = \Omega_{X/k}^q$. This recovers precisely Theorem~\ref{thm:HKR}.
\end{example}
\begin{remark}
Similar to Corollary~\ref{cor:hkrVanishingAndSymmetry}, we get vanishing and symmetry results for $\HH^\bullet$ and $\HH_\bullet$:
\begin{itemize}
\item
$\HH_\bullet([X/G]) \cong\HH_{-\bullet}([X/G])$.\\$\T_g$ is the tangent bundle of\footnote{This can be proved using \cite[Lemma~2]{Cartan15}.} $X^g$. It follows
that for an irreducible component $Y\subseteq X^g$ of dimension $n$, $\Omega^n_g\vert_Y$ is the dualizing bundle
of $Y$. Then the proof works as in Corollary~\ref{cor:hkrVanishingAndSymmetry}.
\item
$\HH_i([X/G]) = 0$ if $i\notin [-\dim X,\dim X]$.\\Indeed, $\H^p(X^g,\Omega^q_g)=0$
unless $p\in [0,\dim X]$ (as $\Omega^q_g$ is quasi-coherent) and $q\in [0,\dim X]$ (as $\T_g$ is a vector bundle of rank $\leq \dim X$).
\item
$\HH_i([X/G]) = 0$ if $i\notin [0,2\dim X]$.\\$\H^{p-c_g}(X^g, \bigwedge^q \T_g\otimes \omega_g)=0$
unless $p \in [c_g,\dim (X^g)+c_g]$ and $q\in [0,\dim X^g]$. In particular, this is zero unless $p+q\in [0,2\dim X]$.
\end{itemize}
\end{remark}
\section{Quotients of curves}
We now specialize Theorem \ref{thm:ACH} to the case of quotients of curves.
\begin{remark}
Let $g\geq 3$ be an integer. Most curves of genus $g$ have a trivial automorphism group:
Indeed, by \cite{Poonen2000}, there exists a curve of genus $g$ having a trivial automorphism group.
Since the moduli stack of curves of genus $g$ is irreducible, the generic curve will have a trivial automorphism group,
so the curves with non-trivial automorphism group are contained in a closed substack of smaller dimension.
In other words, there are only few curves of genus $\geq 3$ to which Theorem~\ref{thm:ACH} can be applied with a non-trivial group $G$.
\end{remark}
\begin{example}[Toy example]\label{example:P1byC2}
Throughout this section, we will illustrate everything using one of the simplest examples:
We consider
$C=\mathbb P^1$ with coordinate ring $k[X,Y]$,
and the group $G=\langle \sigma\rangle\cong\mathbb Z/2\mathbb Z$,
where $\sigma$ is the involution $\sigma([x:y]) = [-x:y]$.
We will later see that the quotient stack $[C/G]$ is isomorphic to a weighted projective line.
\end{example}
\begin{lemma}\label{lem:invariantAffineCovering}
Let $X$ be a projective $k$-scheme, $P\subseteq X$ a finite set of points
and $G\leq \Aut(X)$ a finite group of $k$-scheme automorphisms.
Then there exists a $G$-invariant affine open subset $U\subseteq X$
such that $P\subseteq U$.
In particular, $X$ can be covered by $G$-invariant affine open subsets.
\end{lemma}
\begin{proof}
Replacing $P$ by $\bigcup_{g\in G} gP$, we may assume $P$ is $G$-invariant.
Embed $X\hookrightarrow \mathbb P^N_k$ for some $N>0$.
Recall that any finite subset of $\mathbb P^N_k$ is contained in an affine open subset
(this follows basically from the prime avoidance lemma). Restricting to $X$,
we find an affine open subset $V\subseteq X$ such that $P\subseteq V$.

Since projective schemes are separated, $U:=\bigcap_{g\in G} gV$ is an affine open subset
of $X$. As we assumed that $P$ is $G$-invariant, $P\subseteq U$. It is clear by construction
that $U$ is $G$-invariant.
\end{proof}
\begin{example}
Continuing Example~\ref{example:P1byC2}, we can cover $\mathbb P^1 = D(X)\cup D(Y)$.
Here, $X$ and $Y$ denote the standard coordinate functions.
Then $D(X)$ and $D(Y)$
are affine open and $G$-invariant.
\end{example}
\begin{lemma}\label{lem:localAutomorphism}
Let $(R,\mathfrak m)$ be a noetherian local $k$-algebra with residue field $k$ and let moreover $G\leq \Aut_{k\LAlg}(R)$ be a group
of $k$-algebra automorphisms of $R$. Moreover denote by $\Aut_{k\LVect}(\mathfrak m/\mathfrak m^2)$ the group
of $k$-vector space automorphisms of $\mathfrak m/\mathfrak m^2$.

There is a natural map of groups
\begin{align*}
G\rightarrow \Aut_{k\LVect}(\mathfrak m/\mathfrak m^2).
\end{align*}
If $G$ is a torsion group, this map is injective.
\end{lemma}
\begin{proof}
Each $\varphi\in \Aut_{k\LAlg}(R)$ must preserve $\mathfrak m$ (being the unique maximal ideal) and thus $\mathfrak m^2$.
Hence it restricts to an automorphism of $\mathfrak m/\mathfrak m^2$.
This gives a natural map of groups $\Aut_{k\LAlg}(R)\rightarrow \Aut_{k\LVect}(\mathfrak m/\mathfrak m^2)$.

It remains to show injectivity, under the assumption that $G$ is torsion. So suppose that $\varphi\in \Aut_{k\LAlg}(R)$ is a finite-order
automorphism inducing
the identity map on $\mathfrak m/\mathfrak m^2$.
i.e.\ for all $x\in \mathfrak m$, $\varphi(x)-x \in \mathfrak m^2$. We would like to show $\varphi=\id_R$.

\textbf{Claim 1}. For all $n\geq 1$, $\varphi$ induces the identity map on $\mathfrak m^n/ \mathfrak m^{n+1}$.

Indeed, for $x_1,\dotsc,x_n\in \mathfrak m$, $\varphi(x_1\cdots x_n) \in (x_1+\mathfrak m^2)\cdots (x_n+\mathfrak m^2)\subseteq x_1\cdots x_n + \mathfrak m^{n+1}$.

\textbf{Claim 2}. For all $n\geq 2$, $\varphi$ induces the identity map on $\mathfrak m/\mathfrak m^n$.

This claim is proved inductively, being clear for $n=2$. Now in the inductive step, let $x\in \mathfrak m$,
suppose $\varphi(x)-x\in\mathfrak m^n$. We want to show $\varphi(x)-x\in \mathfrak m^{n+1}$.
Put $y := \varphi(x)-x$. By Claim~1 and $y\in\mathfrak m^n$, we have $\varphi(y)\equiv y \pmod{\mathfrak m^{n+1}}$.
Then one inductively shows $\varphi^m(x) \equiv x+my\pmod{\mathfrak m^{n+1}}$ for all integers $m\geq 1$.

As $\varphi$ is assumed to be of finite order, there exists some $m\geq 1$ with $\varphi^m(x)=x$, hence $my\in\mathfrak m^{n+1}$.
As $m$ is invertible in $k$ (using $\fieldchar k=0$), $y\in \mathfrak m^{n+1}$. Claim~2 is proved.

\textbf{Claim 3}. $\varphi$ induces the identity map on $\mathfrak m$.

Indeed, for $x\in \mathfrak m$, we have $\varphi(x)-x\in\bigcap_{n\geq 1}\mathfrak m^n$.
By Krull's intersection theorem, we have $\bigcap_{n\geq 1}\mathfrak m^n=0$, so $\varphi(x)-x=0$.

Now we can finish the proof of $\varphi=\id_{R}$ as follows:
Consider the ring homomorphism $k\rightarrow R$ inducing the $k$-algebra structure,
and denote its image also by $k$. It follows that $R=k\oplus \mathfrak m$ as $k$-vector space.
$\varphi$ being a morphism of $k$-algebras induces the identity map on $k$. Then $\varphi$ is
the identity map on $R$ by Claim~3.
\end{proof}
\begin{corollary}\label{cor:localAutomorphism}
Let $(R,\mathfrak m)$ be a regular one-dimensional noetherian local $k$-algebra and $G\leq \Aut_{k\LAlg}(R)$ be a finite
group of $k$-algebra automorphisms. Then $G$ is cyclic.
\end{corollary}
\begin{proof}
The assumptions imply that $\mathfrak m/\mathfrak m^2$ is one-dimensional,
so by Lemma~\ref{lem:localAutomorphism}, $G$ can be embedded into
$\Aut_{k\LVect}(\mathfrak m/\mathfrak m^2)\cong k^\times$, the multiplicative group of $k$. Recall
that any finite subgroup of $k^\times$ is cyclic.
\end{proof}
\begin{proposition}\label{prop:curveFiniteGroup}
Let $C$ be a smooth and irreducible projective curve over $k$ and $G\leq \Aut(C)$ a finite group of automorphisms.
\begin{enumerate}[(a)]
\item
The set
\begin{align*}
\{c\in C \mid \exists g \in G\setminus\{1\}: g(c)=c\}
\end{align*}
is a finite set of closed points of $C$.
\item For each point $c\in C$, the stabilizer $G_c = \{g\in G\mid g(c)=c\}$ is cyclic.
\item If $g\in G\setminus\{1\}$ and $c\in C$ satisfy $g(c)=c$ and $(\mathcal O_{C,c},\mathfrak m_c)$ is the local ring at $c$, then
$g$ acts on $\mathfrak m_c/\mathfrak m_c^2$ by multiplication by some scalar $\mu\in k^\times\setminus\{1\}$.
\end{enumerate}
\end{proposition}
\begin{proof}
For each $g\in G\setminus\{1\}$, the equalizer $C^g$ of $\begin{tikzcd}C\ar[r,"\id",shift left=0.5ex]\ar[r,"g"',shift right=0.5ex]&C\end{tikzcd}$
is a closed subscheme of $C$ (as $C$ is separated), and $C^g \neq C$. As $C$ is smooth, it is reduced, hence $C^g$ is topologically a properly contained
closed subset of $C$. Then each irreducible component of $C^g$ must be zero-dimensional (as $\dim C=1$),
hence $C^g$ is a union of finitely many closed points of $C$. The set in (a) therefore equals $\bigcup_{g\in G\setminus\{1\}} C^g$,
proving (a).

Now pick a point $c\in C$ and denote by $(\mathcal O_{C,c},\mathfrak m_c)$ its local ring. Consider the stabilizer $G_c$. Each $g\in G_c$ induces
an automorphism of $\mathcal O_{C,c}$ as $k$-algebra. We claim that the induced group homomorphism
$G_c \rightarrow \Aut_{k\LAlg}(\mathcal O_{C,c})$ is injective:

Suppose $g\in G_c$ induces the identity map on $\mathcal O_{C,c}$, and let $U\subseteq C$ be a $G$-invariant affine open neighbourhood of $c$ (by Lemma~\ref{lem:invariantAffineCovering}).
Write $U\cong \Spec A$. Then $A$ is an integral domain (as $C$ is integral),
so the map $A\rightarrow\mathcal O_{C,c}$ is injective (as it is a localization map). So if $g$ induces the identity map on $\mathcal O_{C,c}$,
it induces the identity map on $A$ and thus on $U$. 
We see $U\subseteq C^g$. However, since $U$ is dense in $C$ and $C^g$ is closed, we get the contradiction $C=C^g$.

Therefore the map $G_c\rightarrow \Aut_{k\LAlg}(\mathcal O_{C,c})$ is injective. By Corollary~\ref{cor:localAutomorphism},
$G_c$ is cyclic. Moreover, Lemma~\ref{cor:localAutomorphism} shows that $G_c\rightarrow \Aut_{k\LVect}(\mathfrak m_c/\mathfrak m_c^2)\cong k^\times$ is injective,
proving (c).
\end{proof}
\begin{example}\label{example:P1byC2stabilizers}
Let $C=\mathbb P^1,G=\langle\sigma\rangle$ be as in Example~\ref{example:P1byC2}.
The fixed points of $\sigma$ is $c_1=[1:0]$ and $c_2=[0:1]$.
Since $G$ is cyclic, it is trivial that all stabilizers are cyclic.

By Proposition~\ref{prop:curveFiniteGroup}, for each of the two fixed points $c_1,c_2$,
the involution $\sigma$ acts on $\mathfrak m_{c_i}/\mathfrak m_{c_i}^2$ by multiplication by $-1$
(as the scalar is $\neq 1$ but $\sigma^2=\id$ acts by multiplication by $1$).
Geometrically, $\sigma$ is a reflection near each of the two points, so it is somewhat intuitive that
tangent vectors get \enquote{flipped}.
\end{example}
With this description of finite group actions on curves, we can simplify Theorem~\ref{thm:ACH}.
For the remainder of this section, let $C$ be a smooth irreducible projective curve over $k$ and $G\leq \Aut(C)$ a finite group of
$k$-scheme automorphisms.
\begin{lemma}\label{lem:coinvariantsVanish}
Let $g\in G\setminus\{1\}$ and $c\in C^g$. Let $(\mathcal O_{C,c},\mathfrak m_c)$ be the local ring and $\mu\in k^\times\setminus\{1\}$ be the scalar
such that $g$ acts on $\mathfrak m_c/\mathfrak m_c^2$ by multiplication by $\mu$.

With the notation as in Theorem~\ref{thm:ACH}, we have:
\begin{enumerate}[(a)]
\item $c_g=1$.
\item $C^g$ is a reduced closed subscheme of $C$,
where the points of $C^g$ are precisely those $c\in C$ such that $g(c)=c$.
\item The stalk of $N_{C^g/C}$ at $c$ is isomorphic to $\Hom_{\mathcal O_{C,c}}(\mathfrak m_c/\mathfrak m_c^2,\mathcal O_{C,c})$
where $g$ acts by multiplication by $\mu^{-1}$.
\item If $\Omega_C^1$ and $\T_C$ denote the (co)tangent bundles
and $(\Omega_C^1)_c$ resp. $(\T_C)_c$ their stalks at $c$,
$g$ acts on $\frac{(\Omega_C^1)_c}{\mathfrak m_c(\Omega_C^1)_c}$ by multiplication by $\mu$. It acts on 
$\frac{(\T_C)_c}{\mathfrak m_c (\T_C)_c}$ by multiplication by $\mu^{-1}$.
\item $\T_g=0$ and $\Omega_g=0$.
\end{enumerate}
\end{lemma}
\begin{proof}
\begin{enumerate}[(a)]
\item This is clear as $C^g$ is a zero-dimensional subscheme of $C$.
\item Let $A$ be the set of points $c\in C$ such that $g(c)=c$. Then $A$ is a finite set of closed points (Proposition~\ref{prop:curveFiniteGroup}),
so $A$ is closed. If we equip $A$ with the reduced subscheme structure, $g$ acts trivially on $A$ (it acts trivially on the points,
and each point is isomorphic to $\Spec k$, which has no non-trivial $k$-scheme automorphisms).
Hence $A\hookrightarrow C$ factors over $C^g$ by definition of $C^g$, showing $A\subseteq C^g$ as closed subschemes.

Now if $c\in C^g$ is any point of $C^g$, $c$ is a closed point of $C$ (as $\dim C^g=0$) and preserved by $g$ (by definition of $C^g$).
Hence the points of $C^g$ are precisely the closed points of $C$ fixed by $g$. It remains to show that $C^g$ is reduced.
Pick some $c\in C^g$ and note that $\mathcal O_{C,c}\rightarrow \mathcal O_{C^g,c}$ is surjective. Denote its kernel by $I$.
Then $I\subseteq \mathfrak m_c$ (as $\mathcal O_{C^g,c}\neq 0$). On the other hand, for all $x\in \mathcal O_{C,c}$,
we must have $gx-x \in I$ (by definition of $C^g$). Pick $x\in \mathfrak m_c\setminus \mathfrak m_c^2$,
then $gx = \mu x +\mathfrak m_c^2$, so that $I\ni gx-x\in \mathfrak m_c\setminus \mathfrak m_c^2$.
We see $\mathfrak m_c\supseteq I\not\subseteq \mathfrak m_c^2$. As $\dim_k \mathfrak m_c/\mathfrak m_c^2=1$,
$I=\mathfrak m_c$ so that $\mathcal O_{C^g,c} \cong k$. This shows $C^g$ is reduced.
\item Denote by $\mathcal I$ the sheaf of $\mathcal O_C$-ideals on $C$ defining the closed subscheme $C^g$.
The stalk of $\mathcal I$ at $c$ is given by $\mathcal I_c = \mathfrak m_c$ (from part (b)).
Hence the stalk of the conormal sheaf $\mathcal I/\mathcal I^2$ at $c$ is $\mathfrak m_c/\mathfrak m_c^2$,
and so the stalk of the normal sheaf is the dual, $\Hom_{\mathcal O_{C^g,c}}(\mathfrak m_c/\mathfrak m_c^2,\mathcal O_{C^g,c})$.
Since $\mathcal O_{C^g,c}\cong k$ with trivial $g$ action and $g$ acts on $\mathfrak m_c/\mathfrak m_c^2$ by multiplication by $\mu$,
$g$ acts on its dual by multiplication by $\mu^{-1}$.
\item By \cite[Proposition~II.8.7]{Hartshorne1977}, there exists a natural isomorphism
\begin{align*}
\mathfrak m_c/\mathfrak m_c^2 \rightarrow \Omega_{\mathcal O_{C,c}/k}^1 \otimes_{\mathcal O_{C,c}} k.
\end{align*}
Now since $k=\mathcal O_{C,c}/\mathfrak m_c$, we get a natural isomorphism
$\mathfrak m_c/\mathfrak m_c^2\cong \frac{\Omega_{\mathcal O_{C,c}/k}^1}{\mathfrak m_c\Omega_{\mathcal O_{C,c}/k}^1}$.
Now observe $\Omega_{\mathcal O_{C,c}/k}^1\cong (\Omega^1_C)_c$ naturally.
As $g$ acts on $\mathfrak m_c/\mathfrak m_c^2$ by multiplication by $\mu$,
it does so on $\frac{(\Omega_C^1)_c}{\mathfrak m_c(\Omega^1_C)_c}$.
Dualizing as before, we get the claim for
$\frac{(\T_C)_c}{\mathfrak m_c(\T_C)_c}$.
\item For any $x\in (\T_C)_c$, $gx = \mu^{-1} x + \mathfrak m_c (\T_C)_c$ by part (d). 
Let $U\subseteq (\T_C)_c$ be the submodule generated by $\{gx-x : x\in (\T_C)_c\}$.
As $gx-x = (\mu^{-1}-1)x + \mathfrak m_c(\T_C)_c$ and $\mu^{-1}\neq 1$,
we conclude $U + \mathfrak m_c (\T_C)_c = (\T_C)_c$. By Nakayama's lemma,
$U=(\T_C)_c$, so that the $g$-coinvariants of $(\T_C)_c$ are zero. Hence $\T_g=0$.
By definition, $\Omega_g=0$.\qedhere
\end{enumerate}
\end{proof}
With this Lemma, we can simplify Theorem~\ref{thm:ACH} for the case of curves:
\begin{proposition}\label{prop:cohomologyCurveQuotient}
With the same assumptions on $C,G$ as before, we have:
\begin{align*}
&\HH^0([C/G]) \cong k
\\&\HH^1([C/G]) \cong \H^0(C,\T_C)_G \oplus \H^1(C,\mathcal O_C)_G,
\\&\HH^2([C/G]) \cong \H^1(C,\T_C)_G,
\\&\HH^n([C/G])=0\text{ if }n\notin\{0,1,2\}.
\end{align*}
\end{proposition}
\begin{proof}
We start with Theorem~\ref{thm:ACH}:
\begin{align*}
\HH^\bullet([C/G]) \cong \left(\bigoplus_{g\in G} \bigoplus_{p+q=\bullet} \H^{p-c_g}(C^g, \bigwedge^q \T_g\otimes \omega_g)\right)_G.
\end{align*}
The inner direct sum splits naturally into the parts for $g=1$ and $g\neq 1$:
\begin{align*}
\HH^\bullet([C/G]) \cong\left(\bigoplus_{p+q=\bullet} \H^{p}(C, \bigwedge^q \T_C)\right)_G
\oplus
\left(\bigoplus_{1\neq g\in G} \bigoplus_{p+q=\bullet} \H^{p-1}(C^g, \bigwedge^q \T_g\otimes \omega_g)\right)_G.
\end{align*}
We claim that the second summand is zero:
Indeed, $\T_g=0$, so the summands for $q\neq 0$ vanish. Observe that
$C^g$ is an affine scheme,
so $\H^{p-1}(C^g,\omega_g)=0$ unless $p=1$, in which case it equals $\bigoplus\limits_{c\in C^g} (\omega_g)_c$.
Now $\omega_g = N_{C^g/C}$, and the stalk at each $c\in C^g$ vanishes when taking $g$-coinvariants (Lemma~\ref{lem:coinvariantsVanish}).
This shows that the second summand is zero, so that
\begin{align*}
\HH^\bullet([C/G]) 
\cong\left(\bigoplus_{p+q=\bullet} \H^{p}(C, \bigwedge^q \T_C)\right)_G
=\bigoplus_{p+q=\bullet} \H^p(C,\bigwedge^q \T_C)_G.
\end{align*}
For $\bullet=0$, the only direct summand is $\H^0(C,\mathcal O_C)_g$.
As $C$ is an integral $k$-scheme, $\H^0(C,\mathcal O_C)$ is a $k$-algebra and a domain.
Moreover, since $\mathcal O_C$ is coherent, $\H^0(C,\mathcal O_C)$ is finite-dimensional,
so that $\H^0(C,\mathcal O_C)\cong k$ as $k$-algebras (using that $k$ is algebraically closed). Since $G$ induces $k$-algebra-automorphisms
on $\H^0(C,\mathcal O_C)\cong k$, but $k$ has only the identity as $k$-algebra-automorphism, 
$\H^0(C,\mathcal O_C)_G\cong k$. Therefore, $\HH^0([C/G])\cong k$.

For $\HH^0([C/G])$, we get the summands for $(p,q) = (0,1),(1,0)$.
For $\HH^2([C/G])$, only the summand for $(p,q)=(1,1)$ survives,
since $\H^p(C,\mathcal F)=0$ for any quasi-coherent sheaf $\mathcal F$ if $p\geq 2>\dim C$,
and $\bigwedge^q \T_C=0$ if $q\geq 2$. Similarly $\HH^n([C/G])=0$ if $n\geq 3$.
\end{proof}
\begin{example}\label{example:P1byC2cohomology}
In the case of Example~\ref{example:P1byC2}, we get
\begin{align*}
&\HH^0([C/G])\cong k
\\&\HH^1([C/G]) \cong \H^0(\mathbb P^1, \T_{\mathbb P^1})_\sigma \oplus \H^1(\mathbb P^1,\mathcal O_{\mathbb P^1})_\sigma
\underset{\H^1(\mathbb P^1,\mathcal O_{\mathbb P^1})=0}=\H^0(\mathbb P^1, \T_{\mathbb P^1})_\sigma.
\\&\HH^2([C/G])\cong \H^1(\mathbb P^1,\T_{\mathbb P^1})_\sigma \underset{\H^1(\mathbb P^1,\T_{\mathbb P^1})=0}=0
\end{align*}
We have $\H^0(\mathbb P^1,\T_{\mathbb P^1})\cong \H^0(\mathbb P^1,\mathcal O_{\mathbb P^1}(2))\cong k^3$,
but it is not a priori clear how $\sigma$ acts on this.

Abstractly, $\H^0(\mathbb P^1,\T_{\mathbb P^1})\cong\H^0(\mathbb P^1,\mathcal O_{\mathbb P^1}(2))$
corresponds to homogeneous degree $2$ polynomials in $k[X,Y]$.
Hence a basis is given by $X^2,Y^2,XY$.

Let $c_1=[1:0],c_2=[0:1]$ be the two fixed points of $\sigma$.
By Example~\ref{example:P1byC2stabilizers}, $\sigma$ acts on $(\T_{\mathbb P^1})_{c_1}$ and
on $(\T_{\mathbb P^1})_{c_2}$ by multiplication by $-1$.
As $X\in \H^0(\mathbb P^1,\T_{\mathbb P^1})$ vanishes in $c_2$ but not in $c_1$,
we must have $\sigma(X) \equiv -X \pmod{Y}$. Similarly,
$\sigma(Y)\equiv -Y\pmod{X}$.

By the (very easy) representation theory of the cyclic group $G=\langle\sigma\rangle$,
the representation $\H^0(\mathbb P^1,\T_{\mathbb P^1})$ splits
as a direct sum $V_1\oplus V_2\oplus V_3$,
where $\sigma$ acts on $V_1$ and $V_2$ by multiplication by $-1$
and on $V_3$ by multiplication by some scalar $\lambda\in\{-1,1\}$.

If $\lambda=1$, the module of coinvariants
$\H^0(\mathbb P^1,\T_{\mathbb P^1})_\sigma$ is one-dimensional, otherwise
it is zero.

Now let $x\in\mathbb P^1$
be any closed point distinct from $c_1,c_2$. Then $Y:=\{x,\sigma(x)\}$ forms
a closed $\sigma$-invariant subscheme. It is easy to check that
$\H^0(Y,\T_{\mathbb P^1})$ contains a non-zero section that is \emph{not} multiplied by $-1$
when applying $\sigma$. Hence $\sigma$ cannot multiply every section in $\H^0(\mathbb P^1,\T_{\mathbb P^1})$
by $-1$ (as $\T_{\mathbb P^1}$ is generated by its global sections). We conclude $\lambda=1$.
Thus $\HH^1([C/G])\cong k$.
\end{example}
\begin{proposition}\label{prop:homologyCurveQuotient}
Denote by
\begin{align*}
G\setminus C = \{Gc \mid c\in C\text{ closed point}\}
\end{align*}
the set of $G$-orbits in $C$. The Hochschild homology of $[C/G]$ is given by
\begin{align*}
&\dim \HH_0([C/G]) = 2 + \sum_{Gc\in G\setminus C} (\# G_c-1),
\\&\HH_1([C/G]) \cong \HH_{-1}([C/G]) \cong \H^0(C,\Omega_C^1)^G \cong \H^1(C,\mathcal O_C)^G,
\\&\HH_n([C/G])=0\text{ if }n\notin\{-1,0,1\}.
\end{align*}
\end{proposition}
\begin{proof}
We again start with
\begin{align*}
\HH_\bullet([C/G]) \cong&\left(\bigoplus_{g\in G} \bigoplus_{q-p=\bullet} \H^p(C^g, \Omega_g^q)\right)^G
\\\cong &
\left(\bigoplus_{q-p=\bullet} \H^p(C, \Omega_C^q)\right)^G
\oplus
\left(\bigoplus_{1\neq g\in G} \bigoplus_{q-p=\bullet} \H^p(C^g, \Omega_g^q)\right)^G
\end{align*}
Now $\Omega_g=0$ and $C^g$ is affine for $g\neq 1$, so $\H^p(C^g, \Omega_g^q)=0$ unless $p=q=0$,
in which case $\H^p(C^g, \Omega_g^q)\cong \mathcal O_{C^g,c}\cong k$.

Similarly, $\H^p(C, \Omega_C^q)=0$ unless $p,q\in\{0,1\}$.
This shows $\HH_n([C/G])=0$ for $n\notin\{-1,0,1\}$ and
\begin{align*}
\HH_1([C/G]) \cong \H^0(C,\Omega_C)^G,\qquad
\HH_{-1}([C/G])\cong \H^1(C,\mathcal O_C)^G.
\end{align*}
These are isomorphic by Serre duality.

So only the computation of $\HH_0([C/G])$ remains.
We get
\begin{align*}
\HH_0([C/G]) \cong \H^0(C,\mathcal O_C)^G \oplus \H^1(C,\Omega^1_C)^G \oplus \left(\bigoplus_{1\neq g\in G} \bigoplus_{c\in C^g} k\right)^G.
\end{align*}
As in the proof of Proposition~\ref{prop:cohomologyCurveQuotient}, $\H^0(C,\mathcal O_C)\cong k$ with trivial $G$-action.
By Serre duality, also $\H^1(C,\Omega_C^1)^G\cong k$.

A $k$-basis of $V:=\bigoplus\limits_{1\neq g\in G} \bigoplus\limits_{c\in C^g} k$ is given by
\begin{align*}
B = \{(g,c) \mid 1\neq g\in G, c\in C^g\}.
\end{align*}
The action of $G$ on $V$ is given by sending for $h\in G$ the basis vector $(g,c)$ to $(hgh^{-1},h(c))$.
So any element of $v\in V^G$ is of the form $v=\sum\limits_{b\in B}\lambda_b b$ with $gv = v$,
i.e.\ $\lambda_{gb}=\lambda_b$ for all $g\in G$. Hence a $k$-basis of $V^G$
is given by
\begin{align*}
\Bigl\{\sum_{b'\in Gb} b'\mid Gb\in G\setminus B\text{ a $G$-orbit}\Bigr\}.
\end{align*}
We conclude
\begin{align*}
\dim V^G =\#\{Gb\in G\setminus B\text{ a $G$-orbit}\} = \# (G\setminus B).
\end{align*}
Let $A := \{c\in C : G_c\neq \{1\}\}$ be the set of points of $C$ with non-trivial stabilizer; this is a finite set closed of closed points of $C$ by Proposition~\ref{prop:curveFiniteGroup}.
Then $B = \{(g,c)\in G\times A \mid g(c)=c\}$. The action of $G$ on $C$ restricts to an action of $A$,
so let $R\subseteq A$ be a class of representatives for the $G$-orbits.
Every element of $B$ is conjugate, via the $G$-action,
to an element of the form $(g,c)$ with $c\in R$ and $g\in G_c\setminus\{1\}$.
Moreover, two such elements $(g,c)$ and $(g',c')$ are conjugate
iff $c=c'$ and there exists some $h\in G$ with $h(c)=c'$ and $g'=h^{-1}gh$.
i.e.\ $(c,g)$ and $(c',g')$ are conjuate via the $G$-action
iff $c=c'$ and $g$ and $g'$ are conjugate in $G_x$. We see
\begin{align*}
\# (G\setminus B) = \sum_{c\in R} \#(\text{$G_c$-congugacy classes in $G_c\setminus\{1\}$}).
\end{align*}
As $G_c$ is cyclic (Proposition~\ref{prop:curveFiniteGroup}), this equals
\begin{align*}
\dim_k V^G = \sum_{c\in R} \#G_c - 1.
\end{align*}
Now if $c'$ is an element of the $G$-orbit of $c$, which we denote by $Gc$,
then $G_{c'}$ and $G_c$ are conjugate in $G$, so in particular $\#G_{c'}=\#G_c$.
Hence the number $\#G_c$ depends only on the $G$-orbit of $c$.
Therefore
\begin{align*}
\dim_k V^G = \sum_{Gc\in G\setminus A}\# G_c-1
\underset{\#G_c=1\text{ if }c\notin A} = \sum_{Gc\in G\setminus A}\# G_c-1.
\end{align*}
Since $\HH_0([C/G]) \cong k^2 \oplus V^G$, this finishes the proof.
\end{proof}
\begin{example}\label{example:P1byC2homology}
In the situation of Example~\ref{example:P1byC2},
with $c_1=[1:0], c_2=[0:1]$, we get
\begin{align*}
&\dim\HH_0([C/G]) \underset{\text{E\ref{example:P1byC2stabilizers}}}= 2+(\#G_{c_1}-1)+(\#G_{c_2}-1) = 2+1+1=4.
\\&\HH_{\pm 1}([C/G])\cong \H^1(\mathbb P^1,\mathcal O_{\mathbb P^1})^G\underset{\H^1(\mathbb P^1,\mathcal O_{\mathbb P^1})=0}=0.
\end{align*}
\end{example}
\section{Iterated root stacks}
Let $C,G$ be as before.
Using the propositions just proved, we can calculate the Hochschild
(co)homology of the quotient stack $[C/G]$.
By \cite{Geraschenko2015} or equivalently \cite[Lemma~5.3.10]{Voight2019},
$[C/G]$ can be identified as an iterated root stack.
The aim of this section is to describe the Hochschild (co)homology of $[C/G]$ in terms of this iterated root stack.
\begin{definition}
Let $C,G$ be as before.
\begin{enumerate}[(a)]
\item A \emph{$G$-invariant morphism} is a morphism of $k$-schemes $C\xrightarrow\varphi D$ such that $\varphi g=\varphi$ for all $g\in G$.
\item A $G$-invariant morphism $C\xrightarrow\varphi D$ is called \emph{quotient scheme} if it is universal in the following sense:
For each $G$-invariant morphism $C\xrightarrow\psi D'$, there exists a unique morphism of $k$-schemes $D\xrightarrow\delta D'$ with $\psi=\delta\varphi$.

If it exists, we denote the quotient scheme by $C\xrightarrow\pi C/G$.
\end{enumerate}
\end{definition}
\begin{remark}
A quotient scheme always exists; it is locally given by $\Spec A\rightarrow \Spec(A^G)$ for $G$-invariant affine open subsets of $C$.
It is unique up to unique isomorphism.
The fibres $\pi^{-1}(x)$ for closed points $x\in C/G$ are precisely $G$-orbits of points in $C$.
\end{remark}
\begin{example}\label{example:P1byC2coarse}
In the case of Example~\ref{example:P1byC2},
the quotient $\mathbb P^1/\langle \sigma\rangle$
is isomorphic to $\mathbb P^1$, the quotient map being
\begin{align*}
\mathbb P^1\rightarrow\mathbb P^1, [x:y]\mapsto [x^2:y^2].
\end{align*}
\end{example}
\begin{proposition}\label{prop:quotientRootStack}
Let $C\xrightarrow\pi C/G$ be the quotient scheme. Then $[C/G]$ is the iterated root stack
obtained from $C/G$ by attaching weight $G_c$ to the point $\pi(c)$ for all $c\in C$.
\end{proposition}
\begin{proof}
Denote this iterated root stack by $\mathcal X$ (it is well-defined, since if $\pi(c)=\pi(c')$, then $c$ and $c'$ are
in the same $G$-orbit, so $G_c$ and $G_{c'}$ are conjugate, and only finitely many points in $C$ have non-trivial stabilizer;
cf. Proposition~\ref{prop:curveFiniteGroup}).
Then the coarse moduli space of $\mathcal X$ is $C/G$. By \cite[Example~5.3.6]{Voight2019}, this is also the coarse moduli space
of $[C/G]$. The result now follows from \cite[Lemma~5.3.10]{Voight2019}.
\nocite{Voight2019}
\end{proof}
\begin{example}In the situation of Example~\ref{example:P1byC2coarse},
we conclude that $[C/G]=[\mathbb P^1/\langle \sigma\rangle]$ is isomorphic to the weighted projective
line $\mathbb P^1\langle \infty,1; 2,2\rangle$. We see that the Examples \ref{example:P1byC2cohomology}
and \ref{example:P1byC2homology} coincide with Proposition~\ref{prop:cohomologyCanonical2}.
\end{example}
We may expect that the Hochschild (co)homology of $[C/G]$ depends
only on the quotient scheme $C\xrightarrow\pi C/G=:Q$ and the
points $\pi(c)\in Q$ with $G_c\neq \{1\}$ together with the order $\#G_c$, as this data
fully determines $[C/G]$.

Proposition~\ref{prop:homologyCurveQuotient} gives such a description for $\HH_0([C/G])$.
In the remainder of this section, we will prove the following theorem:
\begin{theorem}\label{thm:summaryRootStack}
Let $X$ be a smooth irreducible projective curve, $x_1,\dotsc,x_m\in X$ be distinct points and $e_1,\dotsc,e_m\in\mathbb Z_{\geq 2}$.
Let $\mathcal X$ be the iterated root stack obtained from $X$ by attaching the weight $e_i$ to $x_i$ for $i=1,\dotsc,m$.

Denote by $g$ the genus of $X$ and define
\begin{align*}
e := e_1+\cdots + e_m,\quad a:=\begin{cases}3,&g=0\\1,&g=1\\0,&g\geq 2\end{cases},\quad
r:=\min\{a,m\},\quad d := \begin{cases}0,&g=0\\1,&g=1\\3g-3,&g\geq 2\end{cases}.
\end{align*}
The Hochschild (co)homology of $\mathcal X$ in degree $i\in \mathbb Z$ has dimension
\begin{align*}
&\dim\HH_i(\mathcal X)=\begin{cases}2+e-m,&i=0\\
g,&i=\pm 1\\
0,&\text{otherwise}\end{cases}\\
&\dim\HH^i(\mathcal X)=\begin{cases}1,&i=0\\
g+a-r,&i=1\\
d+m-r,&i=2\\
0,&\text{otherwise}.\end{cases}
\end{align*}
\end{theorem}
\begin{example}\label{example:calculationWPLFromAbstractResult}
Let us apply the theorem to weighted projective lines.
Here, we may choose $X=\mathbb P^1$ so that $g=0$. Then $x_1,\dotsc,x_m$
are the points where we attach the weights to and
$e_1,\dotsc,e_m$ are the weights. The Hochschild (co)homology of $\mathcal X=\mathbb P^1\langle x_1,\dotsc,x_m;e_1,\dotsc,e_m\rangle$ in degree $j\in\mathbb Z$
is given by
\begin{align*}
&\dim\HH_j(\mathcal X)=\begin{cases}2+\sum\limits_{i=1}^m (e_i-1),&j=0\\
0,&j\neq 0\end{cases}\\
&\dim\HH^j(\mathcal X)=\begin{cases}1,&j=0\\
3-r=\max\{0,3-m\},&j=1\\
m-r=\max\{0,m-3\},&j=2\\
0,~&\text{otherwise}\end{cases}
\end{align*}
This exactly reproduces Corollary~\ref{cor:happelWPL}.
\end{example}
\begin{remark}
\begin{enumerate}[(a)]
\item For a smooth projective curve $X$, the HKR-theorem \ref{thm:HKR} and Serre duality
imply $\HH_i(X)\cong \HH_{-i}(X)$ for all integers $i$, $\HH_i(X)=0$ if $i\notin\{-1,0,1\}$
and $\HH^i(X)=0$ if $i\notin\{0,1,2\}$ (Corollary~\ref{cor:hkrVanishingAndSymmetry}).

By Theorem~\ref{thm:summaryRootStack}, the same properties are satisfied for iterated root stacks obtained from smooth projective curves.
\item If $\mathcal X$ is a weighted projective line, the previous observation shows $\HH_1(\mathcal X)\cong \HH_{-1}(\mathcal X)$.
If we moreover use that $\mathcal X$ is derived equivalent to an algebra (Proposition~\ref{prop:wplTiltingObject}),
then $\HH_{-1}(\mathcal X)=0$ follows immediately.

Conversely, if $X$ is a smooth projective curve of genus $g>0$ and $\mathcal X$ an iterated root stack
obtained from $X$, the same argument as above shows that $\mathcal X$ cannot be derived equivalent to an algebra.
Therefore, the techniques in Chapter~\ref{chapter-3} cannot be used to compute the Hochschild cohomology of $\mathcal X$.
\item We again see that the dimension of the Hochschild cohomology groups
depends only on the number of exceptional points and the genus, and not on the actual points
or their weights. Similarly, the Hochschild homology depends only on the genus and the weights,
and not the points. After the proof of the theorem, we will indicate why this is the case.
\item The statement on Hochschild homology can alternatively be proved as follows:

By Proposition~\ref{prop:rootStackSemiOrthogonal}, the derived category of a root stack $\tilde X$
associated to the stack $X$, the effective Cartier divisor $D\subseteq X$ and $r\geq 1$ allows a semi-orthogonal decomposition
\begin{align*}
\D^\bounded(\tilde X) =\Bigl\langle \D^\bounded(X),\underbrace{\D^\bounded(D),\dotsc,\D^\bounded(D)}_{r-1\text{ times}}\Bigr\rangle.
\end{align*}
By \cite{Kuznetsov2009}, it follow that $\HH_\bullet(\tilde X)\cong \HH_\bullet(X)\oplus \HH_\bullet(E)^{r-1}$.
If $E$ is a point, we see $\HH_0(\tilde X)\cong \HH_0(X)\oplus k^{r-1}$ and $\HH_i(\tilde X)\cong \HH_i(X)$ for $i\in\mathbb Z\setminus\{0\}$.
Then the claim on Hochschild homology of iterated root stacks in Theorem~\ref{thm:summaryRootStack} follows
by induction on the number of exceptional points, the inductive start being the HKR decomposition \ref{thm:HKR}.
\item
The numbers $d$ and $a$ in Theorem~\ref{thm:summaryRootStack} have the following background,
motivating their names: $d$ is the dimension of the coarse moduli space of curves of genus $g$,
and $a$ is the dimension of the automorphism group of a generic curve. We will not use these facts. Note the relation $d-a = 3g-3$.
\end{enumerate}
\end{remark}
We first prove the theorem for the quotient stack $\mathcal X=[C/G]$.
\begin{lemma}\label{lem:differentialsOfQ}
Denote by $\Omega^1_{Q/k}$ the cotangent sheaf of the quotient scheme $C\xrightarrow\pi C/G=:Q$.
Let $c\in C$ be a closed point.
\begin{enumerate}[(a)]
\item $Q$ is a smooth projective curve and $C\xrightarrow\pi Q$ is a finite and flat morphism.
\item There is a canonical exact sequence of sheaves of $\mathcal O_C$-modules
\begin{align*}
0\rightarrow \pi^\ast\Omega^1_{Q/k}\rightarrow\Omega^1_{C/k}\rightarrow \Omega^1_{C/Q}\rightarrow 0
\end{align*}
\item Let $u\in \mathcal O_{C,c}$ be a uniformizing element (i.e.\ a generator of the maximal ideal).
Then the stalk of $\Omega^1_{C/k}$ is a free $\mathcal O_{C,c}$-module with generator $du$.
\item
The stabilizer of $c$, $G_c$, acts in a canonical way on $\mathcal O_{C,c}$, and $\mathcal O_{Q,\pi(c)}$ can be
identified with the subalgebra of $G_c$-invariants.
\item Let $e:=\# G_c$. Then $t:=u^e \in \mathcal O_{Q,\pi(c)}$ is a uniformizing element.
In particular, $dt$ generates $\Omega_{Q,\pi(c)}$. The map $\pi^\ast\Omega^1_{Q/k}\rightarrow\Omega^1_{C/k}$
in (b) is given by $dt\mapsto dt$.
\end{enumerate}
\end{lemma}
\begin{proof}
\begin{enumerate}[(a)]
\item By Lemma~\ref{lem:invariantAffineCovering}, we may cover $C$ by $G$-invariant affine open subsets $U$.
Fix such a $U$ and an isomorphism $U\cong \Spec A$ for a $k$-algebra $A$. By construction of $Q$, $\pi(U)\subseteq Q$ is an affine open subset
with $\pi(U)\cong A^G$, the subalgebra of $G$-invariants.

Note that $A$ is a Dedekind domain. Let $L$ be the fraction field of $A$ and $K:= L^G$ be the subfield of $G$-invariants.
Then $L/K$ is a Galois extension with Galois group $G$. By standard algebraic number theory,
$A^G = K\cap A$ is a Dedekind domain and $A^G\hookrightarrow A$ a finite and flat ring extension ($A$ is a free $A^G$-module).

We conclude that $\pi(U)\cong\Spec A$ is a regular one-dimensional $k$-scheme and $U\xrightarrow\pi \pi(U)$ is finite and flat.
Glueing these together, we see that $C\xrightarrow\pi Q$ is finite and flat and $Q$ is a smooth curve.

If $C=\Proj S$ for a graded $k$-algebra $S$, $Q$ can be realized as $\Proj S^G$, so $Q$ is again projective.
\item We are now in the situation of \cite[Section~IV.2]{Hartshorne1977}. In particular, we may apply \cite[Proposition~IV.2.1]{Hartshorne1977}.
\item Note that $\Omega^1_{C/k}$ is a discrete valuation ring. Then this is a standard fact (cf.~\cite[Section~IV.2]{Hartshorne1977}).
\item Pick a $G$-invariant affine open subset of $c$ (Lemma~\ref{lem:invariantAffineCovering}). Localizing at the maximal ideal defining $c$,
we see that $G_c$ acts on $\mathcal O_{C,c}$. Then the fact that $\mathcal O_{Q,\pi(c)}$ is the subalgebra of invariants follows from the construction of $Q$.
\item Since $\mathcal O_{Q,\pi(c)}\hookrightarrow \mathcal O_{C,c}$ is a finite extension of discrete valuation rings,
there is a uniformizing element of $\mathcal O_{Q,\pi(c)}$ of the form $u^\alpha$ for some $\alpha>0$.

Recall that the action of $G$ on $\mathfrak m_c/\mathfrak m_c^2\cong k$ is faithful.
Hence a generator $g\in G$ acts on $\mathfrak m_c/\mathfrak m_c^2$ by multiplication by some primitive $e$-th root of unity $\zeta$.
Thus $gu \cong \zeta u\pmod{\mathfrak m_c^2}$. Then $gu^\alpha\cong \zeta^\alpha u^\alpha\pmod{u^{\alpha+1}\mathcal O_{C,c}}$.
So if $u^\alpha$ is $G$-invariant, $e\mid \alpha$. In particular, $\alpha\geq e$.

On the other hand, the Galois-argument of (a) shows that $\alpha\leq e$. Hence $t=u^e$ is a uniformizing element of $\mathcal O_{Q,\pi(c)}$.
The claim on $\pi^\ast\mathcal O_{Q,\pi(c)}\rightarrow\mathcal O_{C,c}$ follows from the construction of the sequence in (b)
given in \cite{Hartshorne1977}.
\qedhere\end{enumerate}
\end{proof}
\begin{lemma}\label{lem:abstractDifferentialAction}
Let $g\in G$. $g^\ast\Omega^1_{C/k}$ is a sheaf of $\mathcal O_{C}$-modules satisfying the same universal property
as $\Omega^1_{C/k}$, so there is a canonical isomorphism $\Omega^1_{C/k}\xrightarrow\cong g^\ast\Omega^1_{C/k}$.
We similarly get an isomorphism $\Omega^1_{C/Q}\xrightarrow\cong g^\ast\Omega^1_{C/Q}$.

Morever, we have a canonical isomorphism
\begin{align*}
g^\ast(\pi^\ast\Omega^1_{Q/k}) \cong (\pi g)^\ast \Omega^1_{Q/k} = \pi^\ast\Omega^1_{Q/k}.
\end{align*}
\begin{enumerate}[(a)]
\item
These give rise to a commutative diagram
\begin{align}\label{diag:differentialAction}
\begin{tikzcd}[ampersand replacement=\&]
0\ar[r]\&\pi^\ast\Omega^1_{Q/k}\ar[d,"\cong"]\ar[r]\&\Omega^1_{C/k}\ar[d,"\cong"]\ar[r]\&\Omega^1_{C/Q}\ar[r]\ar[d,"\cong"]\&0\\
0\ar[r]\&g^\ast\pi^\ast\Omega^1_{Q/k}\ar[r]\&g^\ast\Omega^1_{C/k}\ar[r]\&g^\ast\Omega^1_{C/Q}\ar[r]\&0
\end{tikzcd}
\end{align}
where the top row comes from Lemma~\ref{lem:differentialsOfQ} part (b), and the bottom row is $g^\ast$ applied to the top row.
\item
Let $c\in C$ with uniformizer $u\in \mathcal O_{C,c}$, $e=\#G_c$ and $t:=u^e$.
Pick a $G$-invariant affine open subset $U\cong \Spec A$ of $c$. Then the induced
action of $g$ on $A$ identifies $\mathcal O_{C,c}$ with $\mathcal O_{C,gc}$.

The stalk of\footnote{By the stalk of something related to the category
of $\mathcal O_C$-modules (an object, a morphism, a diagram),
we mean its image under the stalk functor $\mathcal O_C\LMod\rightarrow\mathcal O_{C,c}\LMod$.} the diagram (\ref{diag:differentialAction}) at $c$ can be identified with
\begin{align}\label{diag:localDifferentialAction}
\begin{tikzcd}[ampersand replacement=\&, column sep=4em]
0\ar[r]\&\mathcal O_{C,c} dt\ar[d,"{dt\mapsto d(gt)}"]\ar[r,"{dt\mapsto dt}"]\&\mathcal O_{C,c}du\ar[d,"{du\mapsto d(gu)}"]\ar[r]\&\frac{\mathcal O_{C,c}du}{\mathcal O_{C,c}dt}\ar[r]
\ar[d,"{dt\mapsto d(gt)}"]\&0\\
0\ar[r]\&\mathcal O_{C,gc}d(gt)\ar[r,"{d(gt)\mapsto d(gt)}"]\&\mathcal O_{C,gc}d(gu)\ar[r]\&\frac{\mathcal O_{C,gc}d(gu)}{\mathcal O_{C,gc}d(gt)}\ar[r]\&0
\end{tikzcd}
\end{align}
If moreover $gc=c$, we have $gt=t$.
\end{enumerate}
\end{lemma}
\begin{proof}
It suffices to prove (b), as commutativity of (\ref{diag:localDifferentialAction}) is obvious, and commutativity of (\ref{diag:differentialAction}) may be checked locally.

Proof of (b): The identification in the top (and hence also the bottom) row are from Lemma~\ref{lem:differentialsOfQ}. The vertical isomorphisms
come directly from the construction of the corresponding vertical isomorphisms in (\ref{diag:differentialAction}).

To the last claim: Assume $g\in G_c$. Then $gt=t$, i.e.\ $t$ is $g$-invariant, comes from Lemma~\ref{lem:differentialsOfQ} part (e).
\end{proof}
\begin{lemma}\label{lem:relativeDifferentialsVanish}
Let $c\in C$ be a closed point, $u\in \Omega_{C,c}$ a uniformizing element, $g$ a generator of the cyclic group $G_c$ and $e:=\# G_c$.
Then the stalk of $\Omega^1_{C/Q}$ has a $k$-basis
$du,udu,\dotsc,u^{e-2}du$. The action of $g$ with respect to this basis has matrix
\begin{align*}
\begin{pmatrix}
\zeta&\ast&\hdots&\ast\\
0&\zeta&&\vdots\\
\vdots&&\ddots&\ast\\
0&\hdots&0&\zeta
\end{pmatrix}.
\end{align*}
for some primitive $e$-th root of unity $\zeta$.
\end{lemma}
\begin{proof}
This follows directly from $\Omega^1_{C/Q}\cong \frac{\mathcal O_{C,c} du}{\mathcal O_{C,c} d(u^e)}$ and $gu\cong \zeta u\pmod{u^2}$
for some primitive $e$-th root of unity $\zeta$ ($G_c$ acting faithfully on $\mathfrak m_c/\mathfrak m_c^2$).
\end{proof}
\begin{lemma}\label{lem:invariantsInverseImage}
Let $M$ be a quasi-coherent sheaf of $\mathcal O_Q$-modules and $i\geq 0$.
\begin{enumerate}[(a)]
\item
$G$ acts on $\H^i(C,\pi^\ast M)$ as follows:
For $g\in G$, we have an isomorphism of functors
\begin{align*}
\H^0(C,g^\ast(-))\cong \H^0(C,(g^{-1})_\ast(-)) = \H^0(C,(-)),
\end{align*}
hence $\H^i(C,g^\ast(-))\cong \H^i(C,-)$. Moreover, we have an isomorphism
\begin{align*}
g^\ast\pi^\ast M \cong (\pi g)^\ast M = \pi^\ast M.
\end{align*}
These give rise to an isomorphism
\begin{align*}
g: \H^i(C,\pi^\ast M)\cong\H^i(C,g^\ast\pi^\ast M)\xrightarrow{\H^i(C,g^\ast\pi^\ast M\xrightarrow{\cong}\pi^\ast M)}\H^i(C,\pi^\ast M).
\end{align*}
\item There is an isomorphism, natural in $M$
\begin{align*}
\H^i(C,\pi^\ast M)^G\cong \H^i(Q,M)\cong \H^i(C,\pi^\ast M)_G.
\end{align*}
This isomorphism is compatible with restriction to $G$-invariant open subsets of $C$.
\end{enumerate}
\end{lemma}
\begin{proof}
(a) is trivial to prove.

For (b), note that any additive functor $kG\LMod\rightarrow k\LVect$ is exact: $kG$ is semi-simple by Maschke's theorem (we assume $\fieldchar k=0$ throughout),
so any short exact sequence in $kG\LMod$ is split. In particular, taking (co)invariants is an exact functor.
Moreover, coinvariants are naturally isomorphic to invariants, as both functors correspond to anihilating
all non-trivial simple $kG$-modules.

Recall that $\pi^\ast : \Qcoh(\mathcal O_Q)\rightarrow \Qcoh(\mathcal O_C)$ is exact by the flatness in Lemma~\ref{lem:differentialsOfQ} part (a).

By homological algebra, it follows that $\H^i(C,\pi^\ast\blank)_G$ is the $i$-th right derived
functor of $\H^0(C,\pi^\ast(\blank))_G$. Hence it suffices to show that $\H^0(C,\pi^\ast(\blank))_G\cong \H^0(Q,\blank)$.

Let $C=\bigcup_{i\in I} U_i$ be a cover of $C$ by $G$-invariant affine open subsets (cf.\ Lemma~\ref{lem:invariantAffineCovering}).
Write $U_i\cong \Spec A_i$, so $A_i$ is a Dedekind domain.
As in the proof of Lemma~\ref{lem:differentialsOfQ}, $Q=\bigcup_{i\in I}\pi(U_i)$
is an affine open covering of $Q$ where $\pi(U_i)\cong\Spec A_i^G$.

For each of these $U_i$, we have
\begin{align*}
\H^0(U_i,\pi^\ast M)_G \cong \left(\H^0(\pi(U_i),M)\otimes_{A_i^G} A_i\right)_G.
\end{align*}
The action of $G$ on $\H^0(\pi(U_i),M)\otimes_{A_i^G} A_i$ comes from the action on $A_i$.
By right exactness of the tensor product, we conclude
\begin{align*}
\left(\H^0(\pi(U_i),M)\otimes_{A_i^G} A_i\right)_G\cong~&\H^0(\pi(U_i),M)\otimes_{A_i^G} (A_i)_G
\\\cong~&\H^0(\pi(U_i),M)\otimes_{A_i^G} (A_i)^G\cong\H^0(\pi(U_i),M).
\end{align*}
It is not too hard to check that this isomorphism is compatible with restrictions.
We hence get the commutative diagram with exact rows (using exactness of coinvariants again)
\begin{align*}
\begin{tikzcd}[ampersand replacement=\&]
0\ar[r]\ar[d]\&\H^0(C,\pi^\ast M)_G\ar[r]\ar[d,dotted]\&\bigoplus_{i\in I} \H^0(U_i,\pi^\ast M)_G\ar[d,"\cong"]\ar[r]\&\bigoplus_{i,j\in I} \H^0(U_i\cap U_j,\pi^\ast M)_G\ar[d,"\cong"]\\
0\ar[r]\&\H^0(Q,\pi^\ast M)\ar[r]\&\bigoplus_{i\in I} \H^0(\pi(U_i), M)\ar[r]\&\bigoplus_{i,j\in I} \H^0(\pi(U_i)\cap \pi(U_j),M).
\end{tikzcd}
\end{align*}
This induces the dotted isomorphism which was the one we searched for.
\end{proof}
\begin{remark}
It is \emph{not} in general true that $\H^0(C,\pi^\ast M)\cong kG\otimes_k \H^0(Q,M)$ as $kG$-modules,
even though this holds locally. For example if $M=\mathcal O_Q$, we have $\pi^\ast M\cong \mathcal O_C$,
so $\H^0(C,\pi^\ast \mathcal O_Q)\cong k$ with trivial $G$-action.
\end{remark}
\begin{proposition}\label{prop:HH1}
\begin{align*}
\HH_1([C/G])\cong \HH_{-1}([C/G])\cong \H^0(Q,\Omega^1_{Q/k}).
\end{align*}
\end{proposition}
\begin{proof}
We apply Proposition~\ref{prop:homologyCurveQuotient}. So we claim that $\H^0(C,\Omega^1_{C/k})^G\cong \H^0(Q,\Omega^1_{Q/k})$.

The action of $G$ on $\H^0(C,\Omega^1_{C/k})$ is given as follows: For $g\in G$, consider the isomorphisms
$\Omega^1_{C/k}\cong g^\ast\Omega^1_{C/k}$ from Lemma~\ref{lem:abstractDifferentialAction}
and $\H^0(C,\blank)\cong \H^0(C,g^\ast(\blank))$. Then the coposite isomorphism
\begin{align*}
\\\H^0(C,\Omega^1_{C/k})\xrightarrow{H^0(C,\Omega^1_{C/k}\cong g^\ast\Omega^1_{C/k})}
\H^0(C,g^\ast\Omega^1_{C/k})\xrightarrow\cong \H^0(C,\Omega^1_{C/k})
\end{align*}
is the action of $g$ on $\H^0(C,\Omega^1_{C/k})$.

We similar get an action on $\H^0(C,\Omega^1_{Q/k})$ and $\H^0(C,\Omega^1_{C/Q})$.

Left exactness of $\H^0(C,\blank)$ and (\ref{diag:differentialAction}) give the commutative diagram with exact rows
\begin{align*}
\begin{tikzcd}[ampersand replacement=\&]
0\ar[r]\&\H^0(C,\pi^\ast \Omega^1_{Q/k})\ar[r]\ar[d,"g"]\&\H^0(C,\Omega^1_{C/k})\ar[r]\ar[d,"g"]\&\H^0(C,\Omega^1_{C/Q})\ar[d,"g"]\\
0\ar[r]\&\H^0(C,\pi^\ast\Omega^1_{Q/k})\ar[r]\&\H^0(C,\Omega^1_{C/k})\ar[r]\&\H^0(C,\Omega^1_{C/Q})
\end{tikzcd}.
\end{align*}
Recall that taking $G$-invariants is a left exact functor (in our situation, it is even exact).
This shows that we get an exact sequence of $k$-vector spaces
\begin{align*}
0\rightarrow \H^0(C,\pi^\ast\Omega^1_{Q/k})^G\rightarrow \H^0(C,\Omega^1_{C/k})^G\rightarrow \H^0(C,\Omega^1_{C/Q})^G.
\end{align*}

\textbf{Claim.} $\H^0(C,\Omega^1_{C/Q})^G=0$: By Lemma~\ref{lem:relativeDifferentialsVanish}, $\Omega^1_{C/Q}$
is supported only at the (finitely many!) point $c\in C$ with non-trivial stabilizer $G_c$.
Moreover, for each such $c\in C$, the action of $G$ on $\Omega^1_{C/Q}$
restricts to an action of $G_c$ on $(\Omega^1_{C/Q})_c$ with vanishing module of invariants (by Lemma~\ref{lem:relativeDifferentialsVanish}).
Hence $\H^0(C,\Omega^1_{C/Q})^G\cong \Bigl(\bigoplus\limits_{c\in C} (\Omega^1_{C/Q})_c\Bigr)^G=0$.

From the claim and the above short exact sequence
\begin{align*}
\H^0(C,\Omega^1_{C/k})^G\cong \H^0(C,\pi^\ast\Omega^1_{Q/k})\underset{\text{L\ref{lem:invariantsInverseImage}}}\cong \H^0(Q,\Omega^1_{Q/k}).
\end{align*}
This finishes the proof.
\end{proof}
We just described the Hochschild homology in terms of $Q$ and the orders of stabilizers.
We aim to do the same for Hochschild cohomology.
\begin{lemma}\label{lem:tangentSES}
There is a short exact sequence of quasi-coherent $\mathcal O_C$-modules
\begin{align}\label{diag:tangentSES}
0\rightarrow \T_C\rightarrow \pi^\ast \T_Q\rightarrow \mathcal F\rightarrow 0,
\end{align}
where $\mathcal F$ is a sheaf of $\mathcal O_C$-modules supported at those points $c\in C$ where $G_c\neq\{1\}$.

The short exact sequence (\ref{diag:tangentSES}) can locally described as follows: For $c\in C$, let $u\in \mathcal O_{C,c}$ be a uniformizing element,
$e = \#G_c$ and $t:=u^e\in \mathcal O_{Q,\pi(c)}$.
Then $(\T_C)_c=\Hom_{\mathcal O_{C,c}}((\Omega^1_{C/k})_c,\mathcal O_{C,c})$
is a free $\mathcal O_{C,c}$-module with basis $\delta^u : (\Omega^1_{C/k})_c\rightarrow\mathcal O_{C,c}, du\mapsto 1$.

Similarly, $(\T_Q)_{\pi(c)}$ is a free $\mathcal O_{Q,\pi(c)}$-module
with basis $\delta^t:dt\mapsto 1$.

The stalk of (\ref{diag:tangentSES}) at $c$ is given by
\begin{align*}
0\rightarrow \mathcal O_{C,c}\delta^u\xrightarrow{\delta^u\mapsto eu^{e-1}\delta t} \mathcal O_{C,c}\delta^t\rightarrow\frac{\mathcal O_{C,c}\delta^t}{\mathcal O_{C,c}u^{e-1}\delta^t}
\rightarrow 0.
\end{align*}
\end{lemma}
\begin{proof}
Apply the left-exact functor $\SHom(\blank,\mathcal O_C)$ to the short exact sequence from Lemma~\ref{lem:differentialsOfQ}.
This gives an exact sequence
\begin{align*}
0\rightarrow\SHom(\Omega^1_{C/Q},\mathcal O_C)\rightarrow\SHom(\Omega^1_{C/k},\mathcal O_C)
\rightarrow\SHom(\pi^\ast\Omega^1_{Q/k},\mathcal O_C).
\end{align*}
Since $\Omega^1_{C/Q}$ is a torsion sheaf, 
$\SHom(\Omega^1_{C/Q},\mathcal O_C)=0$. We hence get a short exact sequence of the form
\begin{align}\label{diag:tangentSes2}
0\rightarrow\SHom(\Omega^1_{C/k},\mathcal O_C)\rightarrow \SHom(\pi^\ast\Omega^1_{Q/k},\mathcal O_C)\rightarrow\mathcal F\rightarrow 0.
\end{align}
We claim that this is the desired sequence: Since $\Omega^1_{Q/k}$ and $\mathcal O_Q$ are locally free $\mathcal O_Q$-modules,
\begin{align*}
\SHom(\pi^\ast\Omega^1_{Q/k},\mathcal O_C)
\cong
\SHom(\pi^\ast\Omega^1_{Q/k},\pi^\ast\mathcal O_Q)
\overset{\footnotemark}\cong\pi^\ast\SHom(\Omega^1_{Q/k},\mathcal O_Q)
=\pi^\ast \T_Q.
\end{align*}
\footnotetext{Locally this means $\Hom_A(M,N)\otimes_A B\cong \Hom_B(M\otimes_A B,N\otimes_A B)$
which holds for free $A$-modules $M,N$.}
Since by definition $\SHom(\Omega^1_{C/k},\mathcal O_C)= \T_C$, (\ref{diag:tangentSes2}) can indeed by identified as
(\ref{diag:tangentSES}).

We now look for the description at the level of stalks:
The stalk of $\pi^\ast \Omega^1_{Q/k}\rightarrow \Omega^1_{C/k}$ at $c$ is by Lemma~\ref{lem:differentialsOfQ}
given by $\varphi : dt\mapsto dt=eu^{e-1}du$.

For the dualized map $\T_C\rightarrow\pi^\ast \T_Q$,
recall that the stalk of $\T_C$ at $c$ is given by $\mathcal O_{C,c}\delta^u$,
and the stalk of $\pi^\ast \T_Q$ at $c$ is given by
$\mathcal O_{C,c}\otimes_{\mathcal O_{Q,\pi(c)}}\mathcal O_{Q,\pi(c)}\delta^t$,
so it is a free $\mathcal O_{C,c}$-module with basis $1\otimes\delta^t$,
which we identify with $\delta^t$.

The dual map of $\varphi$ sends $\delta^u$ to $\delta^u\circ\varphi$,
where $\delta^u\circ\varphi(dt) = \delta^u(eu^{e-1}du) = eu^{e-1}\delta^t(dt)$.
Hence the stalk  of $\T_C\rightarrow\pi^\ast \T_Q$ at $c$ is the map which sends $\delta^u$ to $eu^{e-1}\delta^t$.

We conclude that the stalk of (\ref{diag:tangentSES}) at $c$ is given by
\begin{align*}
0\rightarrow\mathcal O_{C,c}\delta^u\rightarrow\mathcal O_{C,c}\delta^t\rightarrow\mathcal F_c\rightarrow 0.
\end{align*}
Then the description of $\mathcal F_c$ is immediate. In particular, $\mathcal F_c=0$ if $e=1$.
\end{proof}
\begin{lemma}\label{lem:tangentInvariance}
Consider the following $G$-actions:
\begin{itemize}
\item $G$ acts on $\T_C=\SHom(\Omega^1_{C/k},\mathcal O_C)$
by the corresponding actions on $\Omega^1_{C/k}$ resp.\ $\mathcal O_C$ (cf.\ Lemma~\ref{lem:abstractDifferentialAction}).
\item $G$ acts on $\pi^\ast \T_Q$ as in Lemma~\ref{lem:invariantsInverseImage}.
\item $G$ acts on $\mathcal F$, which is a direct sum of skyscraper sheaves,
by sending, for $g\in G$ and $c\in C$ with $u,t,e$ as in Lemma~\ref{lem:tangentSES}
\begin{align*}
x\delta^t+\mathcal O_{C,c}u^{e-1}\delta^t \mapsto (gx)\partial^{gt}+\mathcal O_{C,gc}(gu)^{e-1}\delta^{gt}
\end{align*}
where $g:\mathcal O_{C,c}\xrightarrow\cong \mathcal O_{C,\pi(c)}$ is the induced isomorphism.
\end{itemize}
Then the short exact sequence (\ref{diag:tangentSES}) consists of $G$-equivariant morphisms.
\end{lemma}
\begin{proof}
Well-definedness of the $G$-actions is easy to verify.
To check $G$-invariance, inspect the sequence locally (as described in Lemma~\ref{lem:tangentSES})
and proceed similar to Lemma~\ref{lem:abstractDifferentialAction}.
\end{proof}
\begin{lemma}\label{lem:FCoinvariants}
Let $q_1,\dotsc,q_m\in Q$ be those points of $Q$
which are of the form $\pi(c)$ with $c\in C$ having $G_c\neq\{1\}$.
Let $\mathcal I\subseteq \mathcal O_Q$ be the sheaf of ideals
defining the reduced closed subscheme $\{q_1,\dotsc,q_m\}\subseteq Q$
and $\mathcal F' := \frac{\T_Q}{\mathcal I\T_Q}$.
\begin{enumerate}[(a)]
\item
The stalk of $\mathcal F'$ at each $q_i$ is isomorphic to $k$, and the stalk at each $q\in Q\setminus\{q_1,\dotsc,q_m\}$
is zero.

Pick a $k$-basis element $\partial^i\in(\mathcal F')_{q_i}$ for each $i=1,\dotsc,m$.
\item
Consider the $G$-action on $\H^0(C,\mathcal F)$ as in Lemma~\ref{lem:tangentInvariance}.
There is an isomorphism $\H^0(C,\mathcal F)_G\rightarrow \H^0(Q,\mathcal F')$
which can be described as follows:

For $c\in C$ with $G_c\neq 0$, let $e,u,t$ be as in Lemma~\ref{lem:tangentSES}.
Any element in $\mathcal F_c$ can be written as a sum of elements of the
form $\lambda u^a\delta^t$ for $\lambda\in k$ and $0\leq a<e-1$.
Moreover, let $i\in \{1,\dotsc,m\}$ be the index such that $q_i=\pi(c)$.

Then $\H^0(C,\mathcal F)_G\rightarrow \H^0(Q,\mathcal F')$ sends
\begin{align*}
\lambda u^a \delta^t \mapsto \begin{cases}\lambda \partial^i,&a=0\\
0,&a\neq 0.
\end{cases}
\end{align*}
\item
It is possible to choose the basis elements $\partial^i$ in (a) such that there is a commutative diagram
\begin{align*}
\begin{tikzcd}[ampersand replacement=\&]
\H^0(C,\pi^\ast \T_Q)_G\ar[r]\ar[d,"\cong"]\&\H^0(C,\mathcal F)_G\ar[d,"\cong"]\\
\H^0(Q,\T_Q)\ar[r]\&\H^0(Q,\mathcal F')
\end{tikzcd}
\end{align*}
where the top morphism is $\H^0(C,\blank)_G$ applied to (\ref{diag:tangentSES}),
the left morphism comes from Lemma~\ref{lem:invariantsInverseImage}, the right morphism comes from (b)
and the bottom morphism comes from the construction of $\mathcal F'$ as a quotient of $\T_Q$.
\end{enumerate}
\end{lemma}
\begin{proof}
\begin{enumerate}[(a)]
\item Just note that $\T_Q$ is locally free, that the stalk of $\mathcal I$ at $q_i$ is the maximal
ideal of $\mathcal O_{Q,q_i}$ and the stalk of $\mathcal I$ at each point $q\in Q\setminus\{q_1,\dotsc,q_m\}$
is $\mathcal O_{Q,q}$.
\item Pick preimages $c_i\in C$ such that $\pi(c_i)=q_i$ for $i=1,\dotsc,m$.
Each $c\in C$ with $G_c\neq\{1\}$ is in the orbit of precisely one $c_i$.
Thus
\begin{align*}
\H^0(C,\mathcal F)_G\cong \Bigl(\bigoplus_{\substack{c\in C\\G_c\neq\{1\}}}\mathcal F_c\Bigr)_G
\cong\bigoplus_{i=1}^m (\mathcal F_{c_i})_{G_{c_i}}.
\end{align*}
Let $i\in\{1,\dotsc,m\}$ and $g\in G_{c_i}$ be a generator. Let $u_i,t_i,e_i$ be as in Lemma~\ref{lem:tangentSES}.
We have $gu_i\equiv \zeta u_i\pmod{u_i^2}$ for some primitive $e_i$-th root of unity $\zeta$.
As $g(dt_i)=dt_i$, we conclude $g(\delta^{t_i})=\delta^{t_i}$.
Hence $g(\lambda u_i^a \delta^{t_i})=\zeta^a \lambda u_i^a \delta^{t_i}$.
Therefore $\lambda u_i^a \delta^{t_i}$ vanishes when taking $G_{c_i}$-coinvariants
iff $a\neq 0$ (as $a<e-1<e$).
Thus $(F_{c_i})_{G_{c_i}}\cong k \delta^t$ and the remaining statements of (b) follow.
\item Let $\sigma\in \H^0(C,\pi^\ast \T_Q)$. Its image in $\H^0(C,\mathcal F)$
is obtained by restricting $\sigma$ to $(\pi^\ast \T_Q)_{c}$ for those $c\in C$ with $G_c\neq\{1\}$
and then applying $(\pi^\ast \T_Q)_c\rightarrow \mathcal F_c$.

Let $\tilde\sigma\in \H^0(Q,\T_Q)$ be the corresponding global section by Lemma~\ref{lem:invariantsInverseImage}.
By the compatibility with restriction in Lemma~\ref{lem:invariantsInverseImage},
the stalk of $\tilde\sigma$ at $\pi(c)$ is the stalk of $\sigma$ at $c$
(using $(\T_Q)_{\pi(c)}\hookrightarrow \mathcal O_{C,c}\otimes_{\mathcal O_{Q,\pi(c)}} (\T_Q)_{\pi(c)} = (\pi^\ast \T_Q)_c$
which is an inclusion as $\mathcal O_{Q,\pi(c)}\hookrightarrow \mathcal O_{C,c}$
and $(\T_Q)_{\pi(c)}$ is flat). Hence if $\partial^i$ is chosen to be the image of $\partial^{t_i}\in (\T_Q)_{q_i}$,
the $\delta^{t_i}$ in $(\T_Q)_{q_i}$ and the $\delta^{t_i}$ in $\mathcal F_{c_i}$
have the same image in $\mathcal F'_{q_i}$.
This shows commutativity.
\qedhere\end{enumerate}
\end{proof}
\begin{proposition}\label{prop:cohomologyRootStack}
Let $\mathcal F',m$ be as in Lemma~\ref{lem:FCoinvariants}, and denote by $r$ the rank
of the $k$-linear map $\H^0(Q,\T_Q)\rightarrow \H^0(\mathcal F')$.
Then
\begin{align*}
\dim \HH^1([C/G]) =&\dim\H^0(Q,\T_Q)+\dim \H^0(Q,\Omega^1_{Q/k})-r.
\\\dim \HH^2([C/G]) =&\,m-r+\dim \H^1(Q,\T_Q).
\end{align*}
\end{proposition}
\begin{proof}
We use Proposition~\ref{prop:cohomologyCurveQuotient}. First note that
$\H^1(C,\mathcal O_C)_G$ is, by Serre duality and the fact that group invariants and coinvariants are isomorphic in our situation,
isomorphic to $\H^0(C,\Omega^1_{C/k})^G$. In Proposition~\ref{prop:HH1},
we computed $\H^0(C,\Omega^1_{C/k})^G$ to be isomorphic to $\H^0(Q,\Omega^1_{Q/k})$.

Now the short exact sequence of $\mathcal O_C$-modules (\ref{diag:tangentSES})
induces a long exact sequence
\begin{align*}
0\rightarrow \H^0(C,\T_C)\rightarrow\H^0(C,\pi^\ast \T_Q)\rightarrow \H^0(C,\mathcal F)\rightarrow
\H^1(C,\T_C)\rightarrow \H^1(C,\pi^\ast \T_Q)\rightarrow\H^1(C,\mathcal F).
\end{align*}
Note that $\mathcal F$ has finite support, so there exists an affine open subset
$U\subseteq C$ such that $\mathcal F$ is supported on $U$ (Lemma~\ref{lem:invariantAffineCovering}).
Thus $\H^1(C,\mathcal F)=\H^1(U,\mathcal F)=0$.

Now taking coinvariants in the above exact sequence yields an exact sequence (as taking coinvariants is exact)
\begin{align*}
0\rightarrow \H^0(C,\T_C)_G\rightarrow \H^0(C,\pi^\ast \T_Q)_G\rightarrow \H^0(C,\mathcal F)_G
\rightarrow \H^1(C,\T_C)_G\rightarrow \H^1(C,\pi^\ast \T_Q)_G\rightarrow 0
\end{align*}
By Lemma~\ref{lem:invariantsInverseImage}, $\H^\bullet(C,\pi^\ast \T_Q)_G = \H^\bullet(Q,\T_Q)$.
Moreover, $\dim \H^0(C,\mathcal F)_G = m$ by Lemma~\ref{lem:FCoinvariants}.
By Lemma~\ref{lem:FCoinvariants}, we moreover see that
the rank of $\H^0(C,\pi^\ast \T_Q)_G\rightarrow\H^0(C,\mathcal F)_G$ is $r$.

Using exactness and linear algebra, we thus calculate
\begin{align*}
\dim \H^0(C,\T_C)_G =&\dim \im\left(\H^0(C,\T_C)_G\rightarrow \H^0(C,\pi^\ast \T_Q)_G\right)
\\=&\dim \ker\left(\H^0(C,\pi^\ast \T_Q)_G\rightarrow \H^0(C,\mathcal F)_G\right)
\\=&\dim \H^0(C,\pi^\ast \T_Q)_G - \dim\im\left(\H^0(C,\pi^\ast \T_Q)_G\rightarrow \H^0(C,\mathcal F)_G\right)
\\=&\dim \H^0(Q,\T_Q)-r.\\
\dim \H^1(C,\T_C)_G =&\dim\ker\left(\H^1(C,\T_C)_G\rightarrow \H^1(C,\pi^\ast \T_Q)_G\right)
\\=&\dim \im\left(\H^0(C,\mathcal F)_G\rightarrow \H^1(C,\T_C)_G\right)
+\dim \H^1(C,\pi^\ast \T_Q)_G
\\=&\dim \H^0(C,\mathcal F)_G - \dim\ker\left(\H^0(C,\mathcal F)_G \rightarrow \H^1(C,\T_C)_G\right) +\dim \H^1(Q,\T_Q)
\\=&m-\dim\im\left(\H^0(C,\pi^\ast \T_Q)_G\rightarrow \H^0(C,\mathcal F)_G\right)+\dim \H^1(Q,\T_Q)
\\=&m-r+\dim \H^1(Q,\T_Q).
\end{align*}
This finishes the proof.
\end{proof}
In the above Proposition, everything already depends only on $Q$ and the orders of stabilizers.
We can simplify the situation a bit, expressing everything using the genus of $Q$:
\begin{theorem}\label{thm:summaryQuotient}
With $C,G,Q$ as above, denote
\begin{itemize}
\item $g_Q := \dim \H^0(Q,\Omega^1_{Q/k})$ the genus of $Q$.
\item $q_1,\dotsc,q_m$ the points of $Q$ which admit a preimage $c\in C$ such that $G_c\neq\{1\}$.
\item $e_i := \# G_c$ for $c\in C$ a preimage of $q_i$.
\item $e := e_1+\cdots + e_m$.
\item $a := \begin{cases}3,&g_Q=0\\
1,&g_Q=1\\
0,&g_Q\geq 2\end{cases}$.
\item $r := \min\{m,a\}$.
\item $d := \begin{cases}0,&g_Q=0\\
1,&g_Q=1\\
3g_Q-3,&g_Q\geq 2\end{cases}$.
\end{itemize}
Then $a = \dim \H^0(Q,\T_Q)$, $r$ is the same number as in Proposition~\ref{prop:cohomologyRootStack},
$d=\dim\H^1(Q,\T_Q)$. The Hochschild (co)homology of $[C/G]$ in degree $i\in\mathbb Z$ is given by
\begin{align*}
\dim\HH_i([C/G]) =&\begin{cases}2+e-m,&i=0\\
g_Q,&i=\pm 1\\
0,~&\text{otherwise}\end{cases}\\
\dim\HH^i([C/G]) =&\begin{cases}1,&i=0\\
a+g_Q-r,&i=1\\
m-r +d,&i=2\\
0,~&\text{otherwise}.\end{cases}
\end{align*}
\end{theorem}
\begin{proof}
The calculation of $\dim\HH_{\pm 1}([C/G])$ is Proposition~\ref{prop:HH1}.
The other Hochschild homology groups are calculated in Proposition~\ref{prop:homologyCurveQuotient}.
The claims on $\HH^n([C/G])$ for $n\neq 1,2$ are in Proposition~\ref{prop:cohomologyCurveQuotient}.

Believing for a moment the claims on $a,d$ and $r$,
the claims for $\HH^1([C/G])$ and $\HH^2([C/G])$ follow from Proposition~\ref{prop:cohomologyRootStack}.

By Serre duality, $d=\dim \H^0(Q,\Omega^1_{Q/k}\otimes \Omega^1_{Q/k})$.

Note that if we can show $a = \dim\H^0(Q,\T_Q)$, we get at least
\begin{align*}
\rk(\H^0(Q,\T_Q)\rightarrow \H^0(Q,\mathcal F'))
\leq\min\{\dim \H^0(Q,\T_Q),\dim \H^0(Q,\mathcal F)\}
=\min\{a,m\}=r.
\end{align*}
It remains to show $r\leq \rk(\H^0(Q,\T_Q)\rightarrow \H^0(Q,\mathcal F'))$,
$a=\dim\H^0(Q,\T_Q)$ and $d=\dim\H^1(Q,\T_Q)$.
In order to prove this, we consider the cases $g_Q=0, g_Q=1$ and $g_Q\geq 2$ separately.
\begin{itemize}
\item Case $g_Q=0$: Then $Q\cong \mathbb P^1$ by \cite[Example~IV.1.3.5]{Hartshorne1977}.
By \cite[II.8.20.1]{Hartshorne1977}, $\Omega^1_{Q/k}\cong \mathcal O_Q(-2)$.
Hence $\T_Q\cong \mathcal O_Q(2)$. Writing $\mathbb P^1 = \Proj k[X,Y]$
with the usual grading, $\H^0(Q,\T_Q)$ corresponds to the space of homogeneous degree $2$ polynomials
in $k[X,Y]$. A basis for this space is given by $X^2,Y^2,XY$.
We conclude that $a = 3 = \dim\H^0(Q,\T_Q)$.

We next have to show
$\rk(\H^0(Q,\T_Q)\rightarrow \H^0(Q,\mathcal F'))\geq r$.
Since the automorphism group of $\mathbb P^1$ acts triply transitive,
we may choose our coordinates so that $q_1 = [1:0]$ (or $m=0$),
$q_2 = [0:1]$ (or $m\leq 1$) and $q_3 = [1:1]$ (or $m\leq 2$).

The homogeneous degree $2$-polynomials
\begin{align*}
f_1 = X^2 - XY, f_2 = Y^2-XY, f_3 = XY
\end{align*}
satisfy the following property: For $i,j\in\{1,\dotsc,\min\{3,m\}\}$,  $f_i$ vanishes in $q_j$ if $i\neq j$, but $f_i$ does not vanish in $q_i$.
Hence the images of the $f_i$ for $1\leq i\leq \min\{3,m\}$ (seeing $f_i$ as global section of $\T_Q$)
in $\H^0(Q,\mathcal F')$ are linearly independent.
This shows that the rank of $\H^0(Q,\T_Q)\rightarrow \H^0(Q,\mathcal F')$ is at least $\min\{3,m\}=r$.

Finally, we want to show $\dim\H^0(Q,\Omega^1_{Q/k}\otimes \Omega^1_{Q/k})=d=0$:
Indeed, $\Omega^1_{Q/k}\otimes\Omega^1_{Q/k} \cong\mathcal O_{Q}(-4)$, which has no global sections.
\item Case $g_Q=1$, i.e.\ $Q$ is an elliptic curve. By \cite[Example~IV.1.3.6]{Hartshorne1977},
$\Omega^1_{Q/k}\cong\mathcal O_Q$. Hence $\T_Q\cong\mathcal O_Q$.
Thus $\dim\H^0(Q,\T_Q) = 1 = a$.

If $m=0$, we get $0\leq \rk(\H^0(Q,\T_Q)\rightarrow \H^0(Q,\mathcal F'))\leq r = 0$, so that equality holds.

Now suppose that $m>0$, so that $r=1$. We show that the map $\H^0(Q,\T_Q)\rightarrow \H^0(Q,\mathcal F')$
is non-zero. Let $f\in \H^0(Q,\T_Q)$ be the global section corresponding to 
$1\in k$ under the isomorphism $k\cong \H^0(Q,\mathcal O_Q)\cong \H^0(Q,\mathcal \T_Q)$. We want to show
that the image of $f$ in the stalk $(\T_Q)_{q_1}$ is not in $\mathfrak m_{q_1}(\T_Q)_{q_1}$.
This follows since the image of $1$ in the stalk $\mathcal O_{Q,q_1}$ is not in $\mathfrak m_{q_1}$.

Finally, we want to show $\dim\H^0(Q,\Omega^1_{Q/k}\otimes \Omega^1_{Q/k})=d=1$:
This follows since $\Omega^1_{Q/k}\cong \mathcal O_Q$, so that
$\Omega^1_{Q/k}\otimes \Omega^1_{Q/k}\cong\mathcal O_Q$.
\item Case $g_Q\geq 2$: Let $K_Q$ be the canonical divisor of $Q$. By \cite[Lemma~IV.1.2]{Hartshorne1977}, we must have
$l(-K_Q) = 0$, as otherwise $\deg (-K_Q)\geq 0$, hence $\deg(K_Q)\leq 0$ but by \cite[Example~IV.1.3.3]{Hartshorne1977},
$\deg(K)=2g_Q-2>0$.

The line bundle corresponding to the divisor $-K_Q$ is $\T_Q$, hence $\dim\H^0(Q,\T_Q)=0=a$.
We trivially conclude $\rk(\H^0(Q,\T_Q)\rightarrow \H^0(Q,\mathcal F'))=0$.

Finally, we want to calculate $\dim \H^0(Q,\Omega^1_{Q/k}\otimes \Omega^1_{Q/k}) = l(2K_Q)$.
By Riemann-Roch \cite[Theorem~IV.1.3]{Hartshorne1977} and $l(-K_Q)=0$, we have
\begin{align*}
l(2K_Q) =&\deg(2K_Q) + 1-g_Q = 2\deg K_Q + 1-g_Q\\ =&2(2g_Q-2)+1-g_Q = 3-3g_Q=d,
\end{align*}
proving the last claim.
\end{itemize}
In any case, the stated formulae hold, so the proof is complete.
\end{proof}
\begin{remark}
The genus $g_C$ of $C$ and $g_Q$ satisfy Hurwitz' formula \cite[Corollary~IV.2.4]{Hartshorne1977}:
\begin{align*}
2 g_C-2 = (\# G)(2g_Q-2) + e-m.
\end{align*}
Therefore it satisfies to know $g_C$, $\# G$ and the orders of point stabilizers in $C$ in order to calculate everything in Theorem~\ref{thm:summaryQuotient}.
\end{remark}
\begin{remark}
In the notation of Theorem~\ref{thm:summaryQuotient},
$[C/G]$ is isomorphic to the iterated root stack obtained
from $Q$ by attaching the weight $e_i$ to $q_i$ for $i=1,\dotsc,m$.

We see that the Hochschild (co)homology of this iterated root stack depends only
on $g_Q$, $m$ and $e_1+\cdots + e_m$.

Hence we gave a description of the Hochschild (co)homology
of a weighted smooth projective curve \emph{provided} it can be realized
as a quotient stack $[C/G]$. This last assumption is satisfied in most cases:
\end{remark}
\begin{proposition}[{\cite[Corollary~9]{Lenzing2016}}]\label{prop:abstractLenzing}
Let $X$ be a smooth irreducible projective curve,
$x_1,\dotsc,x_m\in X$ pairwise distinct closed points and $e_1,\dotsc,e_m\in\mathbb Z_{\geq 1}$
be integers. Let $\mathcal X$ be the iterated root stack obtained from $X$ by attaching the weight $e_i$
to the point $x_i$ for $i=1,\dotsc,m$.

Assume $\mathcal X$ is not isomorphic to a weighted projective line of the form
$\mathbb P^1\langle y_1,y_2; a_1,a_2\rangle$ for two distinct weights
$a_1,a_2\in\mathbb Z_{\geq 1}$. Then there exist
a smooth irreducible projective curve $C$ and a finite subgroup $G\leq \Aut(C)$ such that
$\mathcal X\cong [C/G]$ as stacks.\rightqed
\end{proposition}
\begin{proof}[Proof of Theorem~\ref{thm:summaryRootStack}]
In case $\mathcal X\cong\mathbb P^1\langle y_1,y_2; a_1,a_2\rangle$
for some points $y_1,y_2$ and distinct weights $a_1,a_2$,
the theorem follows from Example~\ref{example:calculationWPLFromAbstractResult}
by \emph{using} Proposition~\ref{prop:cohomologyCanonical2}.

In all other cases, we may write $\mathcal X\cong[C/G]$ for a smooth irreducible projective curve
$C$ and a finite group of $k$-scheme automorphisms $G\leq \Aut(C)$ by Proposition~\ref{prop:abstractLenzing}.
Then the result follows from Theorem~\ref{thm:summaryQuotient}.
\end{proof}
\begin{remark}
We now give a second answer why, in the situation of Theorem~\ref{thm:summaryRootStack},
the Hochschild cohomology (and $\HH_{\pm 1}$) of the iterated root stack $\mathcal X$ is determined by the genus of the curve $X$
and the number of exceptional points (in those cases where Proposition~\ref{prop:abstractLenzing} applies to $\mathcal X$).
Namely by Propositions \ref{prop:cohomologyCurveQuotient} and \ref{prop:homologyCurveQuotient}, the Hochschild cohomology
and $\HH_{\pm 1}$ are determined by
$\H^\bullet(C,T_C)_G$ and $\H^0(C,\Omega^1_{C/k})_G$. While $C$ and $\H^0(C,\Omega^1_{C/k})$ may depend on the exceptional points and their weights,
$\H^0(C,\Omega^1_{C/k})_G$ depends only on $X$. The reason for this is the vanishing of $\H^0(C,\Omega^1_{C/Q})^G$ in the proof of Proposition~\ref{prop:HH1} (note that $X=Q$
in our situation).

We calculated $\H^\bullet(C,T_C)_G$ using the short exact sequence from Lemma~\ref{lem:tangentSES}.
Here, the information about the exceptional points is encoded in the sheaf $\mathcal F$ and the information
about the original curve $X$ is encoded in $\pi^\ast T_Q$. $\H^\bullet(C,\mathcal F)$ depends
only on the weights of the exceptional points but not their location in $X$, and the morphism $\H^0(C,\mathcal F)\rightarrow \H^1(C,T_C)$
a priori depends on their locations in $X$. However, we saw that $\H^\bullet(C,\mathcal F)_G$ depends only on the number of exceptional points
and, using that the automorphism group is sufficiently transitive in the proof of Theorem~\ref{thm:summaryQuotient},
the map $\H^0(C,\mathcal F)\rightarrow \H^1(C,T_C)$ has always maximal rank. This explains how, in the end, the Hochschild cohomology is determined
by the genus of $X$ and the number of exceptional points alone.
\end{remark}
\section{Realizing a class of weighted projective lines}\label{sec:projectiveCompanion}
Consider the weighted projective line $X$ associated
to the points $x_1,\dotsc,x_n\in\mathbb P^1$ and the weights $a_1,\dotsc,a_n\in\mathbb Z_{\geq 1}$.
In Proposition~\ref{prop:abstractLenzing}, we saw that $X$ can be expressed as a quotient stack of a smooth projective curve $C$ by a finite group action in most cases.

In this section, we give an explicit construction of such a curve $C$ and a group $G$
for an infinite class of weighted projective lines, and show how to calculate
the Hochschild (co)homology using just Propositions \ref{prop:cohomologyCurveQuotient} and \ref{prop:homologyCurveQuotient}.

We assume again the Convention~\ref{conv:wpl}.
\begin{definition}
Let $S$ be the $\mathbb N_0$-graded $k$-algebra defined as follows:

As algebra, $S = k[X_1,\dotsc,X_n]/I$ where $I$ is the ideal generated by $X_i^{a_i} - X_2^{a_2} + x_i X_1^{a_1}$. Elements of $k$ are defined to be of degree $0$,
and we declare the degree of $X_i$ to be $\delta_i := \frac{\lcm(a_1,\dotsc,a_n)}{a_i}$.
(so we get the same algebra as in Definition~\ref{def:GeigleLenzing} but with different grading group).

We define the $k$-scheme $C := \Proj S$ (cf.~\cite[Proposition~II.2.5]{Hartshorne1977}) and call it the \emph{projective companion} of $X$.
\end{definition}
\begin{remark}
In his paper, Lenzing does not define a projective companion but rather a \emph{twisted companion} $Y$.
$Y$ is a smooth one-dimensional stack whose category of coherent sheaves is given by $\catmod^{\mathbb Z} S / \catmod^{\mathbb Z}_0 S$
with $S$ defined as above \cite[Proposition~11]{Lenzing2016}.
$X$ is always a quotient stack of $Y$ by the finite group $G=\frac{\mu_{a_1}\times\cdots\times \mu_{a_n}}{\mu_{\lcm(a_1,\dotsc,a_n)}}$, where $\mu_{a}$
is the group of $a$-th roots of unity in $k$ (ibid).
Moreover, Lenzing gives suficient conditions, in terms of the weights $a_i$, for $Y$ to be a projective curve, so that $X$ is the quotient of an explicit projective curve by an explicit group.

The approach with the projective companion developed here is somewhat orthogonal to this: $C$ is always a projective curve which comes with a natural morphism to $\mathbb P^1$,
expressing $\mathbb P^1$ as a quotient scheme of $C$ by a natural action of $G=\frac{\mu_{a_1}\times\cdots\times \mu_{a_n}}{\mu_{\lcm(a_1,\dotsc,a_n)}}$. Moroever, under certain conditions on the weights $a_i$,
$C$ is smooth and the quotient stack $[C/G]$ is indeed the weighted projective line $X$.

For our purposes, the projective companion is more suitable for two reasons: First, an explicit description of $C$ as scheme simplifies the calculations
of the cohomology of the (co)tangent sheaves, as used in Propositions \ref{prop:cohomologyCurveQuotient} and \ref{prop:homologyCurveQuotient}.
Second, the conditions which we need for $C$ to be smooth and for $[C/G]$ to be $X$ are weaker than the conditions
Lenzing requires in \cite[Proposition~11]{Lenzing2016} for $Y$ to be a scheme. Thus the usage of the twisted companion allows
to cover more examples of weighted projective lines.
\end{remark}
\begin{lemma}[Properties of the algebra $S$]\label{lem:propertiesS}
\begin{enumerate}[(a)]
\item
$S$ is a two-dimensional noetherian domain.
\item Let $k[X,Y]$ denote the polynomial ring with grading $\deg X := \deg Y := \lcm(a_1,\dotsc,a_n)$.
Then the graded ring map $\Phi : k[X,Y]\rightarrow S, X\mapsto X_1^{a_1}, Y\mapsto X_2^{a_2}$
is a finite ring extension.
\end{enumerate}
\end{lemma}
\begin{proof}
\begin{enumerate}[(a)]
\item
By\footnote{In the cited paper, the $L(\mathbf a)$-graded algebra is considered; however,
the results cited here are independent of the grading.} \cite[Proposition~1.3]{Geigle1987}, $S$ is a two-dimensional domain.
\item It is an extension as it is easy to check that no polynomial in $k[X_1,X_2]\setminus\{0\}$ is in $I$.
The extension is integral since $X_1,X_2$ and $X_i = X_2^{a_2} - x_i X_1^{a_1}$ are integral over $k[X,Y]$.
It is an extension of finite type as $S$ is of finite type even over $k$. Hence the extension is finite.
\qedhere\end{enumerate}
\end{proof}
\begin{lemma}[Standard affine covering]\label{lem:standardAffine}
Consider the open subsets $U = D(X_1)$ and $V = D(X_2)$ of $C$.
\begin{enumerate}[(a)]
\item $U$ and $V$ are affine, and $C=U\cup V$.
Moreover\footnote{$S_{(f)}$ denotes the degree zero part of the localization $S_f$
for homogeneous $f\in S_+$.}, $U\cong \Spec S_{(X_1)}$, $V\cong \Spec S_{(X_2)}$ and $U\cap V\cong \Spec S_{(X_1X_2)}$.
Under these isomorphisms,
the inclusions $U\cap V\hookrightarrow U$ resp.\ $U\cap V\hookrightarrow V$ are $\Spec$ applied to the inclusions
\begin{align*}
S_{(X_1)}\hookrightarrow S_{(X_1X_2)},\qquad
S_{(X_2)}\hookrightarrow S_{(X_1X_2)}.
\end{align*}
\item
The finite ring extension $\Phi : k[X,Y]\rightarrow S$ induces a finite and surjective morphism of $k$-schemes
$C\rightarrow\mathbb P^1$. In particular, $\dim C=1$.
\end{enumerate}
\end{lemma}
\begin{proof}
\begin{enumerate}[(a)]
\item The isomorphisms
$U\cong \Spec S_{(X_1)}$, $V\cong \Spec S_{(X_2)}$ and $U\cap V\cong \Spec S_{(X_1X_2)}$
are by \cite[Proposition~2.5]{Hartshorne1977}. In particular, $U$ and $V$ are affine.
Then the identifications of $U\cap V\hookrightarrow U,V$ are immediate.

Let $\mathfrak p\in C$ be a closed point not contained in $U$ nor $V$.
Then $\mathfrak p\subseteq S$ is a prime ideal containing $X_1$ and $X_2$.
Therefore, $\mathfrak p$ contains $X_i^{a_i} = X_2^{a_2} - x_i X_1^{a_1}$ for $i=3,\dotsc,n$,
hence it contains all $X_i$ for $i=1,\dotsc,n$. Thus $S_+\subseteq \mathfrak p$,
contradicting the definition of $C=\Proj S$. The contradiction shows $C=U\cup V$.
\item
The existence of the morphism follows from
 \citestacks{01MY} and $C=D(\Phi(X))\cup D(\Phi(Y))$ by (a).
Finiteness follows from the fact that finite ring morphisms are stable under base change \citestacks{02JK}, hence under localizations.
Surjectivity follows from \citestacks{05DR}.
Then the statement about the dimension of $C$ follows from \citestacks{0ECG}.\qedhere
\end{enumerate}
\end{proof}
\begin{lemma}[Projective embedding]\label{lem:projectiveEmbedding}
\begin{enumerate}[(a)]
\item There exists a natural number $d_0$ such that for all positive multiples $d$ of $d_0$,
the $\mathbb N_0$-graded algebra $S_{(d)}$ defined by
\begin{align*}
S_{(d)} = \bigoplus_{i=0}^\infty S_{id},
\end{align*}
where the degree $i$-homogeneous part of $S_{(d)}$ is $S_{id}$,
is generated as $k$-algebra by $S_d$.
\item For any $d$ as in (a),
there is an isomorphism of $k$-schemes
\begin{align*}
\Proj S\xrightarrow\cong \Proj S_{(d)}
\end{align*}
sending a prime ideal $\mathfrak p\in\Proj S$ to $\mathfrak p\cap S_{(d)}\in \Proj S_{(d)}$.
In particular, $C$ is a projective $k$-scheme.
\end{enumerate}
\end{lemma}
\begin{proof}
(a) is \citestacks{0EGH} and the isomorphism in (b) follows from \citestacks{0B5J} and its proof.
Then the fact that $\Proj S_{(d)}$ is projective is immediate.
\end{proof}
\subsection{Points of $C$}
The aim of this subsection is to describe the closed points of $C$.
\begin{lemma}\label{lem:equivalentPoints}
Let $(\lambda_1,\dotsc,\lambda_n),(\mu_1,\dotsc,\mu_n)\in k^n\setminus\{0\}$ be two vectors. Put $\overline a:=\lcm(a_1,\dotsc,a_n)$. The following are equivalent:
\begin{enumerate}[(a)]
\item There exists some $z\in k^\times$ such that for $i=1,\dotsc,n$:
\begin{align*}
z^{a_i}\lambda_i = \mu_i.
\end{align*}
\item For $i,j\in \{1,\dotsc,n\}$, denote the greatest common divisor of $a_i$ and $a_j$ by $(a_i,a_j)$.
Then
\begin{align*}
\mu_j^{\frac{a_j}{(a_i,a_j)}} \lambda_i^{\frac{a_i}{(a_i,a_j)}} =
\mu_i^{\frac{a_i}{(a_i,a_j)}} \lambda_j^{\frac{a_j}{(a_i,a_j)}}.
\end{align*}
\end{enumerate}
\end{lemma}
\begin{proof}
(a) $\Longrightarrow$ (b): Writing $\mu_i = z^{\overline a/a_i} \lambda_i$, we get
\begin{align*}
\mu_j^{\frac{a_j}{(a_i,a_j)}} \lambda_i^{\frac{a_i}{(a_i,a_j)}} =
\lambda_j^{\frac{a_j}{(a_i,a_j)}} z^{\frac{\overline a}{(a_i,a_j)}} \lambda_i^{\frac{a_i}{(a_i,a_j)}} =
\lambda_j^{\frac{a_j}{(a_i,a_j)}} \mu_i^{\frac{a_i}{(a_i,a_j)}}.
\end{align*}
(b) $\Longrightarrow$ (a):
First suppose that there was some index $i$ for which $\lambda_i\neq 0$ but $\mu_i=0$.
Pick an index $j$ such that $\mu_j\neq 0$ (by choice of $(\mu_1,\dotsc,\mu_n)$).
Then by (b), $\lambda_j^{\frac{a_j}{(a_i,a_j)}}\mu_i^{\frac{a_i}{(a_i,a_j)}}\neq 0$,
contradicting $\mu_i=0$.
By symmetry, we observe $\lambda_i=0\iff\mu_i=0$.

In order to prove (a), we may omit those indices for which $\lambda_i=\mu_i=0$ (because condition (a) holds trivially for them).
Let us therefore assume $\lambda_i\neq 0\neq \mu_i$ for all $i$.

Now condition (b) can be restated as
\begin{align*}
\left(\frac{\lambda_i}{\mu_i}\right)^{\frac{a_i}{(a_i,a_j)}} = 
\left(\frac{\lambda_j}{\mu_j}\right)^{\frac{a_j}{(a_i,a_j)}}\tag{$\ast$}
\end{align*}
for all $i,j$. Pick some $z\in k^\times$ such that $\frac{\mu_1}{\lambda_1} = z^{\overline a/a_1}$ (as $k$ is algebraically closed).
Replacing $(\mu_1,\dotsc,\mu_n)$ by $(z^{\overline a/a_1}\mu_1,\dotsc,z^{\overline a/a_n}\mu_n)$, we may assume $\mu_1=\lambda_1$.
The condition $(\ast)$ for $i=1$ states that $\frac{\lambda_j}{\mu_j}$ is an $\frac{a_j}{(a_1,a_j)}$-th root of unity.

Put $A := (a_1\cdots a_n)^2$ and let $\zeta\in k^\times$ be an $A$-th root of unity.
For $i=1,\dotsc,n$, $\frac{\lambda_i}{\mu_i}$ is an $A$-th root of unity,
so there exists some exponent $e_i\in \mathbb Z$ with $\frac{\lambda_i}{\mu_i} = \zeta^{e_i}$.
Moreover, as $\zeta$ is an $A$-th root of unity and $\zeta^{e_i}=\frac{\lambda_i}{\mu_i}$ is an $a_i$-th root of unity,
$\overline a$ divides $e_i$ by choice of $A$.

We can restate $(\ast)$ as
\begin{align*}
&e_i \frac{a_i}{(a_i,a_j)} \equiv e_j \frac{a_j}{(a_i,a_j)}\pmod A.
\\&\iff \frac{e_i a_i}{\overline a}\equiv \frac{e_j a_j}{\overline a} \pmod{\frac{A(a_i,a_j)}{\overline a}}.
\end{align*}
By elementary number theory (using a general version of the chinese remainder theorem),
this shows that there exists a number $t$, unique modulo $A$,
such that $\frac{e_i a_i}{\overline a}\equiv t\pmod{\frac{A a_i}{\overline a}}$.
Thus $e_i\equiv \frac{\overline a}{a_i} t\pmod A$, and hence
\begin{align*}
\frac{\lambda_i}{\mu_i} = \zeta^{e_i} =\zeta^{t\overline a/a_i}= (\zeta^t)^{\overline a/a_i}.
\end{align*}
We may choose $z=\zeta^{-t}$ to fulfil (a).
\end{proof}
\begin{definition}
\begin{enumerate}[(a)]
\item We call two vectors $(\lambda_1,\dotsc,\lambda_n),(\mu_1,\dotsc,\mu_n)\in k^n\setminus\{0\}$ fulfilling the equivalent conditions
of Lemma~\ref{lem:equivalentPoints} \emph{equivalent}. The equivalence class of $(\lambda_1,\dotsc,\lambda_n)$
will be denoted by $[\lambda_1:\cdots:\lambda_n]$.
\item Define $P$ to be the set of tuples $(\lambda_1,\dotsc,\lambda_n)\in k^n\setminus\{0\}$ such that
\begin{align*}
\lambda_i^{a_i} = \lambda_2^{a_2} - x_i\lambda_1^{a_1}\text{ for }i=3,\dotsc,n
\end{align*}
modulo the equivalence relation from (a).
\end{enumerate}
\end{definition}
The aim of this subsection is to identify $P$ with the set of closed points of $C$.
\begin{lemma}\label{lem:pointsOfC}
Let $i\in\{1,\dotsc,n\}$. Consider the affine open subset $\Spec S_{(X_i)}\cong D(X_i)\subseteq C$,
where $S_{(X_i)}$ is again the degree $0$ part of the localization $S_{X_i}$.
\begin{enumerate}[(a)]
\item The canonical morphism of affine schemes
\begin{align*}
\Spec S_{X_i}\rightarrow \Spec S_{(X_i)} \cong D(X_i)
\end{align*}
is surjective.
\item For all points $[\lambda_1:\cdots:\lambda_n]\in P$,
the ideal $I([\lambda_1:\cdots:\lambda_n])$ generated by those homogeneous elements $f\in S$
which satisfy $f(\lambda_1,\dotsc,\lambda_n)=0$ is a non-zero homogeneous prime ideal of $S$.
It is independent of the equivalence class of $(\lambda_1,\dotsc,\lambda_n)$.
\item Each maximal ideal of $S_{X_i}$ has the form $(X_1-\lambda_1,\dotsc,X_n-\lambda_n)$
for a uniquely determined vector $(\lambda_1,\dotsc,\lambda_n)\in k^n$ where $\lambda_i\neq 0$
and for $j=3,\dotsc,n$:
\begin{align*}
\lambda_j^{a_j} = \lambda_2^{a_2} - x_i \lambda_1^{a_1}.
\end{align*}
\end{enumerate}
\end{lemma}
\begin{proof}
\begin{enumerate}[(a)]
\item The fact $\Spec S_{(X_i)} \cong D(X_i)$ is \cite[Proposition~2.5]{Hartshorne1977}.
Let $I\subseteq \Spec S_{(X_i)}$ be any ideal, and let $S_{X_i} I\subseteq S_{X_i}$
be the ideal of $S_{X_i}$ generated by $I$. We claim that the degree $0$ part of $S{X_i} I$ is $I$, i.e.\ $S_{(X_i)}\cap (S_{X_i} I) = I$:
The inclusion $\supseteq$ is clear. Now let $a\in S_{(X_i)}\cap (S_{X_i} I)$.
As $a\in S_{X_i} I$, we can write
\begin{align*}
a=f_1 i_1+\cdots+f_r i_r
\end{align*}
for $f_1,\dotsc,f_r\in S_{X_i}$ and $i_1,\dotsc,i_r\in I$.
Decomposing the $f_j$, we may assume each $f_j$ is a homogeneous element of $S_{X_i}$.
Considering the degree zero part of our above equation,
\begin{align*}
a=\sum_{\substack{j=1\\\deg(f_j)=0}}^r f_j i_j.
\end{align*}
As $\deg(f_j)=0$ implies $f_j\in S_{(X_i)}$, we see $a\in I$.

Now let $\mathfrak p \in \Spec S_{(X_i)}$. If $\mathfrak p=0$, a preimage is $0\in \Spec S_{X_i}$.
Otherwise $\mathfrak p$ is a maximal ideal of $S_{(X_i)}$ (as $\dim S_{(X_i)}=1$, cf. Lemma~\ref{lem:standardAffine}).
Since $S_{X_i}\mathfrak p$ is a proper ideal of $S_{X_i}$ by the above calculation,
it is contained in a maximal ideal $\mathfrak m$ of $S_{X_i}$, which is our desired preimage.
\item The fact that $I:= I([\lambda_1:\cdots:\lambda_n])$ is a homogeneous prime ideal is immediate
(if $fg\in I$, then $f(\lambda_1,\dotsc,\lambda_n)g(\lambda_1,\dotsc,\lambda_n)=0$,
so $f\in I$ or $g\in I$).
Since $\lambda_i^{a_i} X_j^{a_j} - \lambda_j^{a_j} X_i^{a_i}\in I$ for all $i,j$, we
see $I\neq 0$.

Finally, if $(\lambda_1,\dotsc,\lambda_n) = (z^{\overline a/a_1} \mu_1,\dotsc,z^{\overline a/a_n} \mu_n)$
for some $z\in k^\times$ and $f\in S$ is homogeneous of degree $d$,
$f(\lambda_1,\dotsc,\lambda_n) = z^d f(\mu_1,\dotsc,\mu_n)$, showing independence
of the equivalence class of $\lambda$.
\item
This is a version of Hilbert's Nullstellensatz ($\lambda_i\neq 0$ comes from the localization at $X_i$).
\qedhere\end{enumerate}
\end{proof}
\begin{proposition}\label{prop:weightedNullstellensatz}
There is a one-to-one correspondence between closed points of $C$ and elements of $P$,
given by sending $[\lambda_1:\cdots:\lambda_n]\in P$ to $I([\lambda_1:\cdots:\lambda_n])\in C$.
\end{proposition}
\begin{proof}
The map is well-defined by Lemma~\ref{lem:pointsOfC}.

Injectivity: Suppose that $[\lambda_1:\cdots:\lambda_n],[\mu_1:\cdots:\mu_n]\in P$ describing the same closed point
$I([\lambda_1,\dotsc,\lambda_n]) = I([\mu_1:\cdots:\mu_n])$.
For all $i,j\in \{1,\dotsc,n\}$,
\begin{align*}
\mu_j^{\frac{a_j}{(a_i,a_j)}} X_i^{\frac{a_i}{(a_i,a_j)}} - 
\mu_i^{\frac{a_i}{(a_i,a_j)}} X_j^{\frac{a_j}{(a_i,a_j)}}\in I([\mu_1:\cdots:\mu_n]),
\end{align*}
so this also lies in $I([\lambda_1:\cdots:\lambda_n])$. This means, spelled out, that
\begin{align*}
\mu_j^{\frac{a_j}{(a_i,a_j)}} \lambda_i^{\frac{a_i}{(a_i,a_j)}} =
\mu_i^{\frac{a_i}{(a_i,a_j)}} \lambda_j^{\frac{a_j}{(a_i,a_j)}}.
\end{align*}
By Lemma~\ref{lem:equivalentPoints}, $[\lambda_1:\cdots:\lambda_n] = [\mu_1:\cdots:\mu_n]$.

Surjectivity: Let $\mathfrak p\in C$ be a closed points. Since we can cover $C$ by the $D(X_i)$, there exists some $i\in\{1,\dotsc,n\}$
such that $\mathfrak p \in D(X_i)\cong \Spec S_{(X_i)}$. By Lemma~\ref{lem:pointsOfC},
the image of $\mathfrak p$ in $S_{(X_i)}$ equals
$(X_1-\lambda_1,\dotsc,X_n-\lambda_n)\cap S_{(X_i)}$ for some
$(\lambda_1,\dotsc,\lambda_n)\in k^n\setminus\{0\}$ satisfying $\lambda_j^{a_j} = \lambda_2^{a_2}-x_j \lambda_1^{a_1}$ for all $j\in\{3,\dotsc,n\}$. It follows
that $I([\lambda_1:\cdots:\lambda_n])\subseteq \mathfrak p$.
As $I([\lambda_1:\cdots:\lambda_n])\neq 0$ (Lemma~\ref{lem:pointsOfC}) and $\dim C=1$ (Lemma~\ref{lem:standardAffine}),
$I([\lambda_1:\cdots:\lambda_n]) = \mathfrak p$.
\end{proof}
\subsection{Differentials of $C$}
The aim of this section is to describe the cotangent bundle $\Omega_C^1$
and to prove that $C$ is smooth.
\begin{remark}
As noted before, we will require some condition on the weights $a_1,\dotsc,a_n$ to prove that $C$ is smooth.
It is not clear to the author whether this condition is also necessary. However,
the condition stated here will later be necessary and sufficient in order to obtain the desired
weighted projective line as quotient stack, so we lose no generality imposing the condition already here.

We say that the weight vector $\mathbf a = (a_1,\dotsc,a_n)$ satisfies the \emph{weight condition}
if each weight divides the least common multiple of the other weights, i.e.\
\begin{align*}
a_i \mid \lcm(a_1,\dotsc,a_{i-1},a_{i+1},\dotsc,a_n)\text{ for }i=1,\dotsc,n.
\end{align*}
\end{remark}
\begin{example}
The weight vectors $(2,2), (2,4,4), (2,3,6), (2,3,5,6,10)$ satisfy the weight condition,
whereas $(2,3), (2,4,8)$ and $(2,3,5)$ do not.
\end{example}
\begin{remark}
We claimed before that our weight condition is weaker than the condition required in \cite[Proposition~11]{Lenzing2016}:
There, it is assumed that\footnote{The assumption $n\geq 3$ or $n=2$ and $a_1=a_2$ is not explicitly stated, however,
there are trivial counterexamples for $n=2$:
e.g.\ the weight vector $(2,3)$ seems to satisfy the condition of \cite[Proposition~11]{Lenzing2016}.
However, the statement of that Proposition would be that the corresponding weighted projective
line $X$ is the quotient of the non-weighted projective curve $Y$ by $\frac{\mu_2\times \mu_3}{\mu_6}$, which is the trivial group.
Hence $X\cong Y$ would be non-weighted projective, a contradiction.}
\begin{align*}
\deg X_1 = \frac{\overline a}{a_1},\dotsc,\deg X_n = \frac{\overline a}{a_n}
\end{align*}
are pairwise coprime, and that $n\geq 3$ or $n=2$ and $a_1=a_2$. As $\gcd(\frac{\overline a}{a_i},\frac{\overline a}{a_j}) = \frac{\overline a}{\lcm(a_i,a_j)}$,
this is equivalent to $\overline a$ being the lcm of \emph{any} two weight $a_i,a_j$ with $i\neq j$.
If $n\geq 3$, we can find for each weight $a_h$ some $i,j\in\{1,\dotsc,n\}$ with $h,i,j$ pairwise distinct,
so that \begin{align*}a_h \mid \overline a = \lcm(a_i a_j) \mid \lcm(a_{i'} : i'\neq h),\end{align*}
which is our weight condition.

For $n=3$, our weight condition and Lenzing's are equivalent.
This is not the case for $n\geq 3$, e.g.\ $(2,2,4,4)$ is a weight vector satisfying our weight condition but not Lenzing's.
\end{remark}
\begin{lemma}\label{lem:differentialsOnS}
We define $\tilde \Omega$ as the $S$-module by the following generators and relations:
\begin{align*}
\tilde \Omega := \frac{\bigoplus_{i=1}^n S d X_i}{\left(a_iX_i^{a_i-1} dX_i - a_2 X_2^{a_2-1} dX_2 + x_i a_1 X_1^{a_1-1}\text{ for }i=3,\dotsc,n\right)}.
\end{align*}
We turn $\tilde\Omega$ into a graded $S$-module by the convention $\deg(d X_i) := \deg X_i = \frac{\overline a}{a_i}$.
Define a differential $S\xrightarrow d\tilde \Omega$ by $d(X_i) = dX_i$, extended by the Leibniz rule.

Then $S\xrightarrow d\tilde \Omega$ is a universal derivation, i.e.\ $\tilde \Omega$ is the module of K\"{a}hler differentials.
\end{lemma}
\begin{proof}
Well-definedness of $\tilde \Omega$ and $d$ are easily checked.

By \cite[Example~II.8.2.1]{Hartshorne1977}, the module of K\"{a}hler differentials for
$k[X_1,\dotsc,X_n]$ is $\bigoplus\limits_{i=1}^n k[X_1,\dotsc,X_n] dX_i$ with $d$ similarly as defined above.
Since we can express $S$ as a quotient ring of $k[X_1,\dotsc,X_n]$, the statement follows from \cite[Proposition~II.8.4A]{Hartshorne1977}.
\end{proof}
\begin{lemma}\label{lem:localDifferentialsOnC}
Let $\tilde\Omega$ be as above and
consider the graded morphism of $S$-modules
\begin{align*}
\tilde \Omega \xrightarrow \varphi S, dX_i \mapsto \frac{X_i}{a_i}.
\end{align*}
Denote $\Omega := \ker(\varphi)$.

Let $\mathfrak p$ in $C$ be a closed point with coordinates $[\lambda_1:\cdots:\lambda_n]$.
\begin{enumerate}[(a)]
\item If $w\in S$ is homogeneous of degree $\delta$, we have $\varphi(d(w))=\delta\overline a^{-1} w$.
\item The stalk of the map $S\xrightarrow d \tilde\Omega$ at $\mathfrak p$, i.e.\ $S_{(\mathfrak p)}\xrightarrow d\tilde \Omega_{(\mathfrak p)}$,
factors over $\Omega_{(\mathfrak p)}$.
\item Assume that the weight condition holds. Then $S_{(\mathfrak p)}\xrightarrow d\Omega_{(\mathfrak p)}$ is the module of K\"{a}hler differentials of $S_{(\mathfrak p)}$.
\item Assume that the weight condition holds. Choose indices $i,j\in\{1,\dotsc,n\}$ such that $\lambda_{i'}\neq 0$ for all $i'\neq i$.
There exists a unit $f\in S_{\mathfrak p}$ of degree $\deg(f)=-\frac{\overline a}{a_i} - \frac{\overline a}{a_j}$.
For any such unit, $\Omega_{(\mathfrak p)}$ is a free $S_{(\mathfrak p)}$-module of rank $1$ with generator
\begin{align*}
f(a_i X_j dX_i - a_j X_i dX_j).
\end{align*}
\end{enumerate}
\end{lemma}
\begin{proof}
\begin{enumerate}[(a)]
\item We may assume $w$ has the form $w=X_1^{e_1}\cdots X_n^{e_n}$ with $e_i\geq 0$. Then
\begin{align*}
\varphi(d(w)) =&\varphi\left(\sum_{i=1}^n e_i X_1^{e_1}\cdots X_{i-1}^{e_{i-1}} X_i^{e_i-1} X_{i+1}^{e_{i+1}}\cdots X_n^{e_n} dX_i\right)
=\sum_{i=1}^n \frac{e_i}{a_i} X_e^{e_1}\cdots X_n^{e_n}
\\=&w\sum_{i=1}^n \frac{e_i}{a_i}
=w\sum_{i=1}^n e_i \deg X_i\overline a^{-1} = w\overline a^{-1}\sum_{i=1}^n \deg(X_i^{e_i})
=w \overline a^{-1}\delta.
\end{align*}
\item Let $\frac fg\in S_{(\mathfrak p)}$ with $f,g\in S$ homogeneous (so $g$ does not vanish in $[\lambda_1:\cdots:\lambda_n]$).
Denote $\delta := \deg(f)=\deg(g)$.
Then
\begin{align*}
\varphi(d(f/g)) = \varphi\left(\frac{d(f)}{g} - \frac{f d(g)}{g^2}\right)
=\frac{\varphi(d(f))}g - \frac{f\varphi(d(g))}{g^2}
\underset{\text{(a)}}=\frac{\delta f}{\overline a g} - \frac{f \delta g}{\overline a g^2}=0.
\end{align*}
\item By compatibility with localization, $S_{\mathfrak p}\xrightarrow d\tilde\Omega_{\mathfrak p}$ is the
module of K\"{a}hler differentials of $S_{\mathfrak p}$.

Note that $\gcd(\deg(X_1),\dotsc,\deg(X_n)) = \frac{\overline a}{\lcm(a_1,\dotsc,a_n)} = 1$.
Moreover, by the weight condition, even $\gcd(\deg(X_j):j\neq i)=1$ for all $i=1,\dotsc,n$.
Since at most one $X_i$ can vanish in $[\lambda_1:\cdots:\lambda_n]$,
we have $\gcd(\deg(X_i): X_i\text{ invertible in }S_{\mathfrak p})=1$.

Hence $S_{\mathfrak p}$ contains a unit $t$ of degree $1$.
It follows that $S_{\mathfrak p}\cong S_{(\mathfrak p)}[t,t^{-1}]$ as graded algebras, where
the latter is, by abuse of notation, the ring of Laurent polynomials over $S_{(\mathfrak p)}$ with $\deg(t)=1$.

Note that we have a (non-graded) $k$-algebra homomorphism $S_{(\mathfrak p)}[t,t^{-1}]\rightarrow S_{(\mathfrak p)}$ sending $t$ to $1$.
Then the composition
\begin{align*}
S_{(\mathfrak p)}\hookrightarrow S_{(\mathfrak p)}[t,t^{-1}]\rightarrow S_{(\mathfrak p)}
\end{align*}
is the identity $k$-algebra homomorphism. It follows that the induced map of K\"{a}hler differentials
\begin{align*}
\Omega_{S_{(\mathfrak p)/k}}^1\rightarrow 
\Omega_{S_{(\mathfrak p)[t,t^{-1}]/k}}^1\rightarrow 
\Omega_{S_{(\mathfrak p)/k}}^1
\end{align*}
is the identity map. In particular, the natural map
$\Omega_{S_{(\mathfrak p)}/k}^1\rightarrow 
\Omega_{S_{(\mathfrak p)}[t,t^{-1}]/k}^1\cong \tilde\Omega_{\mathfrak p}$ is injective.
Thus $\Omega_{S_{(\mathfrak p)}/k}^1$ is the $S_{(\mathfrak p)}$-submodule
of $\tilde\Omega_{\mathfrak p}$ generated by the image of $S_{(\mathfrak p)}\xrightarrow d\tilde\Omega_{\mathfrak p}$.
It suffices to show that this image is $\ker(\varphi)\cap \tilde\Omega_{(\mathfrak p)}$.

For $e\in\mathbb Z$ define
\begin{align*}
U_\delta := \{ w\in \tilde \Omega_{\mathfrak p} : \varphi(w) = \delta \overline a^{-1} w\}.
\end{align*}
By (a), this is a graded $S_{\mathfrak p}$-submodule of $\tilde\Omega_{\mathfrak p}$
and $\tilde\Omega_{\mathfrak p}=\bigoplus\limits_{\delta\in \mathbb Z} U_l\delta$.

Note that $\tilde\Omega_{\mathfrak p}$ is, as $S_{\mathfrak p}$-module, generated by
the elements $d(s t^e)$ for $s\in S_{(\mathfrak p)}$ and $e\in \mathbb Z$.
Then $d(s t^e)\in U_e$ by the above calculation. Hence $U_0$ is generated, as $S_{\mathfrak p}$-module,
by $d(s)$ for $s\in S_{(\mathfrak p)}$. Therefore the degree $0$ part of $U_0$ is generated, as $S_{(\mathfrak p)}$-module,
by $d(s)$ for $s\in S_{(\mathfrak p)}$. This shows that $\Omega_{S_{(\mathfrak p)}/k}^1$ is indeed the degree $0$ part of $U_0$, as claimed.
\item The existence of units of arbitrary degree was settled above.

It is simple to see that $s f (a_i X_j dX_i - a_j X_i dX_j) = 0$ for $s\in S_{\mathfrak p}$ implies $s=0$
by the definition of $\tilde \Omega$ (using that $S_{\mathfrak p}$ is a domain).

Finally let $w\in \Omega_{(\mathfrak p)}$ be any element. Write $w=f_1 dX_1 + \cdots + f_n dX_n$
for homogeneous elements $f_i \in S_{\mathfrak p}$ of degree $\deg f_i = -\deg X_i = -\frac{\overline a}{a_i}$.

Using the conditions
\begin{align*}
a_{i'} X_{i'}^{a_{i'}-1} d X_{i'} = a_2 X_2^{a_2-1} dX_2 - x_{i'} X_1^{a_1-1} dX_1,
\end{align*}
for $i'\geq 3$, we claim that it is possible to express each $dX_{\ell}$ using $dX_i$ and $dX_j$, i.e.\
\begin{align*}\Omega_{(\mathfrak p)} = (S_{\mathfrak p})_{-\deg(X_i)} dX_i + (S_{\mathfrak p})_{-\deg(X_j)} dX_j\end{align*}
The claim is proved by case distinction:
\begin{itemize}
\item In case $\{i,j\}=\{1,2\}$, this is clear: Use the condition for $i'=\ell$ and invertibility of $X_{\ell}$ to get $dX_{\ell}\in (S_{\mathfrak p})_{-\deg(X_1)} dX_1
+ (S_{\mathfrak p})_{-\deg (X_2)}dX_2$.
\item In case $1\in \{i,j\}$ but $2\notin \{i,j\}$, use the condition for $i'=\max\{i,j\}$ and invertibility of $X_2$ to get $dX_2$.
Then use $dX_1$ and $dX_2$ and the remining conditions to get all $X_\ell$ for $\ell\geq 3$.
\item In case $2\in\{i,j\}$ but $1\notin\{i,j\}$, use the condition for $i'=\max\{i,j\}$ and invertibility of $X_1$ to get $dX_1$. Proceed as above.
\item In case $\{i,j\}\cap \{1,2\}\neq \emptyset$, subtract the condition for $i'=i$ from the condition for $i'=j$. Thus
\begin{align*}
a_j X_j^{a_j-1} dX_j - a_i X_i^{a_i-1} dX_i = (x_i-x_j) X_1^{a_1-1} dX_1.
\end{align*}
Now use the invertibility of $X_1$ and $x_i\neq x_j$ to get $dX_1$. Proceed as above.
\end{itemize}
In any case, $\Omega_{(\mathfrak p)} = (S_{\mathfrak p})_{-\deg(X_i)} dX_i + (S_{\mathfrak p})_{-\deg(X_j)} dX_j$.
Now if $w=f_i dX_i + f_j dX_j\in \Omega_{(\mathfrak p)}$, we have
\begin{align*}
&0=\varphi(w) = \frac{f_i X_i}{a_i} + \frac{f_j X_j}{a_j}\Longrightarrow f_j = -\frac{a_j X_i}{a_i X_j} f_i
\\&\Longrightarrow w = \frac{f_i}{a_i X_j}\left(a_i X_j dX_i - a_j X_i dX_j\right)
=\frac{f_i}{a_i X_jf} f\left(a_i X_j dX_i - a_j X_i dX_j\right).
\end{align*}
As $\frac{f_i}{a_i X_jf}$ is of degree $0$, we see $w\in S_{(\mathfrak p)}f\left(a_i X_j dX_i - a_j X_i dX_j\right)$.\qedhere
\end{enumerate}
\end{proof}
\begin{corollary}\label{cor:differentialsOnC}
Assume that the weight condition holds.
\begin{enumerate}[(a)]
\item $C$ is smooth.
\item Denote $\omega := (n-2)\overline a - \sum\limits_{i=1}^n \frac{\overline a}{a_i}$. For $i\in \{2,\dotsc,n\}$
and $f = X_2\cdots X_{i-1}X_{i+1}\cdots X_n$, consider the element
\begin{align*}
w_i := \pm (a_1 X_i dX_1 - a_i X_1 dX_i) \prod_{\substack{j=2\\j\neq i}}^n X_j^{1-a_j} \in (\Omega_{f})_{-\omega},
\end{align*}
where the sign is $+1$ for $i\geq 3$ and $-1$ for $i=2$.
Then
\begin{align*}
(S_f)_{\omega}\rightarrow \Omega_{(f)}, g\mapsto g w_i
\end{align*}
is an isomorphism of $S_{(f)}$-modules.
\item For $i\neq i'$, the restrictions of $w_i$ and $w_{i'}$ to $D(X_2\cdots X_n)$ coincide. In particular,
the above isomorphisms glue to an isomorphism of line bundles from $\mathcal O_C(\omega)$ to the sheaf of modules
on $C$ associated to the graded $S$-module $\Omega$, i.e.\ the cotangent sheaf of $C$.
\end{enumerate}
\end{corollary}
\begin{proof}
\begin{enumerate}[(a)]
\item The cotangent sheaf is the sheaf associated to the graded $S$-module $\Omega$.
As $\Omega_{(\mathfrak p)}$ is free of rank one for each closed point $\mathfrak p$,
$C$ is regular. Since we are in characteristic zero, $C$ is already smooth.
\item As $\rho$ and the $X_j$ for $i\neq j\geq 2$ are units, it suffices to show that the following is an isomorphism:
\begin{align*}
(S_f)_{\omega}\rightarrow \Omega_{(f)},g\mapsto g(a_1 X_i dX_1 - a_i X_1 dX_i).
\end{align*}
This can be checked locally at each closed point $\mathfrak p\in D(X_2\cdots X_{i-1} X_{i+2}\cdots X_n)$.
Note that if $f\in (S_{\mathfrak p})_{-\deg(X_1)-\deg(X_i)}$ is a unit,
$(S_{\mathfrak p})_{-\deg(X_1)-\deg(X_i)}$ is a free $S_{(\mathfrak p)}$-module of rank $1$. Now the
claim follows by Lemma \ref{lem:localDifferentialsOnC} part (d).
\item It suffices to show that the restrictions
of each $w_i$ and $w_2$ coincide for $i\geq 3$.
For that notice that if all $X_j$ are invertible for $j\geq 2$,
\begin{align*}
&a_1 X_i dX_1 - a_i X_1 dX_i =
a_1 X_i dX_1 - \frac{X_1}{X_i^{a_i-1}} (a_i X_i^{a_i-1} dX_i)
\\=&
a_1 X_i dX_1 - \frac{X_1}{X_i^{a_i-1}} (a_2 X_2^{a_2-1} dX_2 - x_i a_1 X_1^{a_1-1} dX_1)
\\=&a_1 X_i \left(1 - \frac{x_i X_1^{a_1}}{X_i^{a_i}}\right) dX_1 - \frac{X_2^{a_2-1}}{X_i^{a_i-1}} a_2 X_1 dX_2
\\=&a_i X_i\left(1-\frac{X_i^{a_i} - X_2^{a_2}}{X_i^{a_i}}\right) dX_1- \frac{X_2^{a_2-1}}{X_i^{a_i-1}} a_2 X_1 dX_2
\\=&-\frac{X_2^{a_2-1}}{X_i^{a_i-1}}\left(a_1 X_2 dX_1 - a_2 X_1 dX_2\right).
\end{align*}
Multiplying both sides by $\prod\limits_{\substack{j=2\\j\neq i}}^n X_j^{1-a_j}$, we get the claim.
\qedhere\end{enumerate}
\end{proof}
\begin{remark}\label{rem:genus}Assume that the weight condition holds.

It follows that the genus of $C$ is $\dim_k \H^0(C,\mathcal O_C(\omega)) =\dim_k S_\omega$.
The genus can be computed as follows:

We first express $\mathcal O_C(\omega)$ via a Cartier divisor: Consider the open cover
\begin{align*}
\{U_i = D\left(X_1\cdots X_{i-1}X_{i+1}\cdots X_n\right) \mid i=1,\dotsc,n\}.
\end{align*}
The weight condition implies $\gcd(\deg(X_1),\dotsc,\deg(X_{i-1}),\deg(X_{i+1}),\dotsc,\deg(X_n))=1$.
We therefore may choose indices $e_{i,1},\dotsc,e_{i,n}\in \mathbb Z$
such that
\begin{align*}
e_{i,i} = \begin{cases}-1,&i\in \{1,2\}\\
a_i-1,&i\geq 3\end{cases}
\end{align*}
and
$f_i := X_1^{e_{i,1}}\cdots X_n^{e_{i,n}}\in \Frac(S)$ has degree $0$.
Then $\{(U_i,f_i)\}_{i=1,\dotsc,n}$ defines a Cartier divisor as in \cite[Sec.~II.6]{Hartshorne1977}.

Let $K := S_{(0)}$ be the field of rational functions and denote, by abuse of notation,
the sheaf on $C$ which is constantly $K$ also by $K$. Consider the injective homomorphism of sheaves
\begin{align*}
\mathcal O_C(\omega)\rightarrow K
\end{align*}
given by multiplying a section by $v := X_1 X_2 \prod\limits_{i=3}^n X_i^{1-a_i}$ (which has degree $-\omega$).
Denote the image of this morphism by $\mathcal L$.

We claim that $\mathcal L$ is the line bundle $\mathcal L(D)$ associated to the Cartier divisor\begin{align*}D:=\{(U_i,f_i)\}_{i=1,\dotsc,n}.\end{align*}

i.e.\ for each $i=1,\dotsc,n$ and $\mathfrak p \in U_i$, we would like
to show that the $S_{(\mathfrak p)}$-submodule of $K$ generated by $f_i^{-1}$ equals the submodule
generated by $v (S_{\mathfrak p})_\omega$.
This is clear if $\mathfrak p \in D(X_1\cdots X_n)$, as all the $X_i$ become
units and the submodule in each case is just $\mathcal O_{C,\mathfrak p}=S_{(\mathfrak p)}$.
If $\mathfrak p \in V(X_i)$, the submodule is
\begin{align*}
\begin{cases} X_i S_{(\mathfrak p)},&i=1,2\\
X_i^{1-a_i}S_{(\mathfrak p)},&i\geq 3.
\end{cases}\tag{$\ast$}
\end{align*}
This shows that $\mathcal L(D)\cong \mathcal L \cong \mathcal O_C(\omega)\cong \Omega_{C/k}^1$.
In other words, $D$ is the canonical divisor.
By \cite[Example~IV.1.3.3]{Hartshorne1977}, the degree of this canonical divisor is $2g-2$, where $g$ is the genus of $C$.
Note that by $(\ast)$, the Weyl divisor associated to $D$ is
\begin{align*}
\sum_{i=1}^n \sum_{\mathfrak p \in V(X_i)} \mathfrak p \cdot\begin{cases} 1,&i=1,2\\
1-a_i,&i\geq 3.\end{cases}
\end{align*}
Hence the degree of $D$ is
\begin{align*}
\sum_{i=1}^n \# V(X_i)\begin{cases} 1,&i=1,2\\
1-a_i,&i\geq 3.\end{cases}
\end{align*}
We will later see that $\# V(X_i) = \frac{a_1\cdots a_{i-1}a_{i+1}\cdots a_n}{\overline a}$.
It follows that the degree of $D$ is
\begin{align*}
&\frac{a_1\cdots a_n}{\overline a} \left(\frac 1{a_1} + \frac 1{a_2} + \sum_{i=3}^n \frac 1{a_i} - 1\right)
=\frac{a_1\cdots a_n}{\overline a} \left( (2-n) + \sum_{i=1}^n \frac 1{a_i}\right)
\\=&-\frac{a_1\cdots a_n}{\overline a^2} \omega.
\end{align*}
We conclude that the genus of $C$ is
\begin{align*}
1-\frac 12\left(-\frac{a_1\cdots a_n}{\overline a^2} \omega\right) = 1+\frac{a_1\cdots a_n}{2\overline a^2}\omega.
\end{align*}
\end{remark}
\begin{example}\label{ex:genusZero}
Assume again the weight condition.
By \cite[Example~IV.1.3.5]{Hartshorne1977}, $C$ is isomorphic to $\mathbb P^1$
iff the genus is zero.

If the genus of $C$ is zero, $\frac{a_1\cdots a_n\omega}{2\overline a^2}=-1$. In particular, $\omega < 0$.
Thus $n-2 - \sum\limits_{i=1}^n \frac 1{a_i}<0$.

Conversely, if $n-2 - \sum\limits_{i=1}^n \frac 1{a_i}<0$, we get $\omega<0$ so the genus of $C$
is $\dim S_\omega=0$.
\begin{itemize}
\item Case $n=2$: The weight condition implies $a_1=a_2$, so $\mathbf a=(a_1,a_1)$.
\item Case $n\geq 3$: Then $a_i\geq 2$ implies $\sum\limits_{i=1}^n \frac 1{a_i} \leq \frac n2$, so that $\frac n2-2\geq n-2-\sum\limits_{i=1}^n \frac 1{a_i}<0$.
We see $n<4$.

Now let $n=3$. Then $0>n-2-\sum\limits_{i=1}^n a_i = 1-\sum\limits_{i=1}^3 \frac 1{a_i}$.
The only solutions to this are, up to reordering,
\begin{align*}
\mathbf a \in \{(2,2,a) : a\geq 2\}\cup \{(2,3,3), (2,3,4), (2,3,5)\}.
\end{align*}
Only the weight vector $\mathbf a = (2,2,2)$ satisfies the weight condition.
\end{itemize}
\end{example}
\begin{example}\label{ex:genusOne}
Similarly, $C$ is elliptic (i.e.\ the genus is $1$) iff $\omega = 0$.
By a similar case distinction as above, this is the case for
\begin{align*}
\mathbf a \in \{ (2,3,6), (2,4,4), (3,3,3), (2,2,2,2)\}.
\end{align*}
All of these satisfy the weight condition.
\end{example}
\subsection{The group action on $C$}
As announced before, the natural map $C\rightarrow\mathbb P^1$ from Lemma~\ref{lem:standardAffine} expresses
$\mathbb P^1$ as a quotient scheme of $C$ by a finite group. In this section, we make the group and its action explicit.
\begin{definition}\label{def:groupOnC}
For $a\in \mathbb N$, let $\mu_a \leq k^\times$ be the group of $a$-th roots of unity in $k$.

Consider the action of $\mu_{a_1}\times \cdots \times \mu_{a_n}$ on $S$ by the algebra automorphism
\begin{align*}
(\zeta_1,\dotsc,\zeta_n)\cdot X_i := \zeta_i X_i
\end{align*}
for $i=1,\dotsc,n$. Denote the induced map on $C=\Proj S$ by $\varphi_{(\zeta_1,\dotsc,\zeta_n)}$.
Denote the image of the group homomorphism
\begin{align*}
\mu_{a_1}\times\cdots\times \mu_{a_n}\xrightarrow{\varphi_{(\blank)}} \Aut(C)
\end{align*}
by $G$.
\end{definition}
\begin{lemma}\label{lem:groupActionOnPoints}~
\begin{enumerate}[(a)]
\item If $\mathfrak p\in C$ is a closed points with coordinates $[\lambda_1:\cdots:\lambda_n]$,
$\varphi_{(\zeta_1,\dotsc,\zeta_n)}(\mathfrak p)$ is the closed point in $C$
with coordinates $[\zeta_1^{-1}\lambda_1:\cdots:\zeta_n^{-1}\lambda_n]$.
\item For $i=1,\dotsc,n$, consider the homomorphism $\mu_{\overline a}\rightarrow \mu_{a_i}$,
sending $\zeta$ to $\zeta^{\overline a/a_i}$. Then the induced homomorphism
$\mu_{\overline a}\rightarrow \mu_{a_1}\times\cdots\times \mu_{a_n}$ and
the group homomorphism from Definition~\ref{def:groupOnC} form a short exact sequence
\begin{align*}
0\rightarrow \mu_{\overline a}\rightarrow \mu_{a_1}\times\cdots\times \mu_{a_n}\rightarrow G\rightarrow 0.
\end{align*}
In particular, $G\cong \frac{\mu_{a_1}\times\cdots\times \mu_{a_n}}{\mu_{\overline a}}$.

\item The natural map $C\rightarrow \mathbb P^1$ from Lemma~\ref{lem:standardAffine}
expresses $\mathbb P^1$ as a quotient scheme $C/G$.
\item If $\mathfrak p\in C$ is a closed point with coordinates $[\lambda_1:\cdots:\lambda_n]$
such that $\lambda_i\neq 0$ for all $i$, the stabilizer of $\mathfrak p$ in $G$ is trivial.
\item Let $i\in\{1,\dotsc,n\}$. There exists a unique $G$-orbit of closed points $\mathfrak p\in C$
with coordinates $[\lambda_1:\cdots:\lambda_n]$ such that $\lambda_i=0$.

The image of each such point in $\mathbb P^1$ has the coordinate $x_i$ in $\mathbb P^1$.
\item Assume that the weight condition holds and let $\mathfrak p$ be one of the closed points from (e).
Then the stablilizer of $\mathfrak p$ in $G$ is cyclic of order $a_i$.
\end{enumerate}
\end{lemma}
\begin{proof}
\begin{enumerate}[(a)]
\item Let $f\in S$. Then we have
\begin{align*}
f\in \mathfrak p\iff&f(\lambda_1,\dotsc,\lambda_n)=0\iff \varphi_{(\zeta_1,\dotsc,\zeta_n)}(f) (\zeta_1^{-1}(\lambda_1),\dotsc,\zeta_n^{-1}(\lambda_n)) = 0
\\\iff&\varphi_{(\zeta_1,\dotsc,\zeta_n)} \in \mathfrak p'
\end{align*}
where $\mathfrak p'$ is the closed point of $C$ with coordinates $[\zeta_1^{-1}\lambda_1:\cdots:\zeta_n^{-1}\lambda_n]$.
This shows that $\varphi_{(\zeta_1,\dotsc,\zeta_n)}(\mathfrak p)=\mathfrak p'$.
\item Let $(\zeta_1,\dotsc,\zeta_n)\in \mu_{a_1}\times\cdots\times \mu_{a_n}$.
This element is in the kernel of $\mu_{a_1}\times\cdots\times \mu_{a_n}\rightarrow G$
iff for all closed points $\mathfrak p\in C$,
$\varphi_{(\zeta_1,\dotsc,\zeta_n)}(\mathfrak p)=\mathfrak p$.

Write $\mathfrak p$ with coordinates $[\lambda_1:\cdots:\lambda_n]$.
Then
\begin{align*}
\varphi_{(\zeta_1,\dotsc,\zeta_n)}([\lambda_1:\cdots:\lambda_n]) = [\zeta_1^{-1}\lambda_1:\cdots:\zeta_n^{-1}\lambda_n].
\end{align*}
This equals $[\lambda_1:\cdots:\lambda_n]$ iff there exists some $t\in k^\times$ such that $\zeta_i^{-1}\lambda_i = t^{\overline a/a_i}\lambda_i$ for all $i$.

On the one hand, if $(\zeta_1,\dotsc,\zeta_n)$ has a preimage $\zeta\in\mu_{\overline a}$, we may choose $t=\zeta^{-1}$ for all points
$\mathfrak p$ to establish $\varphi_{(\zeta_1,\dotsc,\zeta_n)}([\lambda_1:\cdots:\lambda_n]) = [\lambda_1:\cdots:\lambda_n]$.

Now conversely suppose that $\varphi_{(\zeta_1,\dotsc,\zeta_n)}\in \Aut(C)$ is the trivial map and pick a point $\mathfrak p=[\lambda_1:\cdots:\lambda_n]\in D(X_1\cdots X_n)$.
Then there exists some $t\in k^\times$ such that $\zeta_i^{-1} = t^{\overline a/a_i}$ for all $i$, as $\lambda_i\neq 0$.
As $t^{\overline a} = \zeta_i^{-a_i} = 1$, $t$ is an $\overline a$-th root of unity.
Now $t^{-1}\in \mu_{\overline a}$ is a preimage of $(\zeta_1,\dotsc,\zeta_n)$ in $\mu_{\overline a}$.

The fact that $\mu_{\overline a}\rightarrow \mu_{a_1}\times\cdots\times \mu_{a_n}$ is injective
follows easily by elementary number theory, using only $\overline a=\lcm(a_1,\dotsc,a_n)$.
\item It suffices to show that two closed points $\mathfrak p,\mathfrak p'\in C$ have the same image in $\mathbb P^1$
iff they are in the same $G$-orbit (surjectivity of $C\rightarrow\mathbb P^1$ was proved in Lemma~\ref{lem:standardAffine}).

In coordinates, write $\mathfrak p=[\lambda_1:\cdots:\lambda_n]$ and $\mathfrak p' = [\lambda_1':\cdots:\lambda_n']$.
Then the respective images in $\mathbb P^1$ have coordinates $\frac{\lambda_2^{a_2}}{\lambda_1^{a_1}}$ respectively
$\frac{(\lambda_2')^{a_2}}{(\lambda_1')^{a_1}}$ (with the convention $1/0=\infty$).

Suppose first that $\mathfrak p' = \varphi_{(\zeta_1,\dotsc,\zeta_n)}(\mathfrak p)$. Then
$\lambda_i' = \zeta_i\lambda_i$ such that $(\lambda_i')^{a_i} = \lambda_i^{a_i}$ for all $i$.
It follows that $\mathfrak p$ and $\mathfrak p'$ have the same image in $\mathbb P^1$.

Now assume conversely that the images of $\mathfrak p$ and $\mathfrak p'$ in $\mathbb P^1$ coincide.
Hence there exists some $t\in k^\times$ such that $(\lambda_1')^{a_1} = t \lambda_1^{a_1}$
and $(\lambda_2')^{a_2} = t \lambda_2^{a_2}$.

Pick some $s\in k$ such that $s^{\overline a}=t$. Put $\tilde{\lambda_i} := s^{\overline a/a_i}\lambda_i'$.
Then $[\lambda_1':\cdots:\lambda_n'] = [\tilde{\lambda_i},\dotsc,\tilde{\lambda_n}]$
and $\tilde{\lambda_1}^{a_1} = \lambda_1^{a_1}$, $\tilde{\lambda_2}^{a_2} = \lambda_2^{a_2}$.

By the equation
\begin{align*}
\lambda_i^{a_i} = \lambda_2^{a_2} - x_i \lambda_1^{a_1}\text{ for }i=3,\dotsc,n,
\end{align*}
we conclude $\tilde{\lambda_i}^{a_i} = \lambda_i^{a_i}$ for all $i$. Therefore, there
exists an $a_i$-th root of unity $\zeta_i\in \mu_{a_i}$ such that $\tilde{\lambda_i} = \zeta_i\lambda_i$.
We conclude
\begin{align*}
\varphi_{(\zeta_1,\dotsc,\zeta_n)}([\lambda_1:\cdots:\lambda_n]) = [\tilde{\lambda_1}:\cdots:\tilde{\lambda_n}] = [\lambda_1':\cdots:\lambda_n'].
\end{align*}
Hence $\varphi_{(\zeta_1,\dotsc,\zeta_n)}(\mathfrak p)=\mathfrak p'$.
\item We saw this proof already in (b).
\item We distingusih the cases $i=1$ and $i\geq 2$ for formal reasons, although they are no actually different:

Case $i=1$: Existence of at least one point $\mathfrak p=[\lambda_1:\cdots:\lambda_n]$ with $\lambda_1=0$:
Set $\lambda_1=0,\lambda_2=1$ and choose, for $i\geq 3$, $\lambda_i$ to satisfy $\lambda_i^{a_i} = \lambda_2^{a_2}-x_i \lambda_i^{a_i}$.

Now if $\mathfrak p=[\lambda_1:\cdots:\lambda_n]$ is any closed point of $C$ satisfying $\lambda_1=0$, it follows that
$\lambda_2\neq 0$ and the image of $\mathfrak p$ in $\mathbb P^1$ is $\frac{\lambda_2^{a_2}}{\lambda_1^{a_1}} = \frac{\lambda_2^{a_2}}0=\infty=x_1$.

Therefore, any two points in $\mathfrak p$ with vanishing $X_1$-coordinate have the same image in $\mathbb P^1$,
so they are $G$-conjugate by (c).

Case $i\geq 2$: Existence of at least one point $\mathfrak p=[\lambda_1:\cdots:\lambda_n]$ with $\lambda_i=0$:
Set $\lambda_1=1$ and $\lambda_i=0$. If $i\geq 3$, choose $\lambda_2$ to satisfy $\lambda_i^{a_i} = \lambda_2^{a_2}-x_i \lambda_1^{a_1}$.
Now for $i\neq j\geq 3$, choose $\lambda_j$ to satisfy $\lambda_j^{a_j} = \lambda_2^{a_2} - x_j \lambda_1^{a_1}$.

Next suppose $\mathfrak p=[\lambda_1:\cdots:\lambda_n]$ is any closed point in $C$ satisfying $\lambda_i=0$.
Note that $0=\lambda_i^{a_i} = \lambda_2^{a_2}-x_i \lambda_1^{a_1}$ holds in any case (for $i=2$, recall $x_2=0$).
Hence the image of $\mathfrak p$ in $\mathbb P^1$ has coordinate \begin{align*}\frac{\lambda_2^{a_2}}{\lambda_1^{a_1}}
=\frac{x_i\lambda_1^{a_1}}{\lambda_1^{a_1}}=x_i.\end{align*}

Again, if $\mathfrak p$ and $\mathfrak p'$ are two points with vanishing $X_i$-coordinate, they have the same image in $\mathbb P^1$,
so they are $G$-conjugate.
\item Let $\mathfrak p=[\lambda_1:\cdots:\lambda_n]$. We claim that the stabilizer of $\mathfrak p$ is the image
of the composition \begin{align*}\mu_{a_i}\hookrightarrow \mu_{a_1}\times\cdots\times \mu_{a_n}\rightarrow G.\end{align*}

On the one hand, an element $(\zeta_1,\dotsc,\zeta_n)$ with $\zeta_j=1$ for all $j\neq i$
will certainly stabilize $\mathfrak p$ as $\zeta_j\lambda_j = \lambda_j$ holds for all $j$ (since $\lambda_i=0$).

Now suppose conversely that $(\zeta_1,\dotsc,\zeta_n)\in \mu_1\times\cdots\times \mu_n$ stabilizes $\mathfrak p$.
By definition, $[\zeta_1^{-1}\lambda_1:\cdots:\zeta_n^{-1}\lambda_n] = [\lambda_1:\cdots:\lambda_n]$,
i.e.\ $\zeta_j^{-1}\lambda_j = t^{\overline a/a_j} \lambda_j$ for some $t\in k^\times$ and all $j$.
Then $t$ is an $\overline a$-th root of unity as before and $t^{\overline a/a_j}\zeta_j=1$ if $i\neq j$. Thus
\begin{align*}
\varphi_{(\zeta_1,\dotsc,\zeta_n)} = \varphi_{t^{\overline a/a_1}\zeta_1,\dotsc,t^{\overline a/a_n} \zeta_n}
=\varphi_{(1,\dotsc,1,t^{\overline a/a_i} \zeta_i,1,\dotsc,1)}
\end{align*}
has a preimage in $\mu_{a_i}$.

This proves that the stabilizer of $\mathfrak p$ is the image of $\mu_{a_i}$ in $G$.
We still have to show that the map $\mu_{a_i}\rightarrow G$ is injective.
Suppose $\zeta_i\in \mu_{a_i}$ is in the kernel. By (a), there exists some $\zeta\in \mu_{\overline a}$
with $\zeta^{\overline a/a_i}=a_i$ and $\zeta^{\overline a/a_j}=1$ for $j\neq i$.

By the weight condition, the induced map 
\begin{align*}
\mu_{\overline a}\rightarrow\mu_{a_1}\times\cdots\times\mu_{a_{i-1}}\times \mu_{a_{i+1}}\times\cdots\times \mu_{a_n}
\end{align*}
still is injective, so $\zeta=1$ and hence $\zeta_i=1$.

Now $\mu_{a_i}\hookrightarrow G$ is the stabilizer of $\mathfrak p$, which is cyclic of order $a_i$.
\qedhere\end{enumerate}
\end{proof}
\begin{corollary}\label{cor:realizingWPL}
Assume that the weight condition holds. Then $[C/G]$ is, as stack, isomorphic
to the weighted projective line $\mathbb P^1\langle x_1,\dotsc,x_n; a_1,\dotsc,a_n\rangle$.
\end{corollary}
\begin{proof}
This follows from Lemma~\ref{lem:groupActionOnPoints} and Proposition~\ref{prop:quotientRootStack}.
\end{proof}
\begin{remark}
By the proof of Lemma~\ref{lem:groupActionOnPoints}, part (f), it is easy to see that the existence of an isomorphism
$[C/G]\cong \mathbb P^1\langle x_1,\dotsc,x_n; a_1,\dotsc,a_n\rangle$ is equivalent
to the weight condition (if the weight condition fails, some stabilizers will not be large enough).
e.g.\ the weight vector $(2,3)$ fails the weight condition,
and indeed, $G\cong \frac{\mu_2\times \mu_3}{\mu_6}\cong\{1\}$,
hence no point in $C$ has a non-trivial stabilizer.
\end{remark}
\begin{remark}[Continuing Remark~\ref{rem:genus}]
To finish the argument of Remark~\ref{rem:genus}, we have to show
that the number of points $\mathfrak p =[\lambda_1:\cdots:\lambda_n]$ such that $\lambda_i=0$
is $\frac{a_1\cdots a_{i-1}a_{i+1}\cdots a_n}{\overline a} = \frac{\abs G}{a_i}$.
This follows, assuming the weight condition, directly from Lemma~\ref{lem:groupActionOnPoints} and the orbit-stabilizer-formula.

Moreover, we note that the Euler characteristic of $C$ is $2-2g=\abs G\frac{\omega}{\overline a}$, where
$g$ is its genus.
By \cite[Theorem~2]{Lenzing2016}, we hence would expect the Euler characteristic
of $[C/G] \cong \mathbb P^1\langle x_1,\dotsc,x_n; a_1,\dotsc,a_n\rangle$
to be $\frac{\omega}{\overline a}$. And indeed, this is the Euler characteristic
by \cite[58]{Voight2019}.
\end{remark}
\begin{example}
In the terminology of \cite[Theorem~7]{Lenzing2016}, we say that the weighted projective line $\mathbb P^1\langle x_1,\dotsc,x_n; a_1,\dotsc,a_n\rangle$
is of parabolic type iff $\mathbf a \in\{(2,3,6), (2,4,4), (3,3,3), (2,2,2,2)\}$.
These are precisely the cases where $C$ is elliptic by Example~\ref{ex:genusOne}.
\end{example}
\begin{lemma}[Group action of differentials]\label{lem:actionOnDifferentials}
Let $(\zeta_1,\dotsc,\zeta_n)\in \mu_1\times\cdots\times\mu_n$ and $g$ be the corresponding group element $g:=\varphi_{(\zeta_1,\dotsc,\zeta_n)}\in G$.
Denote by $\mathcal O_C(\omega)\xrightarrow \psi\Omega^1_{C/k}$ be the isomorphism
of line bundles from Corollary~\ref{cor:differentialsOnC}.

As $g$ is a scheme automorphism, consider the induced isomorphism
$g_\ast\Omega_{C/k}^1\cong \Omega_{C/k}^1$.
Moreover, the action of $(\zeta_1,\dotsc,\zeta_n)$ on $S$
induces an isomorphism $g_\ast \mathcal O_C(\omega)\cong \mathcal O_C(\omega)$.

Then the following diagram of sheaf isomorphisms commutes:
\begin{align*}
\begin{tikzcd}[ampersand replacement=\&]
\mathcal O_C(\omega)\ar[r,"\psi"]\ar[d,"\cdot(\zeta_1\cdots\zeta_n)^{-1}"]\&\Omega^1_{C/k}\ar[r,"\cong"]\&g_\ast\Omega^1_{C/k}\\
\mathcal O_C(\omega)\ar[r,"\cong"]\&g_\ast\mathcal O_C(\omega)\ar[ru,"g_\ast\psi"]
\end{tikzcd}.
\end{align*}
Here, the leftmost map is multiplication of sections by $\frac 1{\zeta_1\cdots\zeta_n}\in k^\times$.
\end{lemma}
\begin{proof}
It suffices to check commutativity on an open covering of $C$.
Let $i\in\{2,\dotsc,n\}$ and consider the open subset $D(f)$ for $f=X_2\cdots X_{i-1}X_{i+1}\cdots X_i$.
Note that $D(f)$ is $G$-invariant. Hence $\H^0(D(f),g_\ast M)=\H^0(D(f),M)$ for all sheaves of $\mathcal O_C$-modules
$M$.

We get the following commutative diagram from \ref{cor:differentialsOnC}:
\begin{align*}
\begin{tikzcd}[ampersand replacement=\&]
\&S_{(f)}\ar[r,hook]\ar[d,"d_{S_{(f)}/k}"]\&S_f\ar[d,"d_{S_f/k}"]\\
\H^0(D(f),\mathcal O_C(\omega)) = (S_f)_{\omega}\ar[r,"w_i\cdot"]\&\Omega_{(f)}\ar[r,hook]\&\Omega_{S_f/k}^1
\end{tikzcd}.
\end{align*}
Now consider the induced action of $\underline \zeta = (\zeta_1,\dotsc,\zeta_n)$:
First, note that the following diagram commutes:
\begin{align*}
\begin{tikzcd}[ampersand replacement=\&,column sep=5em]
S_f\ar[r,"{X_i\mapsto \zeta_i X_i}"]\ar[d,"d"]\& S_f\ar[d,"d"]\\
\Omega_{S_f/k}^1\ar[r,"{sdX_i\mapsto (\zeta\cdot s)\zeta_i dX_i}"]\&\Omega_{S_f/k}^1
\end{tikzcd}.
\end{align*}
By restriction, also the following diagram commutes:
\begin{align*}
\begin{tikzcd}[ampersand replacement=\&,column sep=7em]
S_{(f)}\ar[r,"{X_i\mapsto \zeta_i X_i}"]\ar[d,"d"]\& S_{(f)}\ar[d,"d"]\\
\Omega_{(f)}\ar[r,"{sdX_i\mapsto (\zeta\cdot s)\zeta_i dX_i}"]\&\Omega_{(f)}
\end{tikzcd}.
\end{align*}
Note that the top isomorphism $S_{(f)}\rightarrow S_{(f)}$ is $\H^0(D(f),\blank)$ applied to $g_\ast\mathcal O_C\xrightarrow \cong\mathcal O_C$.
It follows that the bottom isomorphism $\Omega_{(f)}\xrightarrow{sdX_i\mapsto (\zeta\cdot s)\zeta_i dX_i} \Omega_{(f)}$ must be
$\H^0(D(f),\blank)$ applied to the isomorphism $g_\ast \Omega_{C/k}^1\xrightarrow \cong\Omega_{C/k}^1$.

We now identify
\begin{align*}
\begin{tikzcd}[ampersand replacement=\&,column sep=7.4em]
\H^0(D(f),\mathcal O_C(\omega))\ar[r,"{\H^0(D(f),\psi)}"]\ar[d,"\cdot(\zeta_1\cdots\zeta_n)^{-1}"]
\&\H^0(D(f),\Omega^1_{C/k})\ar[r,"{\H^0(D(f),\Omega_{C/k}^1\cong
g_\ast\Omega_{C/k}^1)}"]\&\H^0(D(f),g_\ast\Omega^1_{C/k})\\
\H^0(D(f),\mathcal O_C(\omega))\ar[r,"{\H^0(D(f),\mathcal O_C\cong g_\ast\mathcal O_C)}"]\&\H^0(D(f),g_\ast\mathcal O_C(\omega))\ar[ru,"{\H^0(D(f),g_\ast\psi)}"']
\end{tikzcd}
\end{align*}
as
\begin{align*}
\begin{tikzcd}[ampersand replacement=\&,column sep=8em]
(S_f)_{\omega}\ar[r,"{w_i\cdot}"]\ar[d,"\cdot (\zeta_1\cdots\zeta_n)^{-1}"]
\&\Omega_{(f)}\ar[r,"{sdX_i\mapsto (\zeta^{-1}.s)\zeta_i^{-1} dX_i}"]\&\Omega_{(f)}\\
(S_f)_{\omega}\ar[r,"{s\mapsto (\zeta^{-1}.s)}"]\&(S_f)_{\omega}\ar[ru,"{w_i\cdot}"']
\end{tikzcd}.
\end{align*}
Now commutativity follows from the observation that $\zeta^{-1}\cdot w_i = (\zeta_1\cdots\zeta_n)^{-1}w_i$
by the explicit formula for $w_i$ from Corollary~\ref{cor:differentialsOnC} and $\zeta_i^{a_i}=1$.
\end{proof}
\begin{corollary}\label{cor:actionOnCodifferentials}
In the setting of Lemma~\ref{lem:actionOnDifferentials}, the two isomorphisms
$g_\ast \mathcal O_C\cong\mathcal O_C$ and $g_\ast\Omega_{C/k}^1\cong \Omega_{C/k}^1$ induce
an isomorphism
\begin{align*}
g_\ast \T_C =g_\ast \SHom(\Omega_{C/k}^1,\mathcal O_C)
=\SHom(g_\ast\Omega_{C/k}^1,g_\ast\mathcal O_C)
\cong \SHom(\Omega_{C/k}^1,\mathcal O_C)=\T_C.
\end{align*}
Then the following diagram of sheaf isomorphisms commutes:
\begin{align*}
\begin{tikzcd}[ampersand replacement=\&]
\&g_\ast \mathcal O_C(-\omega)\ar[r,"\cong"]\&\mathcal O_C(-\omega)\ar[d,"\cdot(\zeta_1\cdots\zeta_n)^{-1}"]\\
g_\ast \T_C\ar[ru,"{g_\ast\SHom(\psi,\mathcal O_C)}"]\ar[r,"\cong"]\&\T_C\ar[r,"{\SHom(\psi,\mathcal O_C)}"]\&\mathcal O_C(-\omega)
\end{tikzcd}
\end{align*}
\end{corollary}
\begin{proof}
Apply $\SHom(\blank,\mathcal O_C)$ to Lemma~\ref{lem:actionOnDifferentials} and use $\SHom(\blank,\mathcal O_C)\xrightarrow\cong \SHom(\blank,g_\ast\mathcal O_C)$.
\end{proof}
\subsection{Calculating Hochschild (co)homology}
The stage is set to apply Propositions~\ref{prop:cohomologyCurveQuotient} and \ref{prop:homologyCurveQuotient}.
For this calculation, we assume that the weight vector $\mathbf a = (a_1,\dotsc,a_n)$ satisfies the weight condition.

Denote by $X=\mathbb P^1\langle x_1,\dotsc,x_n; a_1,\dotsc,a_n\rangle$ the weighted projective line.
\begin{lemma}
$\H^0(C,\Omega^1_{C/k})$ decomposes into a direct sum of one-dimensional $G$-invariant subspaces,
and the action of $G$ on each of those is non-trivial.
\end{lemma}
\begin{proof}
First let $g\in G$ corresponding to $(\zeta_1,\dotsc,\zeta_n)\in \mu_{a_1}\times\cdots\times\mu_{a_n}$. We would like to study how it acts on $\H^0(C,\Omega^1_{C/k})$.
We apply $\H^0(C,\blank)$ to Lemma~\ref{lem:actionOnDifferentials}.
Note that $\H^0(C,\blank) = \H^0(C,g_\ast (\blank))$ by definition of $g_\ast$,
so that the action of $g$ on $\H^0(C,\Omega^1_{C/k})$ is given
by applying $\H^0(C,\blank)$ to $\Omega^1_{C/k}\cong g_\ast \Omega^1_{C/k}$.
Denote this induced action also by $g$.
We similarly get an induced action $\H^0(C,\mathcal O_C(\omega))\rightarrow\H^0(C,\mathcal O_C(\omega))$
such that, by Lemma~\ref{lem:actionOnDifferentials}, the following diagram of isomorphisms commutes:
\begin{align*}
\begin{tikzcd}[ampersand replacement=\&]
\H^0(C,\mathcal O_C(\omega))\ar[r,"{\H^0(C,\psi)}"]\ar[d,"\cdot(\zeta_1\cdots\zeta_n)^{-1}"]\&\H^0(C,\Omega^1_{C/k})\ar[r,"g"]\&\H^0(C,\Omega^1_{C/k})\\
\H^0(C,\mathcal O_C(\omega))\ar[r]\&\H^0(C,\mathcal O_C(\omega))\ar[ru,"{\H^0(C,\psi)}"']
\end{tikzcd}.
\end{align*}
Now the map $\H^0(C,\mathcal O_C(\omega))\rightarrow \H^0(C,\mathcal O_C(\omega))$ allows a simpler description:
$\H^0(C,\mathcal O_C(\omega))\cong S_\omega$, and the induced action from $(\zeta_1,\dotsc,\zeta_n)$
comes from the action of $(\zeta_1^{-1},\dotsc,\zeta_n^{-1})$ on $S$.

We get the decomposition
\begin{align*}
\H^0(C,\mathcal O_C(\omega))\cong \bigoplus_{\substack{e_1,\dotsc,e_n\in \mathbb Z_{\geq 0}\\
e_1\deg X_1+\cdots+e_n\deg X_n=\omega\\
e_i<a_i\text{ for }i\geq 3}} k X_1^{e_1}\cdots X_n^{e_n}
\end{align*}
where $(\zeta_1,\dotsc,\zeta_n)$ acts on
$k X_1^{e_1}\cdots X_n^{e_n}$ by multiplication by $\zeta_1^{-e_1}\cdots \zeta_n^{-e_n}$.

Hence the isomorphism $\H^0(C,\psi)$ allows to decompose
\begin{align*}
\H^0(C,\Omega^1_{C/k})=\bigoplus_{\substack{e_1,\dotsc,e_n\in \mathbb Z_{\geq 0}\\
e_1\deg X_1+\cdots+e_n\deg X_n=\omega\\
e_i<a_i\text{ for }i\geq 3}} \underbrace{\H^0(C,\psi)(kX_1^{e_1}\cdots X_n^{e_n})}_{=: V_{e_1,\dotsc,e_n}},
\end{align*}
where $(\zeta_1,\dotsc,\zeta_n)$
acts on $V_{e_1,\dotsc,e_n}$ by multiplication by $\zeta_1^{-e_1-1}\cdots \zeta_n^{-e_n-1}$.

Suppose now that on some $V_{e_1,\dotsc,e_n}$, this action is trivial.
Then $\zeta_1^{-e_1-1}\cdots \zeta_n^{-e_n-1}=1$ for choices of $\zeta_1,\dotsc,\zeta_n$.
Choosing $\zeta_1$ to be a primitive $a_i$-th root of unity and all other $\zeta_i$ to be $1$,
we see $a_1 \mid -e_1-1$, so $a_1\mid e_1+1$. Similarly, $a_i\mid e_i+1$.

On the other hand, the condition $e_1 \deg X_1+\cdots + e_n\deg X_n = \omega$
is equivalent (after division by $\overline a$) to
\begin{align*}
\frac{e_1}{a_1}+\cdots + \frac{e_n}{a_n} = n-2-\sum_{i=1}^n \frac 1{a_i}
\iff \frac{e_1+1}{a_1}+\cdots + \frac{e_n+1}{a_n} = n-2.
\end{align*}
Since the $e_i$ are non-negative, $\frac{e_i+1}{a_i}>0$. As $\frac{e_i+1}{a_i}$ is an integer,
it is $\geq 1$, so we get the contradiction $n-2\geq \sum\limits_{i=1}^n 1 = n$.

This contradiction shows that the action on $V_{e_1,\dotsc,e_n}$ cannot be trivial, finishing the proof.
\end{proof}
\begin{corollary}\label{cor:coinvariantsGlobalSectionsCotangentSheaf}
$\H^0(C,\Omega^1_{C/k})^G\cong \H^0(C,\Omega^1_{C/k})_G\cong 0$.\rightqed
\end{corollary}
\begin{proposition}
The Hochschild homology of $X \cong [C/G]$ is given by
\begin{align*}
\dim\HH_0(X) =~&2-n+\sum_{i=1}^n a_i
\\\HH_i(X)=~&0,\text{ for }i\in \mathbb Z\setminus\{0\}.
\end{align*}
\end{proposition}
\begin{proof}
We use Proposition~\ref{prop:homologyCurveQuotient}. Then $\HH_0([C/G])$ has dimension
\begin{align*}
2+\sum_{Gc\in G\setminus C} (\#G_c-1)
\underset{\text{L\ref{lem:groupActionOnPoints}}}=2+\sum_{i=1}^n (\#G_{\mathfrak p_i}-1) = 2-n+\sum_{i=1}^n \# G_{\mathfrak p_i}
\underset{\text{L\ref{lem:groupActionOnPoints}}}=\sum_{i=1}^n 2-n+\sum_{i=1}^n a_i.
\end{align*}
for some choice of points $\mathfrak p_i\in V(X_i)$ for $i=1,\dotsc,n$.

Moreover $\HH_{\pm 1}([C/G])=0$ follows using Corollary~\ref{cor:coinvariantsGlobalSectionsCotangentSheaf}.
$\HH_i([C/G])=0$ for $i\notin\{0,\pm 1\}$ is directly from Proposition~\ref{prop:homologyCurveQuotient}.
This finishes the calculation.
\end{proof}
\begin{proposition}
The Hochschild cohomology of $X\cong [C/G]$ is given by
\begin{align*}
&\dim \HH^0(X)=1,
\\&\dim \HH^1(X)=\begin{cases}1,&n=2\\0,&n\geq 3\end{cases}
\\&\dim\HH^2(X)=\begin{cases}0,&n\leq 3\\n-3,&n\geq 3\end{cases}
\end{align*}
(and all other cohomology groups vanish).
\end{proposition}
\begin{proof}
We apply Proposition~\ref{prop:cohomologyCurveQuotient}. Hence
\begin{align*}
&\HH^0(X)\cong k,
\\&\HH^1(X)\cong \H^0(C,\T_C)_G\oplus \H^1(C,\mathcal O_C)_G,
\\&\HH^2(X)\cong \H^1(C,\T_C)_G.
\end{align*}
By Serre duality and Corollary~\ref{cor:coinvariantsGlobalSectionsCotangentSheaf},
$\H^1(C,\mathcal O_C)_G\cong \H^0(C,\Omega^1_{C/g})_G\cong 0$.

Now we calculate the $G$-action on $\H^\bullet(C,\T_C)$.
For this, we use the \v{C}ech complex associated to the affine open covering in Lemma~\ref{lem:standardAffine}.
\begin{align*}
0\rightarrow \H^0(U,\T_C)\oplus \H^0(V,\T_C)\rightarrow \H^0(U\cap V,\T_C)\rightarrow 0.
\end{align*}
By Corollary~\ref{cor:actionOnCodifferentials}, 
\begin{align*}
\H^0(U,\T_C)\cong \H^0(U,\mathcal O_C(-\omega))\cong\bigoplus_{\substack{e_1,\dotsc,e_n\in\mathbb Z\\
e_1,e_3,\dotsc,e_n\geq 0\\
e_i<a_i\text{ for }i\geq 3\\
e_1\deg X_1+\cdots+e_n\deg X_n = -\omega}} k X_1^{e_1}\cdots X_n^{e_n}.
\end{align*}
The group element $g\in G$ associated to $(\zeta_1,\dotsc,\zeta_n)\in \mu_1\times\cdots\times\mu_n$
acts on $X_1^{e_1}\cdots X_n^{e_n}$ by multiplication by $\zeta_1^{-e_1}\cdots \zeta_n^{-e_n}$.

We conclude by Corollary~\ref{cor:actionOnCodifferentials} that
\begin{align*}
\H^0(U,\T_C)\cong \bigoplus_{\cdots} \underbrace{k\H^0(U,\SHom(\psi,\mathcal O_C))(X_1^{e_1}\cdots X_n^{e_n})}_{=: W_{U,e_1,\dotsc,e_n}}
\end{align*}
with the same condition on the $e_i$, where $g$ acts on $W_{U,e_1,\dotsc,e_n}$ by multiplication
by $\zeta_1^{1-e_1}\cdots \zeta_n^{1-e_n}$.

We similarly get
\begin{align*}
&\H^0(V,\T_C)\cong \bigoplus_{\substack{e_1,\dotsc,e_n\in\mathbb Z\\
e_2,\dotsc,e_n\geq 0\\
e_i< a_i\text{ for }i\geq 3\\
e_1\deg X_1+\cdots+e_n\deg X_n=\omega}} \underbrace{k\H^0(V,\SHom(\psi,\mathcal O_C))(X_1^{e_1}\cdots X_n^{e_n})}_{=: W_{V,e_1,\dotsc,e_n}}
\\&\H^0(U\cap V,\T_C)\cong \bigoplus_{\substack{e_1,\dotsc,e_n\in\mathbb Z\\
e_3,\dotsc,e_n\geq 0\\
e_i< a_i\text{ for }i\geq 3\\
e_1\deg X_1+\cdots+e_n\deg X_n=\omega}} \underbrace{k\H^0(U\cap V,\SHom(\psi,\mathcal O_C))(X_1^{e_1}\cdots X_n^{e_n})}_{=: W_{U\cap V,e_1,\dotsc,e_n}}.
\end{align*}
Here, $g$ acts on $W_{V,e_1,\dotsc,e_n}$ and on $W_{U\cap V,e_1,\dotsc,e_n}$ by multiplication
by $\zeta_1^{1-e_1}\cdots \zeta_n^{1-e_n}$.

We conclude, using the above \v{C}ech complex, that
\begin{align*}
&\H^0(C,\T_C)\cong \bigoplus_{\substack{e_1,\dotsc,e_n\in\mathbb Z_{\geq 0}\\
e_i<a_i\text{ for }i\geq 3\\
e_1\deg X_1+\cdots+e_n\deg X_n=-\omega}} W_{e_1,\dotsc,e_n},
\\&\H^1(C,\T_C)\cong
\bigoplus_{\substack{e_1,e_2\in\mathbb Z_{<0};~e_3\dotsc,e_n\in\mathbb Z_{\geq 0}\\
e_i<a_i\text{ for }i\geq 3\\
e_1\deg X_1+\cdots+e_n\deg X_n=-\omega}} W_{e_1,\dotsc,e_n},
\end{align*}
where $W_{e_1,\dotsc,e_n}$ is a one-dimensional $k$-vector space
on which $g$ acts by multiplication by $\zeta_1^{1-e_1}\cdots\zeta_n^{1-e_n}$.

Thus $W_{e_1,\dotsc,e_n}$ is a $G$-invariant subspace, and it is trivial as $G$-representation
iff $a_i$ divides $1-e_i$ for $i=1,\dotsc,n$.

Hence
$\dim \H^0(C,\T_C)_G$ equals the number of tuples $(e_1,\dotsc,e_n)$ such that
\begin{itemize}
\item All $e_i$ are in $\mathbb Z_{\geq 0}$.
\item $e_i<a_i$ if $i\geq 3$.
\item $e_i\equiv 1\pmod{a_i}$ for all $i$.
\item $e_1 \deg X_1+\cdots + e_n\deg X_n=-\omega$.
\end{itemize}
The last condition is, after division by $\overline a$, equivalent
to
\begin{align*}
\frac{e_1}{a_1}+\cdots + \frac{e_n}{a_n} = 2-n+\sum_{i=1}^n \frac 1{a_i}
\iff n-2=\frac{1-e_1}{a_1}+\cdots + \frac{1-e_n}{a_n}.\tag{$\ast$}
\end{align*}
The condition $e_i\equiv 1\pmod{a_i}$ shows that $\frac{1-e_i}{a_i}$ is an integer for all $i$.
Since $e_i\geq 0$, $\frac{1-e_i}{a_i}\leq 0$ with equality iff $e_i=1$. 
We see that $(e_1,\dotsc,e_n)=(0,\dotsc,0)$ is the only solution to $(\ast)$ if $n=2$,
and that $(\ast)$ has no solutions for $n\geq 3$.

Therefore, $\HH^1([C/G])\cong \H^0(C,\T_C)_G$ is one-dimensional for $n=2$ and zero for $n\geq 3$.

Similarly, $\dim \H^1(C,\T_C)_G$ equals the number of tuples $(e_1,\dotsc,e_n)$ where
\begin{itemize}
\item All $e_i$ are in $\mathbb Z$, $e_1,e_2<0$ and $e_3,\dotsc,e_n\geq 0$.
\item $e_i<a_i$ for $i\geq 3$.
\item $e_i\equiv 1\pmod{a_i}$ for all $i$.
\item $e_1\deg X_1 + \cdots + e_n\deg X_n=-\omega$.
\end{itemize}
As before, the last condition is equivalent to $(\ast)$.

Moreover, for $i\geq 3$, we have $0\leq e_i<a_i$ and $e_i\equiv 1\pmod a_i$.
Thus $e_i=1$. Then $(\ast)$ simplifies to $n-2=\frac{1-e_1}{a_1}+\frac{1-e_2}{a_2}$.

Writing $\frac{1-e_1}{a_1}=\ell_1,\frac{1-e_2}{a_2}=\ell_2\in \mathbb Z_{>0}$,
we see that $\dim \H^1(C,\T_C)_G$ equals the number of tuples $(\ell_1,\ell_2)\in\mathbb Z_{> 0}^2$
such that $\ell_1+\ell_2=n-2$. By elementary considerations, the number of solutions is
$\max\{0,n-3\}$, finishing the calculation of $\HH^2([C/G])$.
\end{proof}
We reproduced Example~\ref{example:calculationWPLFromAbstractResult} for those weighted projective lines where the weight vector satisfies the weight condition.
\section{Weighted projective lines of non-negative Euler characteristic}\label{sec:nonNegEuler}
Recall that the Euler characteristic of the weighted projective line $X = \mathbb P^1\langle x_1,\dotsc,x_n; a_1,\dotsc,a_n\rangle$
is (by \cite[58]{Voight2019})
\begin{align*}
\chi_X = \sum_{i=1}^n 1-\frac 1{a_i}.
\end{align*}
If $X$ is a quotient stack $X\cong [C/G]$, we have
$\chi_X = \frac{\chi_C}{\# G}$ (by \cite[Theorem~2]{Lenzing2016}), where $\chi_C = 2-2g_C$ is the Euler characteristic
of $C$.

By \cite[Remark~5.4]{Geigle1987}, the complexity of $\coh(X)$ depends on whether
$\chi_X<0, \chi_X=0$ or $\chi_X>0$. Note that if $X\cong [C/G]$, $\chi_X$ is negative/zero/positive
iff $\chi_C$ is negative/zero/positive. In terms of \cite[Theorem~7]{Lenzing2016},
we call the weighted projective line $X$ \emph{spherical/parabolic/hyperbolic} if
$\chi_X>0$ resp.\ $\chi_X=0$ resp.\ $\chi_X<0$.

In this section, we classify how the spherical and parabolic weighted projective lines can be realized
as quotient stacks. We won't give explicit calculations of the Hochschild (co)homology as in section~\ref{sec:projectiveCompanion},
but refer to the abstract calculations in the proof of Theorem~\ref{thm:summaryQuotient}.

The class of hyperbolic weighted projective lines seems to be too large to admit such a classification.
\subsection{Weighted projective lines of positive Euler characteristic}
Suppose that $\chi_X$ is positive. Similar to Example~\ref{ex:genusZero}, we can show that this holds iff the weight vector $\mathbf a=(a_1,\dotsc,a_n)$ satisfies
\begin{align*}
\mathbf a \in \{(a,b)\mid a,b\in \mathbb Z_{\geq 1}\}\cup \{(2,2,a)\mid a\geq 2\}\cup\{(2,3,3), (2,3,4), (2,3,5)\}.
\end{align*}
Suppose that $X\cong [C/G]$ for a smooth and irreducible projective curve $C$ and a finite subgroup of automorphisms
$G\leq \Aut(C)$. Then $\chi_C>0$. Thus $C$ has genus zero, so $C\cong \mathbb P^1$.
\begin{proposition}[Consequence of {\cite[Theorem~13]{Lenzing2016}}]
Each finite subgroup $G\leq \Aut(\mathbb P^1)$ is isomorphic to one of the following:
\begin{itemize}
\item a cyclic group $C_n$ for $n\geq 1$,
\item a dihedral group $D_{2n}$ of order $2n$ for $n\geq 1$,
\item the alternating group $A_4$ or $A_5$ or the symmetric group $S_4$.
\end{itemize}
Moreover, these give the following weighted projective lines:
\begin{align*}
&[\mathbb P^1/C_n]\cong \mathbb P^1\langle x_1,x_2; n,n\rangle
\\&[\mathbb P^1/D_{2n}]\cong \mathbb P^1\langle x_1,x_2,x_3; 2,2,n\rangle
\\&[\mathbb P^1/A_4]\cong \mathbb P^1\langle x_1,x_2,x_3; 2,3,3\rangle
\\&[\mathbb P^1/S_4]\cong \mathbb P^1\langle x_1,x_2,x_3; 2,3,4\rangle
\\&[\mathbb P^1/A_5]\cong \mathbb P^1\langle x_1,x_2,x_3; 2,3,5\rangle
\rightqed
\end{align*}
\end{proposition}
\begin{corollary}
If $n=2$ and $a_1\neq a_2$, the weighted projective line $X=\mathbb P^1\langle x_1,x_2; a_1,a_2\rangle$
is not isomorphic to a quotient stack $[C/G]$ for any smooth projective curve $C$ and finite group $G\leq \Aut(C)$.\rightqed
\end{corollary}
This gives a complete description when and how a weighted projective line of positive Euler characteristic can be realized as a quotient stack.
\subsection{Weighted projective lines of Euler characteristic zero}
Suppose now that $\chi_X=0$. By Example~\ref{ex:genusOne}, this holds iff the weight vector $\mathbf a=(a_1,\dotsc,a_n)$ is one of the following:
\begin{align*}
\mathbf a = (2,3,6),~(2,4,4),~(3,3,3),~(2,2,2,2).
\end{align*}
If $X\cong [C/G]$ for a smooth irreducible projective curve $C$ and a finite group $G\leq \Aut(C)$,
we must have $\chi_C=0$. Hence $C$ will be an elliptic curve.

Recall that for each elliptic curve $C$, the set of points of $C$ comes, after fixing a basepoint $P_0$, with a natural group structure (\cite[Section~IV.4]{Hartshorne1977}).
For a point $P\in C$, let $\tau_P\in \Aut(C)$ be the automorphism sending $Q$ to $P+Q$.
Hence each automorphism of $C$ has the form $\tau_P \sigma$ for a point $P\in C$ and an automorphism $\sigma\in\Aut(C)$ with $\sigma(P_0)=P_0$.
\begin{lemma}\label{lem:fixedPointAutomorphismsEllipticCurves}
Let $G\leq \Aut(C)$ be the group of those automorphisms $\sigma$ with $\sigma(P_0)=P_0$.
Then $G$ is cyclic, and the order of $G$ depends on the $j$-invariant of $C$ as follows:
\begin{align*}
\# G = \begin{cases}6,&j(C)=0\\4,&j(C) = 1728\\2,~&\text{otherwise}.\end{cases}
\end{align*}
\begin{proof}
Finiteness of $G$ is shown in \cite[Corollary~IV.4.7]{Hartshorne1977}, where also the order is stated explicitly.
The fact that $G$ is cyclic follows from Proposition~\ref{prop:curveFiniteGroup}.
\end{proof}
\end{lemma}
\begin{example}\label{example:automorphismsOfEllipticCurves}
The $j$-invariant determines an elliptic curve uniquely up to isomorphism \cite[Theorem~IV.4.1]{Hartshorne1977}.
Hence there exists, up to isomorphism, exactly one elliptic curve with a basepoint preserving automorphism
of order $4$ (resp.\ $6$).

By Example~\ref{ex:genusOne}, the weighted projective lines \begin{align*}\mathbb P^1\langle x_1,x_2,x_3; 2,4,4\rangle
\qquad\text{and}\qquad \mathbb P^1\langle x_1,x_2,x_3; 2,3,6\rangle\end{align*} have projective companions
which are elliptic curves with a basepoint preserving automorphism
of order $4$ resp.\ $6$.

If $C$ is any elliptic curve, a basepoint preserving automorphism of order $2$ is given by the inversion map coming from the group structure of $C$.
Since each cyclic group contains at most one involution, this is moreover the only basepoint preserving automorphism of order $2$.
\end{example}
\begin{proposition}\label{prop:quotientsOfEllipticCurves}
Let $C$ be an irreducible elliptic curve and $G\leq\Aut(C)$ be a finite group of automorphisms.
Pick a basepoint $P_0\in C$ such that the stabilizer $G_{P_0}$ has maximal order, and let $\tau(C) := \{\tau_P \mid P\in C\}\leq\Aut(C)$
be the group of automorphisms which come from the group structure of $C$. Define $N := G\cap \tau(C)$. Then one of the following holds:
\begin{itemize}
\item $G\subseteq \tau(C)$ and $[C/G]$ is an elliptic curve (without stacky points).
\item $G$ is a semi-direct product $N\rtimes \langle \sigma\rangle$ where $\sigma$ is the involution from Example~\ref{example:automorphismsOfEllipticCurves},
there are precisely four orbits of points of $C$ which have non-trivial stabilizer and the quotient is $[C/G]\cong \mathbb P^1\langle x_1,x_2,x_3,x_4;~2,2,2,2\rangle$.
\item $C$ has $j$-invariant $0$, $G$ is a semi-direct product $N\rtimes \langle \varphi^2\rangle$ where $\varphi$ is the basepoint preserving
automorphism of order $6$ from Example~\ref{example:automorphismsOfEllipticCurves} and $[C/G]\cong \mathbb P^1\langle x_1,x_2,x_3;~3,3,3\rangle$.
\item $C$ has $j$-invariant $0$, $G$ is a semi-direct product $N\rtimes \langle \varphi\rangle$ where $\varphi$ is the basepoint preserving
automorphism of order $6$ from Example~\ref{example:automorphismsOfEllipticCurves} and $[C/G]\cong \mathbb P^1\langle x_1,x_2,x_3;~2,3,6\rangle$.
\item $C$ has $j$-invariant $1728$, $G$ is a semi-direct product $N\rtimes \langle \psi\rangle$ where $\psi$
is the basepoint preserving automorphism of order $4$ from Example \ref{example:automorphismsOfEllipticCurves}
and $[C/G]\cong\mathbb P^1\langle x_1,x_2,x_3;~2,4,4\rangle$.
\end{itemize}
\end{proposition}
This proposition was what we were searching for: A complete description how weighted projective lines of Euler characteristic zero
can be represented as quotient stacks of smooth projective curves by a finite group.

For the proof of the proposition, we need a Lemma.
\begin{lemma}\label{lem:factsEllipticCurves}
Let $C$ be an irreducible elliptic curve, $P_0\in C$ a basepoint defining a group structure
and $\tau(C) := \{\tau_P\mid P\in C\}\leq \Aut(C)$.
\begin{enumerate}[(a)]
\item If $P\in C$ is a closed point and $\phi\in \Aut(C)$ an automorphism with $\phi(P_0)=P_0$,
we have $\phi\tau_P\phi^{-1} = \tau_{\phi(P)}$.

In particular, $\tau(C)$ is normal in $\Aut(C)$.
\item If $\phi\neq\id_C$ is an automorphism of $C$, the following are equivalent:
\begin{itemize}
\item $\phi$ has no closed fixed points in $C$.
\item There exists a closed point $P_0\neq P\in C$ with $\phi=\tau_P$.
\end{itemize}
\end{enumerate}
\end{lemma}
\begin{proof}
\begin{enumerate}[(a)]
\item The fact $\phi\tau_P\phi^{-1} = \tau_{\phi(P)}$ follows from \cite[Lemma~IV.4.9]{Hartshorne1977}.

Now let $\tilde\phi\in \Aut(C)$ be any automorphism, write $P := \tilde\phi(P_0)$ and $\phi := \tau_{-P}\tilde\phi$.
Then $\phi(P_0)=P_0$ and $\tilde\phi = \tau_P\phi$. By the previous argument, $\phi$ is in the normalizer of $\tau(C)$.
As $\tau_P$ is in $\tau(C)$, hence in the normalizer of $\tau(C)$, we conclude that $\tilde\phi$ is in the normalizer of $\tau(C)$.
\item First note that if $\tau_P$ has a closed fixed point $Q$, we have $P+Q=Q$ so that $P=P_0$.

Now suppose that $\phi$ is an automorphism without closed fixed points
and consider the endomorphism $f:C\rightarrow C, Q \mapsto \varphi(Q)-Q$. 
This is a morphism of schemes by \cite[Proposition~IV.4.8]{Hartshorne1977}.

By \cite[Proposition~II.6.8]{Hartshorne1977}, $f$ is constant or surjective. By assumption,
$P_0$ is not in the image of $f$, hence $f$ is constant.
If $\{P\} = \im(f)$, we conclude $\varphi=\tau_P$.
\qedhere\end{enumerate}
\end{proof}
\begin{proof}[Proof of Proposition~\ref{prop:quotientsOfEllipticCurves}]
Note that by Lemma~\ref{lem:factsEllipticCurves},
$N$ is normal in $G$. Moreover, the induced injection $G/N\hookrightarrow \Aut(C)/\tau(C)$
shows that $G/N$ is cyclic of order $1,2,4$ or $6$ by Lemma~\ref{lem:fixedPointAutomorphismsEllipticCurves}.

Instead of classifying finite subgroups of $\Aut(C)$, we examine the possibilities for $[C/G]$.

If $G$ acts freely on the points of $C$ (i.e.\ all stabilizers are trivial),
$[C/G]=C/G$ is a scheme (by Proposition~\ref{prop:quotientRootStack}), and $\chi_{C/G}=\frac{\chi_C}{\# G} = 0$
implies that $C/G$ is an elliptic curve without stacky points.

Now suppose that $G$ does not act freely on $C$. By Riemann-Hurwitz, $\chi_{C/G} > \frac{\chi_C}{\# G}=0$, so $C/G$ must have genus $0$,
i.e.\ $C/G\cong \mathbb P^1$. We conclude that $[C/G]$ is a weighted projective line of Euler characteristic zero.
Since $G$ does not act freely on $C$ but $N$ does, we conclude $G\neq N$.
Hence $G/N$ is cyclic of order $2,4$ or $6$.

Write $[C/G]\cong \mathbb P^1\langle x_1,\dotsc,x_n;~a_1,\dotsc,a_n\rangle$. We consider the possibilities for $\mathbf a = (a_1,\dotsc,a_n)$:

\textbf{Case $\mathbf a = (2,2,2,2)$.} This means that there are precisely four $G$-orbits of closed points in $C$ with
stabilizer of order $2$, and all other closed points of $C$ have trivial stabilizer. By choice of $P_0$,
the stabilizer of $P_0$ is non-zero. By Example~\ref{example:automorphismsOfEllipticCurves}, the stabilizer
of $P_0$ is $\{\id,\sigma\}$ for the involution $\sigma$. Thus $\sigma\in G$.

It remains to show that $G$ is generated by $\sigma$ and $N$. Suppose this was not the case.
Then $G/N$ would be cyclic of order $>2$. Hence there would be an automorphism
$\phi\in G$ whose image in $G/N$ has order $>2$. In particular, the order of $\phi$ in $G$ is $>2$.
Since $\phi\notin N$, Lemma~\ref{lem:factsEllipticCurves} shows that $\phi$ has a fixed point. This contradicts
the fact that the point stabilizers all have order $\leq 2$.

\textbf{Case $\mathbf a = (3,3,3)$.} This means that there are precisely three $G$-orbits of closed points in $C$
with stabilizer of order $3$, and all other closed points of $C$ have trivial stabilizer. 

In particular, $C$ has an automorphism of order $3$ which fixed $P_0$ (by choice of $P_0$). By Lemma~\ref{lem:fixedPointAutomorphismsEllipticCurves},
the $j$-invariant of $C$ equals $0$ and this automorphism must be $\varphi^2$.

It remains to show that $G$ is generated by $\varphi^2$ and $N$. If this was not the case,
we would again find an automorphism $\phi\in G$ with $\phi\notin N$ and the order of $\phi$ in $G$ being $>3$.
Then $\phi$ has a closed fixed point by Lemma~\ref{lem:factsEllipticCurves}. This contradicts the fact that all stabilizers of
closed points have order $\leq 3$.

\textbf{Case $\mathbf a=(2,3,6)$.} Then the stabilizer of $P_0$ has order $6$ by choice of $P_0$.
Hence $j(C)=0$ and $\varphi\in G$. The image of $\varphi$ must generate $G/N$ (comparing orders),
so that $G$ is generated by $N$ and $\varphi$.

\textbf{Case $\mathbf a=(2,4,4)$.} Then the stabilizer of $P_0$ has order $4$ by choice of $P_0$.
Hence $j(C)=1728$ and $\varphi\in G$. The image of $\varphi$ must generate $G/N$ (comparing orders),
so that $G$ is generated by $N$ and $\varphi$.
\end{proof}
\newpage
\chapter*{Conclusion and Outlook}
\addcontentsline{toc}{chapter}{Conclusion and Outlook}

In this thesis, we presented two methods to compute the Hochschild homology and cohomology of weighted projective lines.

The first method, due to Happel, makes use of a tilting sheaf, which allows us to equivalently study the Hochschild (co)homology
of the associated canonical algebra. Next to Happel's original approach using one-point extensions, we presented an alternative
approach that describes the Hochschild cohomology of finite-dimensional algebras of global dimension $\leq 2$.

The second method combines results from \cite{Arinkin2014} and \cite{Lenzing2016}.
We showed how the Hochschild cohomology of a weighted projective curve can be computed
in case it can be realized as the orbifold quotient of a smooth projective curve by a finite group of automorphisms.
By \cite{Lenzing2016}, all weighted projective curves except for a small class of weighted projective lines can be realized as such a quotient stack.

Both methods combined yield a description of the Hochschild (co)homology of all weighted projective curves (Theorem~\ref{thm:summaryRootStack}).
Note that an iterated root stack obtained from a curve of genus $\geq 1$ will have non-vanishing Hochschild homology in degree $-1$,
so it cannot be derived equivalent to a finite-dimensional algebra.

There remain, of course, unanswered questions, and possibilities for further research.

In the entire thesis, we focussed on the calculation of the Hochschild (co)homology groups as vector spaces.
We did not describe higher structures, like e.g.\ the Gerstenhaber algebra structure on $\HH^\bullet$.

Moroever, we did not discuss interpretations and applications of Hochschild (co)homology groups, like their usage in deformation theory.

Finally, it is worth noting that the result in \cite{Arinkin2014}, and in some cases also Happel's method, can be applied
to (weighted) surfaces and stacks of higher dimension. However, Lenzing's methods in \cite{Lenzing2016} seem to be limited to 
curves.
\chapter*{Acknowledgements}
\pdfbookmark[chapter]{Acknowledgements}{acknowledgements}
I would like to thank my advisor, Dr.\ Pieter Belmans, for the many insightful and inspiring discussions.
\newpage
\addcontentsline{toc}{chapter}{Bibliography}
\printbibliography
\appendix
\chapter{GAP routines}\label{chap:gap_routines}
We present two routines for the computer algebra system \cite{GAP}
which can be used to calculate the Hochschild cohomology
of canonical algebras.
\begin{lstlisting}
LoadPackage("qpa");
show_hh := function(A)
	#Shows the Hochschild cohomology in degrees >= 1.
	#Terminates when it reaches the global dimension of A.
	local M, N, i, ext, d;
	M:=AlgebraAsModuleOverEnvelopingAlgebra(A);
	#output the Ext^n(M,M) as follows:
	#Since gap has no Ext^n-function,
	#we have to use Ext^1 and iteratively
	#replace the first M by its sysygies.
	N:=M;
	i:=1;
	repeat
		d:=Length(ExtOverAlgebra(N,M)[2]);
		Print(i,": ",d,"\n");
		N:=1stSyzygy(N);
		i:=i+1;
	until IsProjectiveModule(N);
end;
canonical_algebra := function(weights, points)
	#Takes the weights and points to construct a canonical algebra.
	#The first two entries in the array points are ignored; they
	#are assumed to be infinity and zero.
	#Moreover, we assume the length of weights equals the
	#length of points, which are both assumed to be at least 2.
	#Finally, we assume that at most one weight equals 1.
	local pt, step, vertices, numPoints, arrows, Q, vertexName,
	   arrowName, previousVertex, strings, relations, gens,
	   index, string, kQ, relation;
	numPoints := Length(weights);
	vertices := ["0","c"];
	arrows:=[];
	#builds an array of vertices and and array of arrows
	for pt in [1..numPoints] do
		previousVertex:="0";
		for step in [1..weights[pt]] do
			vertexName := StringFormatted("v_{}_{}", pt, step);
			if step = weights[pt] then
				vertexName := "c";
			else
				vertices := Concatenation(vertices, [vertexName]);
			fi;
			arrowName := StringFormatted("a_{}_{}", pt,step);
			arrows := Concatenation(arrows, [[previousVertex, vertexName, arrowName]]);
			previousVertex := vertexName;
		od;
	od;
	Q := Quiver(vertices, arrows);
	kQ := PathAlgebra(Rationals, Q);
	strings:=[];
	relations := [];
	gens := GeneratorsOfAlgebra(kQ);
	#gens first contains the generators corresponding to vertices
	#and then the generators corresponding to arrows
	#Skip over the vertices:
	index := 2;
	for pt in [1..numPoints] do
		index := index+weights[pt]-1;
	od;
	#iteratively build relations.
	for pt in [1..numPoints] do
		string := 1;
		for step in [1..weights[pt]] do
			index := index+1;
			string := string * gens[index];
		od;
		strings := Concatenation(strings,[string]);
		if pt > 2 then
			relation:=string - strings[1] + points[pt] * strings[2];
			relations := Concatenation(relations, [relation]);
		fi;
	od;
	return kQ/relations;
end;
show_hh(canonical_algebra([2,2,2,2],[0,0,1,2]));
\end{lstlisting}
This outputs
\lstset{language=html}
\begin{lstlisting}
1: 0
Dimension vector for syzygy: [ 0, 6, 1, [...], 0 ]
2: 1
Dimension vector for syzygy: [ 0, 2, 0, [...], 0 ]
\end{lstlisting}
This means that the weighted projective line $X$ associated to the points
\begin{align*}x_1=\infty,x_2=0,x_3=1,x_4=2\end{align*} and the weights \begin{align*}a_1=a_2=a_3=a_4=2\end{align*}
has the following Hochschild cohomology groups:
\begin{align*}
&\HH^1(X)\cong 0
\\&\HH^2(X)\cong k
\\&\HH^n(X)\cong 0\text{ for }n\geq 3.
\end{align*}
\end{document}